%% file: article_blm.tex
\documentclass[12pt]{elsart}

\usepackage{graphicx}
\usepackage{psfrag}
\usepackage[dvipsnames]{xcolor}
\usepackage{amssymb}
\usepackage{amsmath}
\usepackage{dsfont}
\usepackage{rotating}
\usepackage{bm}
\usepackage{url}
\usepackage[toc]{appendix}
\newcommand{\Frac}{\displaystyle\frac}
\newcommand{\Int}{\displaystyle\int}
\newcommand{\Sum}{\displaystyle\sum}

\newcommand{\bea}{\begin{equation}}
\newcommand{\eea}{\end{equation}}
\newcommand{\beqa}{\begin{eqnarray}}
\newcommand{\eeqa}{\end{eqnarray}}
\newcommand{\nn}{\nonumber}

\newcommand{\rmd}{{\rm d}}
\newcommand{\rmdt}{{\rm d} t}

\newcommand{\tr}{\text{tr}}
\newcommand{\mddt}[1]{\Frac{\rmd {#1}}{\rmdt}}
\newcommand{\pddt}[1]{\Frac{\partial #1}{\partial t}}
\newcommand{\pdd}[2]{\Frac{\partial #1}{\partial #2}}
\newcommand{\GRAD}[1]{\overrightarrow{\nabla #1}}
\newcommand{\DIV}[1]{\nabla \cdot #1}

\newcommand{\tens}[1]{\mathbb #1}
\renewcommand{\d}{\textnormal{d}}

\newcommand{\CPU}{\textnormal{CPU}}

\newcommand{\dt}{\Delta t}

\newcommand{\Q}{\mathbf{Q}} 
\newcommand{\q}{\mathbf{q}} 
\newcommand{\f}{\mathbf{f}}

\newcommand{\w}{\mathbf{w}}

\newcommand{\lcp}{\ell_{cp}}
\newcommand{\ncp}{\vec{n}_{cp}}
\newcommand{\lncp}{\lcp\ncp}
\newcommand{\lcpf}{\ell_{cpf}}

\newcommand{\Mcp}{\tens{M}_{cp}}
\newcommand{\Mcpf}{\tens{M}_{cpf}}
\newcommand{\Mp}{\tens{M}_{p}}

\newcommand{\aposteriori}{\textit{a posteriori }}
\newcommand{\apriori}{\textit{a priori }}
\newcommand{\BC}{\text{BC}}
\newfont{\numerikEleven}{ecrm1000}
\newfont{\numerikTen}{cmss10}
\newfont{\numerikNine}{cmss9}
\newfont{\numerikEight}{cmss8}
\newtheorem{remark}{Remark}

\begin{document}

\begin{frontmatter}
  
  %=================================================================
  \title{A cell-centered Lagrangian ADER-MOOD finite volume scheme on unstructured meshes for a class of hyper-elasticity models}
  %=================================================================
  
  %=================================================================
  \author[ferrara]{Walter Boscheri} 	\ead{walter.boscheri@unife.it},
  \author[imb]{Rapha{\"e}l Loub{\`e}re}	\ead{raphael.loubere@u-bordeaux.fr},
  \author[cea]{Pierre-Henri Maire}  	\ead{pierre-henri.maire@cea.fr}.
  %=================================================================
  \address[ferrara]{Dipartimento di Matematica e Informatica,	%Via Machiavelli 30,	44121 -
    Ferrara Italy}
  \address[imb]{Institut de Math\'{e}matiques de Bordeaux (IMB), %UMR5251, Universit{\'e} de Bordeaux,  F33400,
    Talence, France}
  \address[cea]{CEA-CESTA, %15 Avenue des Sablières, 33114
    Le Barp, France}
  %=================================================================
  
  %=================================================================
  \begin{abstract}
    % Text of abstract
    In this paper we present a conservative cell-centered Lagrangian finite volume 
    scheme for the solution of the hyper-elasticity equations on unstructured multidimensional grids.
    The starting point of the new method is the Eucclhyd scheme forwarded in \cite{Maire2007,Maire2010,Maire2011}, which is here combined with the \textit{a posteriori} Multidimensional Optimal Order Detection (MOOD)
    limiting strategy to ensure robustness and stability at shock waves with piece-wise linear spatial reconstruction. The ADER (Arbitrary high order schemes using DERivatives) approach \cite{mill,toro3} is adopted to obtain second-order of accuracy in time as well.
    This method has been tested in an hydrodynamics context in \cite{LAM2018} and the present work aims at extending it to the case of hyper-elasticity models.
    Such models are presented in a fully Lagrangian framework and the dedicated Lagrangian numerical
    scheme is derived in terms of nodal solver, GCL compliance, subcell forces and compatible discretization.
    The Lagrangian numerical method is implemented in 3D under MPI parallelization framework allowing to
    handle genuinely large meshes.
    A relative large set of numerical test cases is presented to assess the ability of the method to achieve effective
    second order of accuracy on smooth flows, maintaining an essentially non-oscillatory behavior and
    general robustness across discontinuities and ensuring at least physical admissibility of the solution where appropriate. 
    Pure elastic neo-Hookean and non-linear materials are considered for our benchmark test problems in 2D and 3D. These test cases feature material bending, impact, compression, non-linear deformation and further bouncing/detaching motions.
  \end{abstract}
  %=================================================================
  %=================================================================
\end{frontmatter}

\journal{Journal of Computational Physics}

%-------------------
% TABLE OF CONTENTS 
\tableofcontents
%------------------

% NOTES	
%\textcolor{red}{
%  NOTE: we can add comments in colors using \texttt{$\backslash$raph$\{\}$},
%  \texttt{$\backslash$walter$\{\}$},\texttt{$\backslash$ph$\{\}$}.
%  For instance \raphael{Raph's color}, \walter{Walter's} and \ph{PH's one}. 
%}

%---------------------------------------------%
%---------------------------------------------%
\section{Introduction}
\label{sec.introduction}
%---------------------------------------------%
%\raphael{@ALL: I have removed some citations when they appear in pack of three or four, but I have kept them
%  as comments in the .tex if you wish to put them back.\\}

% ---- Context
This work is concerned with the accurate multi-dimensional simulation of hyper-elasticity models by
updated Lagrangian Finite Volume (FV) numerical scheme.
% ---- Paper LAM2018
Previously we have presented a second-order accurate cell-centered Lagrangian scheme on unstructured mesh
for the hydrodynamics system of conservation laws restricted to 2D in \cite{LAM2018}.
This scheme is constructed upon a subcell discretization, popularized for staggered Lagrangian schemes \cite{Burton90,Caramana1998} and later extended to cell-centered ones \cite{Maire2010,Maire2011},
further associated with a nodal solver relying on total energy conservation and Geometrical Conservation Law (GCL) compliance.
Second-order of accuracy is usually gained by a MUSCL-like approach ---piece-wise linear reconstructions supplemented
with limiters--- and a predictor-corrector, Runge-Kutta or a Generalized Riemann Problem (GRP) time discretization. \\
Contrarily, for the scheme in \cite{LAM2018}, we have relied on ADER (Arbitrary high order schemes using DERivatives)
methodology developed originally from an Eulerian perspective \cite{toro10,Lagrange2D}. %,LagrangeQF,ALEMOOD1,ALEMQF}.
This is supplemented with an \aposteriori MOOD limiting strategy \cite{CDL1} %,CDL2,CDL3}
to stabilize and produce a fail-safe Lagrangian scheme.
We have shown in \cite{LAM2018} that such a cell-centered numerical method is able to perform on classical and demanding
hydrodynamics test cases using unstructured mesh made of simplexes in 1D and 2D. \\
% ---- Here hyper-elasticity
In this work we propose the extension of this numerical method in 3D to solve problems involving elastic materials.
We ought to solve an hyper-elasticity model of PDEs (Partial Differential Equations) \cite{Kluth10,LagrangeHPR,Gil2D_2014,CCL2020}.
Historically hypo-elasticity models \cite{Truesdell55,Bernstein60,Truesdell63} have been sometimes preferred by numericists, see for instance \cite{wil1,Gavriluk08,Maire_elasto,Sambasivan13,cheng_jia_jiang_toro_yu_2017}.
A parallel discussion about hypo- and hyper-elastic models and their resolution can be found for instance in \cite{Peshkov_Boscheri_Loub_hyper_hypo19}.
% --- Novelties
In this article are tackled several issues related to the 3D extension of our ADER Lagrangian scheme, as well as
the increase of complexity in the modeling of hyper-elastic materials.
First the hyper-elastic model demands the resolution of a constitutive law which, in the framework of ADER methodology,
requires some adaptation.
Second the \aposteriori MOOD limiting strategy must consider new admissibility criteria brought by the model
related to involution-like constrain of the materials in order to still ensure the robust and fail-safe characteristics
while maintaining an acceptable accuracy.
Third the boundary conditions (BCs) must be dealt with care as materials may balistically fly but also
impact, bounce, slide, spread, tear apart onto a wall or different materials. 
Fourth, in relation to the points above, 3D Lagrangian numerical simulation code requires extra-care as efficient 3D
simulations demand a well designed parallelization methodology as well as appropriate BCs and robust
limiting strategy. \\
% --- Numerical results
Numerical results involving materials enduring large deformation (bending, twisting, etc.) adopted from \cite{Gil2D_2014,Haider_2018,CCL2020}
will be presented to assess the ability of this updated Lagrangian numerical scheme to simulate such
hyper-elastic situations. \\
For a broad and modern introductions on hypo- or hyper-elasticity we refer the readers in particular to \cite{Kluth10,LagrangeHPR,CCL2020,Peshkov_Boscheri_Loub_hyper_hypo19}.
For 3D cell-centered Lagrangian computations among many works we refer to
\cite{SaltzmanOrg3D,LoubereSedov3D,MaireHD3D,CCL2020}.
% --- Paper organization
After this short introduction the paper presents in details the hyper-elastic model
and the governing equations to be solved.
Then in the third section, the Lagrangian numerical scheme is introduced with emphasis on the
ADER approach, the nodal solver and the \aposteriori limiting strategy.
All numerical tests are gathered in the fourth section.
We propose the numerical results of our simulations for the a large set of 2D and 3D problems
involving materials impacting, detaching, compressing, swinging, twisting, etc.
Conclusions and perspectives are finally drawn in the last section.
%---------------------------------------------%

\input{Hyperelasticity.tex}
%---------------------------------------------%
%---------------------------------------------%

	%---------------------------------------------%
	%---------------------------------------------%
\input{Discretization.tex}

%---------------------------------------------%
%\clearpage
%---------------------------------------------%

\input{Numerics.tex}

\clearpage
%---------------------------------------------%
%---------------------------------------------%
\section{Conclusions and perspectives}
\label{sec.conclusion}
%---------------------------------------------%
%
% Conclusions
%
This paper considers the second-order accurate cell-centered Lagrangian scheme
originally designed for the hydrodynamics system of conservation laws \cite{LAM2018},
and, extends it to solve the hyper-elasticity model for materials in 2D and 3D.
We have focused the first part of the paper on presenting the hyper-elasticity model and its consistency in the Lagrangian frame.
The so-called neo-Hookean model is mostly considered in this work.
Then the numerical method based on a conservative Lagrangian formulation in mass, momentum and total
energy is presented.
It is supplemented with a nodal solver allowing the determination of a vertex velocity used to  
build a consistent discretization between the trajectory equation and the geometrical conservation law.
Second-order of accuracy in space and time is achieved via an ADER procedure which generates a predictor solution that can further be used inside the classical subcell force based Lagrangian scheme with nodal solver. 
Robustness and stability are gained by the use of an \aposteriori MOOD limiting strategy, that is a second-order unlimited candidate solution at $t^{n+1}$ is tested against appropriate detection criteria to determine troubled cells. 
The solution in those cells is discarded and re-computed starting again from valid data at $t^n$ 
but using a second-order TVD like scheme or, ultimately, the fail-safe first-order Godunov parachute scheme. 
The constitutive equation on tensor $\tens{B}$ is solved in time using a second-order Crank-Nicholson scheme. Moreover evolving boundary conditions have been implemented to allow for impacting and detaching of materials onto walls. \\
This numerical scheme has been further implemented in 2D and 3D under MPI protocol for the parallelization.
It has been then tested on unstructured simplicial meshes on a large panel of 2D test cases:
 swinging plate,  elastic vibration of a beryllium plate and  a finite deformation of a cantilever thick beam.
Then, in 3D, we have presented the results for  Blake's problem,  the twisting column,  the rebound of a hollow circular bar and at last the  impact of a jelly-like droplet.
This test suite covers a large amount of situations involving elastic materials and the current Lagrangian numerical scheme has proven to be robust, essentially non-oscillatory and, at the same time maintains an almost optimal precision by a careful utilization of the high order scheme where appropriate and the low order ones in the vicinity of problematic zones.
Moreover its performance in 2D/3D both in terms of robustness, efficiency and compliance with other published results renders this numerical method appealing for future uses and possible coupling with more complex physical models. \\
% Perspectives
A plan for future work involves the introduction of plasticity into this hyper-elasticity model.
Another direction of evolution would be to add some Arbitrary-Lagrangian-Eulerian capability and the possibility to let two elastic materials interacting with each other, for instance impacting, deforming and further detaching from each others.

%=============================================================================
%==========    A C K N O W L E D G M E N T S
{
  %\scriptsize
  \section*{Acknowledgments}
  The material of this research has been partly built during the
  SHARK FV workshops which took place
  on May 2017, 2018, 2019 in Povoa de Varzim, Portugal
  \texttt{www.SHARK-FV.eu/}. \\  
}
%=============================================================================

\appendix
\section{Principal invariants of a tensor}
\label{sec:pit}
Let us consider an invertible second order tensor $\tens{A}$. Its principal invariants are the coefficients of its characteristic polynomial
$$\det(\tens{A}-\lambda \tens{I}_{\text{d}})=\lambda^{3} - I_1(\tens{A})\lambda^{2}+ I_2(\tens{A}) \lambda -I_3(\tens{A}).$$
They are defined by
\begin{subequations}
  \label{eq:inva}
  \begin{align}
    & I_1(\tens{A})=\tr (\tens{A}),\label{eq:inva1}\\
    & I_2(\tens{A})=\frac{1}{2}\left [ \tr^{2}(\tens{A})-\tr(\tens{A}^{2}) \right ],\label{eq:inva2}\\
    & I_3(\tens{A})=\det (\tens{A}). \label{eq:inva3}
  \end{align}
\end{subequations}
Moreover, the Cayley-Hamilton theorem states that $\tens{A}$ satisfies its characteristic equation
\begin{equation}
  \label{eq:ch}
  \tens{A}^{3} - I_1(\tens{A})\tens{A}^{2}+ I_2(\tens{A}) \tens{A}-I_3(\tens{A})\tens{I}_{\text{d}}=0.
  \end{equation}
The derivative of the principal invariants of $\tens{A}$ with respect to itself write
\begin{subequations}
  \label{eq:dinva}
  \begin{align}
    & \frac{\partial I_1(\tens{A})}{\partial \tens{A}}=\tens{I}_{\text{d}},\label{eq:dinva1}\\
    & \frac{\partial I_2(\tens{A})}{\partial \tens{A}}=I_1(\tens{A})\tens{I}_{\text{d}}-\tens{A},\label{eq:dinva2}\\
    & \frac{\partial I_3(\tens{A})}{\partial \tens{A}}=I_3(\tens{A})\tens{A}^{-1} . \label{eq:dinva3}
  \end{align}
\end{subequations}
\section{Boundary conditions (BCs)} \label{app:BCs}
In this cell-centered Lagrangian scheme, boundary conditions are imposed  in the nodal solver \eqref{eq:nodal_solver}.
Let $\mathcal{F}^{BC}(p)$ represent the set of boundary edges (or faces in 3D) surrounding node $p$.
Three type of BCs are considered in this work.
\begin{itemize}	
\item \textit{Traction}: a prescribed traction $\tens{T}_f^{BC}$ on a boundary face $f$ is taken into account as an additional term on the right hand side as
  \begin{equation}
    \vec{v}_{p} =   \Mp^{-1} \left( \Sum_{c\in  \mathcal{C}(p)}\Mcp \vec{v}_c - \Sum_{f\in \mathcal{F}(p)\slash \mathcal{F}^{BC}(p)} \tens{T}_{cpf} \lcpf \vec{n}_{cpf} \right)
    - \sum_{f\in \mathcal{F}^{BC}(p)} \tens{T}_{cpf}^{BC} \, \lcpf \vec{n}_{cpf}
    \label{eqn.pBC}
  \end{equation}
\item \textit{Velocity}: the prescribed velocity $v_{cpf}^{BC}:=\vec{v}_{cpf}^{BC}\cdot \vec{n}_{cpf}^n$ can be interpreted as a traction BCs. The equivalent traction $\tens{T}_{\vec{v}}^{BC}$ is given by
  \begin{eqnarray}
    \tens{T}_{\vec{v}}^{BC} &=& \Frac{\Mp^{-1} \left( \Sum_{f\in \mathcal{F}(p)} \tens{T}_{cpf} \, \lcpf \vec{n}_{cpf} 
      + \Mp \vec{u}_p \right) \cdot \vec{d}_p - \Sum_{f\in \mathcal{F}^{BC}(p)} v_{cpf}^{BC} \, \lcpf}{\Mp^{-1} \vec{d}_p \cdot \vec{d}_p}, \quad
    %\vec{d}_p = \Sum_{f\in \mathcal{F}^{BC}(p)} \lcpf \vec{n}_{cpf},
    \label{eqn.vBC1}
  \end{eqnarray}
  where $\vec{d}_p= \Sum_{f\in \mathcal{F}^{BC}(p)} \lcpf \vec{n}_{cpf}$ represents the corner vector associated to the boundary faces.
  Then, the node velocity is evaluated by considering
  \begin{equation}
    \vec{v}_p = \Mp^{-1} \left(  \Sum_{c\in  \mathcal{C}(p)} \Mcp \vec{v}_c - \sum_{f\in \mathcal{F}(p)} \tens{T}_{cpf} \, \lcpf \vec{n}_{cpf}   \right) - \tens{T}_{\vec{v}}^{BC} \vec{d}_p.
    \label{eqn.vBC}
  \end{equation}
\item \textit{Symmetry}: symmetry BC involves geometric considerations;
  either a symmetry plane defined by an orthonormal basis $(\vec{\tau}_1,\vec{\tau}_2)$,
  or a symmetry line along a direction vector $\vec{\tau}_1$, or even a symmetry point where we simply set $\vec{v}_p=\vec{0}$. 	
  In the case of a symmetry plane then the node velocity writes $\vec{v}_p=\alpha_1 \vec{\tau}_1 + \alpha_2 \vec{\tau}_2$ and the momentum balance equation becomes 
  \begin{equation}
    \Mp\vec{v}_p = \Mp (\alpha_1 \vec{\tau}_1 + \alpha_2 \vec{\tau}_2) = \left( \Mcpf \vec{v}_c- \sum_{f\in \mathcal{F}(p)} \tens{T}_{cpf} \, \lcpf \vec{n}_{cpf}   \right),
    \label{eqn.sBC1}
  \end{equation}
  which is solved by successive projection on $\vec{\tau}_1$ and $\vec{\tau}_2$.
  On a symmetry line one has $\vec{u}_p=\alpha_1 \vec{\tau}_1$, that is $\alpha_2=0$.
  %and the node velocity is obtained by
%	\begin{equation}
%	  \alpha_1\Mp \vec{\tau}_1 = \left(  \Mcpf \vec{v}_c-\sum_{f\in \mathcal{F}(p)} \tens{T}_{cpf} \, \lcpf \vec{n}_{cpf}  \right).
%	\label{eqn.sBC2}
%	\end{equation}
\end{itemize}
For further details and comments on BCs we refer the reader to \cite{Maire2007,MaireHD3D} and \cite{LAM2018}.

%---------------------------------------------%
%---------------------------------------------%
% --- B I B L I O G R A P H Y  ---------------% 
\bibliography{biblio}
\bibliographystyle{plain}
%---------------------------------------------%

\end{document}

%% file: Hyperelasticity.tex
\section{Updated Lagrangian hyperelastic modeling for isotropic materials}
\label{sec.equations} 
%---------------------------------------------%
In this section, we aim at recalling the conservation laws describing the time evolution of isotropic solid materials undergoing large deformations. The resulting conservation laws of mass, momentum and total energy shall be written under the updated Lagrangian form. The isotropic materials under consideration are characterized by an hyperelastic constitutive law. Namely, the Cauchy stress tensor is defined as being the derivative of the free energy with respect to the deformation tensor. In this framework, the material indifference principle and the thermodynamic consistency are automatically satisfied. The interested reader might refer for instance to \cite{Gurtin2010} for further developments on these subtle topics. For the sake of completeness, we recall hereafter some notions of kinematics that shall be useful for writing the conservation laws and the constitutive law.

\subsection{Kinematics} \label{ssec.kinematic} 

% ---------- FIG ---------------------
\begin{figure}[!htbp]
  \begin{center}
    \includegraphics[width=0.6\textwidth]{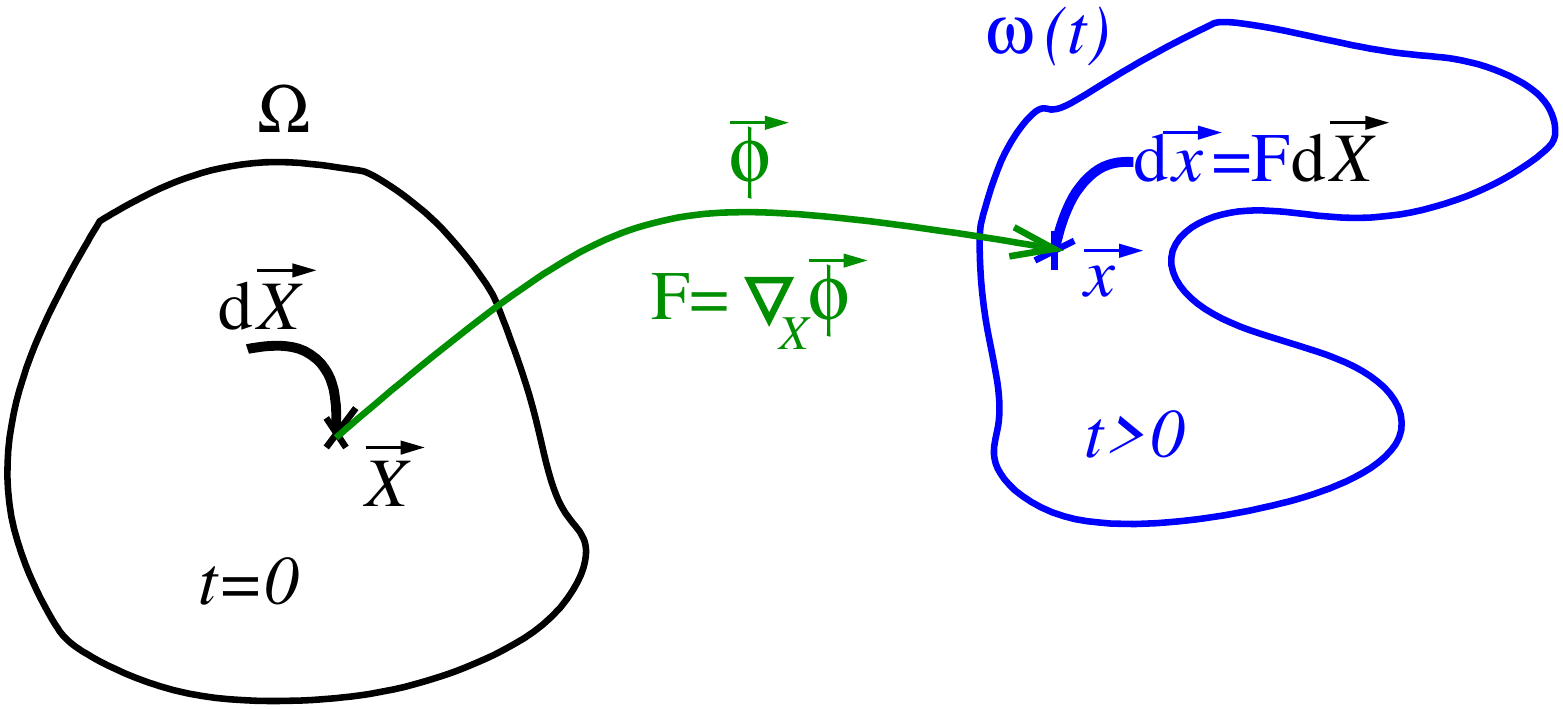}     
    \caption{ Sketch of the Lagrangian-Eulerian mapping relating a material Lagrangian point $\bm{X}$ at $t=0$ and a spatial Eulerian one $\bm{x}$ at $t>0$ through $\bm{\Phi}$, and its Jacobian $\tens{F}(\bm{X},t) = \pdd{\bm{\Phi}}{\bm{X}}( \bm{X},t )$. } 
    \label{fig:mapping}
  \end{center}
\end{figure}
% ---------- FIG ---------------------

\subsubsection{Lagrange-Euler mapping.} \label{sssec:mapping}
We consider a solid body in the $d$-dimensional Euclidean space occupying the region $\Omega$ in its initial configuration at time $t=0$, and the region $\omega(t)$ in its current configuration at time $t>0$. The motion of this body is characterized by the smooth function, $\bm{\Phi}$, that assigns to each material point $\bm{X}$ and time $t$ the spatial point $\bm{x}$ such that
\begin{align*}
& \Omega \;\longrightarrow \omega(t) \\
& \bm{X}\;  \longmapsto \bm{x} = \bm{\Phi}(\bm{X},t).
\end{align*}
This smooth function is the Lagrange-Euler mapping which relates the Lagrangian (material) point $\bm{X}$ to its Eulerian (spatial) counterpart $\bm{x}$. By definition, this mapping satisfies $\bm{\Phi}(\bm{X},0)= \bm{X}$ and its Jacobian, also named the deformation gradient, reads
$$ \tens{F}(\bm{X},t) = \nabla_{X}\bm{\Phi}(\bm{X},t),$$
where the symbol $\nabla_{X}$ denotes the gradient operator with respect to the Lagrangian coordinate. The determinant of the deformation gradient is denoted $J(\bm{X},t)=\det \left( \tens{F}(\bm{X},t) \right)$ and satifies $J(\bm{X},t=0)=1$ since $\tens{F}(\bm{X},t=0)=\tens{I}_{\text{d}}$ where $\tens{I}_{\text{d}}$ is the identity tensor. A continuity argument leads us to assume that $J(\bm{X},t)>0$ for all $t>0$, ensuring as such that $\bm{\Phi}$ is a one-to-one mapping.

A physical quantity can be expressed as well in terms of the Lagrangian coordinate as in terms of the Eulerian coordinate. More precisely, let $G(\bm{X},t)$ denotes the Lagrangian representation of a physical quantity. Its Eulerian representation reads $g(\bm{x},t)$. Obviously, these are two representations of the same physical quantity and, as a consequence, they fulfill the identities
$$g(\bm{x},t) = G \left[ \bm{\Phi}^{-1}( \bm{X},t) , t \right],\; \text{and} \; G(\bm{X},t) = g \left[ \bm{\Phi}( \bm{X},t) , t \right].$$
In what follows, the same notation shall be employed for both descriptions.

Time differentiating the mapping holding $\bm{X}$ fixed allows us to define the kinematic velocity
\begin{equation}
\label{eq:vel_mapping}
\bm{v}( \bm{X},t)= \pdd{\bm{\Phi}}{t}\vert_{\bm{X}}(\bm{X},t).
\end{equation}
Now, time differenting the identity $g(\bm{X},t)=g(\bm{\Phi}(\bm{X},t))$ holding $\bm{X}$ fixed and applying the chain rule yields
\begin{equation}
\label{eq:mat_der}
\pdd{g}{t}\vert_{\bm{X}}(\bm{X},t) = \pdd{g}{t}|_{\bm{x}}(\bm{x},t) + \bm{v}( \bm{X},t)\cdot \nabla_{x} g,
\end{equation}
where $\nabla_{x} g$ is the gradient of g with respect to the Eulerian coordinate, {\it i.e.}, $\nabla_{x}g=\frac{\partial g}{\partial \bm{x}}$. Thus, the Lagrangian time derivative is nothing but the material time derivative which is denoted
$$\mddt{g}(\bm{x},t)=\frac{\partial g}{\partial t}(\bm{x},t)+\bm{v} \cdot \nabla_{x}g.$$
\subsubsection{Measures of deformation}
\label{sssec:deformation}
Let us consider the infinitesimal material fiber $\rmd \bm{X}$ in the initial configuration which maps into $\rmd \bm{x} = \tens{F} \rmd \bm{X}$ through the motion. We express the streching of this infinitesimal fiber as follows
\begin{equation}
  \label{eq:rcgt}
\mathrm{d}\bm{x}\cdot \mathrm{d}\bm{x}-\mathrm{d}\bm{X}\cdot \mathrm{d}\bm{X}=(\tens{C}-\tens{I}_{\text{d}})\mathrm{d}\bm{x}\cdot \mathrm{d}\bm{x},
\end{equation}
where $\tens{C} = \tens{F}^t\tens{F}$ is the right Cauchy-Green tensor. On the other hand, noticing that $\rmd \bm{X}=\tens{F}^{-1}\rmd \bm{x}$, we also express the stretching of the infinitesimal fiber as follows
\begin{equation}
\label{eq:lcgt}
\mathrm{d}\bm{x}\cdot \mathrm{d}\bm{x}-\mathrm{d}\bm{X}\cdot \mathrm{d}\bm{X}=(\tens{I}_{\text{d}}-\tens{B}^{-1})\mathrm{d}\bm{X}\cdot \mathrm{d}\bm{X},
\end{equation}
where $\tens{B} = \tens{F}\tens{F}^{t}$ is the left Cauchy-Green tensor. The right and the left Cauchy-Green tensors are symmetric positive definite and share the same eigenvalues, refer to \cite{Gurtin2010}. These tensors are relevant measures of deformation since for a rigid rotation they boil down to the identity tensor. 
\subsubsection{Geometric conservation law (GCL)}
\label{sssec:GCL}
Time differentiating the deformation gradient, $\tens{F}=\nabla_{X}\bm{\Phi}$, recalling that the partial time derivative of the mapping is the kinematic velocity, $\bm{v}=\frac{\partial \bm{\Phi}}{\partial t}$, leads to the Geometric Conservation Law (GCL) written under total Lagrangian form
\begin{equation}
\label{eq:GCL}
\pddt{\tens{F}} - \nabla_{X}\bm{v} = 0,
\end{equation}
where $\tens{F}(X,0)=\tens{I}_{\text{d}}$. The GCL is supplemented with the compatibility constraint $\nabla_{X} \times \tens{F}=\bm{0}$, which ensures that the solution of the foregoing partial differential equation corresponds to the gradient of a mapping. Here, for any second order tensor $\tens{T}$, $\nabla_{X} \times \tens{T}$ denotes the rotational of $\tens{T}$. It is the tensor defined by $(\nabla_{X} \times T)\bm{a}=\nabla_{X} (\tens{T}^{t}\bm{a})$ for all constant vector $\bm{a}$. We note in passing that the compatibility constraint is an involutive constraint for the GCL. Namely, if this constraint is fulfilled initially, it will be satisfied for all time $t>0$. The satisfaction of this compatibility constraint at the discrete level is the cornerstone on which any proper discretization of the conservation laws written in total Lagrangian form should rely, refer to \cite{Vilar1}.

Introducing the material time derivative and applying the chain rule, we express the GCL under the updated Lagrangian form
\begin{equation}
\label{eq:GCLup}
\mddt{\tens{F}} - (\nabla_{x}\bm{v}) \tens{F} =0.
\end{equation}
Here, the deformation gradient and the velocity are viewed as functions of the spatial coordinate $\bm{x}$. The notation $\tens{L}=\nabla_{x} \bm{v}$ represents the velocity gradient tensor with respect to the current configuration. Employing this notation the updated Lagrangian form of the GCL reads
\begin{equation}
\label{eq:GCLL}
\mddt{\tens{F}} - \tens{L} \tens{F} = 0.
\end{equation}
Bearing this in mind, let us investigate two important consequences of the GCL that will be usefull in the sequel.

The first one is related to the time rate of change of the Jacobian $J=\det \tens{F}$ which is deduced from the GCL thanks to the chain rule
$$\mddt{(\det \tens{F})}=\frac{\partial (\det \tens{F})}{\partial \tens{F}} : \mddt{\tens{F}},\;\text{where}\; \frac{\partial (\det \tens{F})}{\partial \tens{F}}=(\det \tens{F}) \tens{F}^{-t}.$$
Here, the symbol $:$ denotes the inner product between tensors, {\it i.e.}, $\tens{A}:\tens{B}=\tr (\tens{A}^{t}\tens{B})$, where $\tr $ denotes the trace operator. Finally, substituting the GCL \eqref{eq:GCLL} into the foregoing equation yields the partial differential equation satisfied by the Jacobian of the deformation gradient
$$\mddt{J} - J \tr (\tens{L})  = 0.$$
Observing that $\tr (\tens{L})=\nabla_{x} \cdot \bm{v}$ leads to rewrite the time rate of change of the Jacobian as follows
\begin{equation}
\label{eq:jacobian}
\mddt{J}- J \nabla_{x}\cdot \bm{v}  = 0.
\end{equation}
The second one is related to the computation of the time rate of change of the left Cauchy-Green tensor, $\tens{B}=\tens{F}\tens{F}^{t}$, which reads
$$\mddt{\tens{B}}=\mddt{\tens{F}}\tens{F}^{t}+\tens{F}\mddt{\tens{F}^{t}}.$$
Substituting the expression of the time rate of change of $\tens{F}$ provided by the GCL into the foregoing equation leads to
\begin{equation}
  \label{eq:trcb}
  \mddt{\tens{B}}-\tens{L}\tens{B}-\tens{B}\tens{L}^{t}=0.
\end{equation}
The left-hand side of this equation is nothing but the Lie derivative of $\tens{B}$ ortherwise named the Oldroyd rate of $\tens{B}$ \cite{Gurtin2010}.
\subsection{Governing equations}
\label{ssec:eqs}
This section aims at briefly recalling the conservation laws and the constitutive law governing the time evolution of an isotropic material undergoing large deformations. The interested reader might consult \cite{Gurtin2010} for further details regarding their derivation. 
\subsubsection{Conservation laws}
Under the updated Lagrangian representation, the conservation laws of mass, momentum and total energy write
\begin{subequations}
  \label{eq:cl}
  \begin{align}
    &\rho \mddt{\tau}-\nabla \cdot \bm{v}=0,\label{eq:cl1}\\
    &\rho \mddt{\bm{v}}-\nabla \cdot \tens{T}=\bm{0},\label{eq:cl2}\\
    &\rho \mddt{e}-\nabla \cdot (\tens{T}\bm{v})=0. \label{eq:cl3}
  \end{align}
\end{subequations}
Here, the symbol $\mddt{}$ denotes the material derivative defined by \eqref{eq:mat_der}, $\rho >0$ is the mass density and $\tau=\frac{1}{\rho}$ the specific volume. The specific total energy is $e=\varepsilon + \frac12 \bm{v}^2$ where $\varepsilon$ denotes the specific internal energy. The Cauchy stress tensor, $\tens{T}$, is symmetric, {\it i.e.}, $\tens{T}=\tens{T}^{t}$, which ensures the conservation of angular momentum. Let us note that the nabla operator employed in the foregoing equations is expressed in terms of the Eulerian coordinate $\bm{x}$. 
This system of conservation laws written under Lagrangian updated representation is supplemented by the trajectory equation already introduced in \eqref{eq:vel_mapping}, which is rewritten under the form
\begin{equation}
\label{eqn.trajODE}
\mddt{\bm{x}} = \bm{v}(\bm{x}(t),t), \qquad \bm{x}(0)=\bm{X}.
\end{equation}
It is worth pointing out that \eqref{eq:cl1} is obtained by combining the Lagrangian mass conservation equation, $\mddt{(\rho J)}=0$ and the GCL \eqref{eq:jacobian}.
To close the foregoing system of conservation laws, it remains to provide a constitutive law for expressing the Cauchy stress tensor in terms of the deformation and a thermodynamic variable. This will be achieved in the next paragraph introducing the free energy $\Psi$. This thermodynamic potential is related to the specific energy, the absolute temperature, $\theta >0$, and the specific entropy $\eta$ by means of the classical thermodynamic relation
\begin{equation}
  \label{eq:psi}
  \Psi=\varepsilon-\theta \eta.
  \end{equation}
\subsubsection{Constitutive law for isotropic materials}
\label{sssec:frame_invariance}
The constitutive law is derived invoking the frame indifference principle and the compatibility with thermodynamics. This means that the constitutive equations should be invariant under changes of frame and satisfy the second law of thermodynamics \cite{Gurtin2010}. Here, the material under consideration is characterized by a free energy expressed in terms of the left Cauchy-Green tensor and the absolute temperature
$$\Psi \equiv \Psi(\tens{B},\theta).$$ 
Moreover, since this material is isotropic, its constitutive law is invariant under the group of all rotations acting in the spatial configuration. Thus, the theorem of representation of isotropic scalar function \cite{Gurtin2010} leads to the following expression of the free energy
\begin{equation}
  \label{eq:psi_iso}
  \Psi \equiv \Psi(I_1(\tens{B}),I_2(\tens{B}),I_3(\tens{B}),\theta).
\end{equation}
Here, $I_i(\tens{B})$ for $i=1,2,3$ are the principal invariants of the left Cauchy-Green tensor defined in Appendix~\ref{sec:pit}. 

Finally, the constitutive law provides the expressions of the Cauchy stress tensor and the specific entropy in terms of the free energy
\begin{equation}
\label{eq:const_law}
\tens{T}( \tens{B},\theta) = 2 \rho \left( \frac{\partial \Psi}{\partial \tens{B}} \right)_{\theta} \tens{B} , \qquad \text{and} \qquad
\eta( \tens{B},\theta) = - \left( \frac{\partial \Psi}{\partial \theta} \right)_{\tens{B}},
\end{equation}
where $\Frac{\partial \Psi}{\partial \tens{B}}$ is the tensor whose $ij$ component is $\Frac{\partial \Psi}{\partial \tens{B}_{ij}}$. Thanks to \eqref{eq:psi}, we observe that the specific internal energy $\varepsilon$ is also a function of the left Cauchy Green tensor and the temperature, {\it i.e.},  $\varepsilon(\tens{B},\theta) =\Psi(\tens{B},\theta) + \theta \, \eta(\tens{B},\theta)$.

The foregoing generic expression of the Cauchy stress tensor might be investigate further exploiting the isotropy of the material. Indeed, differentiating \eqref{eq:psi_iso} with respect to $\tens{B}$ and applying the chain rule leads to 
$$\left( \frac{\partial \Psi}{\partial \tens{B}} \right)_{\theta}=\left(\frac{\partial \Psi}{\partial I_1}\right)_{\theta} \tens{I}_{\text{d}}+ \left(\frac{\partial \Psi}{\partial I_2}\right)_{\theta} (I_1  \tens{I}_{\text{d}}-\tens{B})+\left(\frac{\partial \Psi}{\partial I_3}\right)_{\theta} I_3 \tens{B}^{-1},$$
where the derivative of the principal invariants of $\tens{B}$ with respect to $\tens{B}$ are recalled in Appendix~\ref{sec:pit}. Susbtituting the foregoing equation into the constitutive law provides us 
\begin{equation}
  \label{eq:cst_iso}
  \tens{T}=2\rho \left \{I_3 \left (\frac{\partial \Psi}{\partial I_3}\right)_{\theta}\tens{I}_{\text{d}}+\left [\left (\frac{\partial \Psi}{\partial I_1}\right)_{\theta}+I_1\left(\frac{\partial \Psi}{\partial I_2}\right)_{\theta} \right] \tens{B}-\left (\frac{\partial \Psi}{\partial I_2}\right)_{\theta}\tens{B}^{2}\right \}.
\end{equation}
This is the general expression of the Cauchy stress tensor for an isotropic hyperelastic material. It is quadratic with respect to the left Cauchy-Green tensor. Let us point out that the Cauchy stress tensor and the left Cauchy-Green tensor commute, {\it i.e.} $\tens{T}\tens{B}=\tens{B}\tens{T}$. This important property is the consequence of the material isotropy.

Its remains to check the consistency of this constitutive law with the second law of thermodynamics. First, differentiating the definition of the free energy \eqref{eq:psi} yields
\begin{align*}
  \theta \rmd \eta&= \rmd \varepsilon-\rmd \Psi -\eta \rmd \theta, \\
  &= \rmd \varepsilon -\frac{\partial \Psi}{\partial \tens{B}}:\rmd \tens{B} -\frac{\partial \Psi}{\partial \theta} \rmd \theta -\eta \rmd \theta,\;\text{since}\;\Psi=\Psi(\tens{B},\theta).
\end{align*}
Susbtituting the constitutive law \eqref{eq:const_law} in the above equation and recalling that $\varepsilon=e-\frac{1}{2}\bm{v}^{2}$ we arrive at the fundamental Gibbs relation
\begin{equation}
  \label{eq:gibbs}
  \theta \rmd \eta=-\frac{1}{2\rho}\tens{T}\tens{B}^{-1}:\rmd \tens{B}-\bm{v}\cdot \rmd \bm{v}+\rmd e.
\end{equation}
The Gibbs relation enables us to compute the time rate of change of entropy as follows
\begin{equation}
  \label{eq:gibbs2}
  \rho \theta \mddt{\eta}=-\frac{1}{2}\tens{T}\tens{B}^{-1}:\mddt{\tens{B}}-\rho \bm{v}\cdot \mddt {\bm{v}}+\rho \mddt {e}.
\end{equation}
  One the one hand, substituting the GCL \eqref{eq:trcb} into the first term of the right-hand side of \eqref{eq:gibbs2} leads to
  \begin{align*}
    \frac{1}{2}\tens{T}\tens{B}^{-1}:\mddt{\tens{B}}=&\frac{1}{2}\tens{T}\tens{B}^{-1}:(\tens{L}\tens{B}-\tens{B}\tens{L}^{t})\\
    =&\tens{T}:\tens{L},\;\text{since}\;\tens{T}\;\text{and}\;\tens{B}\;\text{commute}.
  \end{align*}
  On the other hand, substituting the conservation laws \eqref{eq:cl2} and \eqref{eq:cl3} into the remaining terms of the right-hand side of \eqref{eq:gibbs2} yields
  \begin{align*}
    -\rho \bm{v}\cdot \mddt {\bm{v}}+\rho \mddt {e}=&-\bm{v}\cdot \nabla \cdot (\tens{T})+\nabla \cdot (\tens{T}\bm{v}),\\
    =& \tens{T}:\nabla \bm{v}.
  \end{align*}
Here, we have employed the identity $\nabla \cdot (\tens{T}^{t}\bm{v})=\bm{v}\cdot \nabla \cdot (\tens{T})+\tens{T}:\nabla \bm{v}$. Finally, gathering the foregoing results and observing that $\tens{L}=\nabla \bm{v}$ we arrive at
\begin{equation}
  \label{eq:entrop}
  \rho \theta \mddt{\eta}=0.
  \end{equation}
This shows that system of conservation laws \eqref{eq:cl} is equipped with a supplementary conservation law which states that entropy is conserved along flow trajectories. Thus, constitutive law \eqref{eq:const_law} for isotropic materials is consistent with the second law of thermodynamics. Let us point out that the algebric manipulations which led to this result have been completed under the smoothness assumption of the flow variables. In the presence of discontinuities such as shock waves, the entropy conservation law turns into the entropy inequality
\begin{equation}
  \label{eq:entropineq}
  \rho \theta \mddt{\eta} \geq 0.
  \end{equation}
\subsubsection{Volumetric shear strain decomposition}
\label{ssec:neo-hookean}
We want to study materials that can sustain only limited shear strain but respond elastically to large change in volume. Following \cite{Plohr2012}, we introduce the additive decomposition of the free energy into a volumetric part and a shear part. This in turn provides the additive decomposition of the Cauchy stress into a spherical part, which is nothing but the pressure, and a deviatoric part. To construct this addtive decomposition, we start by introducing the multiplicative decomposition of the deformation gradient tensor, $\tens{F}$, into a volumetric and an isochoric parts. The volumetric part is equal to $J^{\frac{1}{3}}\tens{I}_{\text{d}}$, whereas its isochoric part reads $\overline{\tens{F}}=J^{-\frac{1}{3}}\tens{F}$. This part of the deformation gradient is volume preserving since $\det (\overline{\tens{F}})=1$. This in turn implies that the isochoric part of the left Cauchy-Green tensor reads $\overline{\tens{B}}=J^{-\frac{2}{3}} \tens{B}$. Bearing this decomposition in mind, the expression of the free energy \eqref{eq:psi_iso} turns into
\begin{equation}
  \label{eq:psi_iso_bis}
  \Psi \equiv \Psi(J,I_1(\overline{\tens{B}}),I_2(\overline{\tens{B}}),\theta).
\end{equation}
The dependence of the free energy on $I_3$ is held by $J$ since $I_3(\overline{\tens{B}})=\det(\overline{\tens{B}})=1$.
Now, we decompose this latter expression of the free energy  into
\begin{equation}
\label{eq:decomp_psi}
\Psi = \Psi_v( J, \theta) + \Psi_s(\overline{I}_1,\overline{I}_2,\theta),
\end{equation}
where $\Psi_v$ and $\Psi_s$ denote respectiveley the volumetric and the shear parts of the free energy knowing that $\overline{I}_{1}=I_1(\overline{\tens{B}})$ and $\overline{I}_{2}=I_2(\overline{\tens{B}})$ are the principal invariants of the isochoric part of the left Cauchy-Green tensor $\overline{\tens{B}}$, refer to Appendix~\ref{sec:pit} for the definition of the principal invariants of a tensor.

Finally, substituting the volumetric/shear decomposition of the free energy into the constitutive law \eqref{eq:const_law} and after some algebra we arrive at
\begin{equation}
  \label{eq:T_psi}
  \tens{T}=\rho J \left (\frac{\partial \Psi_v}{\partial J}\right)_{\theta}\tens{I}_{\text{d}} +2\rho \left [\left (\frac{\partial \Psi_s}{\partial \overline{I}_{1}}\right)_{\theta}\overline{\tens{B}}_{0}-\left (\frac{\partial \Psi_s}{\partial \overline{I}_{2}}\right)_{\theta}(\overline{\tens{B}}^{-1})_{0} \right].
 \end{equation}
Here, for a tensor, the superscript $0$ denotes its deviatoric part. Thus, $\tens{T}_{0}$ is the deviatoric part of the Cauchy stress tensor defined by $\tens{T}_{0}=\tens{T}-\frac{1}{3}\tr(\tens{T})\tens{I}_{\text{d}}$ and obviously $\tr(\tens{T}_{0})=0$. Let us note that the foregoing expression of Cauchy stress tensor in terms of $\overline{\tens{B}}^{-1}$ has been obtained thanks to the Cayley-Hamilton theorem, refer to Appendix~\ref{sec:pit}, which allows to write $\overline{\tens{B}}^{-1}=\overline{\tens{B}}^{2}-\overline{I}_{1}\overline{B}+\overline{I}_{2} \tens{I}_{\text{d}}$. Observing \eqref{eq:T_psi}, we arrive at the conclusion that the Cauchy stress decomposes into a spherical and a deviatoric parts which are respectively defined by
\begin{subequations}
  \label{eq:cst}
  \begin{align}
    p=& -\rho J \left(\pdd{\Psi_v}{J} \right)_\theta,\;\text{spherical part}\label{eq:cst1}\\
    \tens{T}_0 =& 2 \rho \left(\pdd{\Psi_S}{\overline{I}_1} \right)_\theta \overline{\tens{B}}_0 - 2\rho \left(\pdd{\Psi_S}{\overline{I}_2} \right)_\theta (\overline{\tens{B}}^{-1})_0,\;\text{deviatoric part}.\label{eq:cst2}
  \end{align}
\end{subequations}
Here, $p=p(J,\theta)$ is nothing but the pressure and we point out that $\tens{T}_{0}=\tens{T}_0(\overline{I}_1,\overline{I}_2,\theta)$.
\begin{remark}[Hyperelasticity versus hypoelasticity]
  Hyperelasticity relies on the definition of a free energy which allows to express the deviatoric part of the Cauchy stress in terms of the deviatoric part of the left Cauchy-Green tensor. This framework provides a constitutive law fulfilling
  \begin{itemize}
  \item The material frame indifference principle;
  \item The thermodynamic consistency with the second law.
  \end{itemize}
On the other hand, for hypoelasticity, refer for instance to \cite{Maire_elasto}, the constitutive law is written under incremental form. Namely, the time rate of change of the deviatoric stress is expressed in terms of the deviatoric part of the strain rate tensor. The enforcement of the principle of material frame indifference relies on the use of a somewhat arbitrary objective stress rate such as the Jaumann rate, refer to \cite{Gurtin2010}. Moreover, the use of objective stress rate makes appear non conservative terms which render the mathematical analysis of discontinuous solutions quite delicate. This framework does not allow the fulfillment of thermodynamic consistency. Indeed, for smooth elastic flows the entropy is not conserved.
\end{remark}
According to the constitutive law \eqref{eq:const_law} the volumetric/shear decomposition of the free energy also induces a similar additive decomposition of the specific entropy $\eta=\eta_v+\eta_s$ where
\begin{subequations}
  \label{eq:eta_Psi}
  \begin{align}
    \eta_v(J,\theta)=&-\left (\frac{\partial \Psi_v}{\partial \theta}\right)_{J},\;\text{volumetric part}\label{eq:eta_Psi1}\\
    \eta_s(\overline{I}_1,\overline{I}_2,\theta)=&-\left (\frac{\partial \Psi_s}{\partial \theta}\right)_{\overline{I}_1,\overline{I}_2},\;\text{shearing part}.\label{eq:eta_Psi2}
  \end{align}
\end{subequations}
Gathering the foregoing results and recalling that, $\varepsilon = \Psi + \theta \eta$, leads to 
$$\varepsilon =  \Psi_v + \Psi_s + \theta ( \eta_v+\eta_s ) = (\Psi_v + \theta  \eta_v)  + (\Psi_s + \theta  \eta_s).$$
Thus, it is natural to introduce the volumetric and the shearing parts of the specific internal energy as follows
\begin{subequations}
  \label{eq:eta_epsi}
  \begin{align}
    \varepsilon_v(J,\theta)=&\Psi_v(J,\theta)+\theta \eta_v(J,\theta),\label{eq:eta_epsi1}\\
    \varepsilon_s(\overline{I}_1,\overline{I}_2,\theta)=& \Psi_s(\overline{I}_1,\overline{I}_2,\theta)+\theta \eta_s(\overline{I}_1,\overline{I}_2,\theta).\label{eq:eta_epsi2}
  \end{align}
\end{subequations}
\begin{remark} [About other thermodynamic potentials]
The thermoelastic response of the material could have been defined choosing internal energy, $\varepsilon\equiv \varepsilon(\tens{B},\eta)$, as a thermodynamic potential to further derive the constitutive law, refer for instance to \cite{Gavriluk08,Kluth10}. However, as noticed in \cite{Plohr2012}, such a choice is inappropriate because it would imply that the absolute temperature $\theta$ (which is an intensive thermodynamic quantity) is a sum of volumetric/shear contributions. Moreover, the choice of the absolute temperature as an independent variable is more convenient since the notion of stress depending on temperature is more familiar, mostly because the temperature can easily be measured with classical devices such as thermometers. 
\end{remark}
\subsubsection{Examples of constitutive laws}
Let us point out that the volumetric/shear decomposition allows us to define separately the pressure by introducing an hydrodynamic equation of state characterized by the volumetric free energy $\Psi_v=\Psi_v(J,\theta)$. The pressure and the internal energy are expressed by means of classical thermodynamic relations
\begin{equation}
  \label{eq:eoshyd}
  p(\tau,\theta)=-\rho^{0}\left (\frac{\partial \Psi_v}{\partial J}\right)_{\theta},\;\;\varepsilon_v(J,\theta)=\Psi_v(J,\theta)-\theta \left(\frac{\partial \Psi_v}{\partial \theta}\right)_{J},
\end{equation}
where $\rho^{0}>0$ denotes the initial mass density of the solid. In what follows, for numerical applications, we shall make use of the volumetric free energy
\begin{equation}
\label{eq:EOS_neoHook}
\Psi_v = \frac{\mu}{4\rho^0} \left( (J-1)^2 + (\log J)^2 \right),
\end{equation}
which leads to the pressure $p=-\frac{\mu}{2}(J-1+\frac{\log J}{J})$ and the volumetric internal energy $\varepsilon_v=\Psi_v$. 
Apart from this equation of state, we shall also utilize the stiffened gas equation of state, which writes under the incomplete form
\begin{equation}
 \label{eq:EOS_stiffened_gas}
\varepsilon_v = \frac{p+\gamma p_\infty}{ (\gamma - 1) \rho},
\end{equation}
where $\gamma$ and $p_\infty$ are material-dependent parameters. More generaly, one can utilizes his favorite equation of state regardless of the shearing free energy choice. However, one shall always choose at least a convex equation of state to ensure the hyperbolicty of the hydrodynamic part of the system of conservation laws.

Regarding the shear part of the free energy we use the family of rank-one convex stored energies proposed by \cite{Gavrilyuk2015}
\begin{equation}
\label{eq:familiy_energies}
\Psi_s (\overline{I}_{1},\overline{I}_{2})= \Frac{\mu}{4\rho^0}\left [-2a (\overline{I}_{1}-3)+\frac{(1+a)}{3} (\overline{I}_{2}^{2}-9) \right],
\end{equation}
where is an adjustable parameter. For $a\in [-1,\frac{1}{2}]$, it is shown in \cite{Gavrilyuk2015} that the resulting system of conservation laws is hyperbolic. For the numerical applications, we shall consider the particular case $a=-1$ which corresponds to neo-Hookean materials. In this case, the shear part of free energy reads $\Psi_s=\frac{\mu}{2\rho^0}( \overline{I}_1 - 3)$ and thus the deviatoric part of the Cauchy stress tensor is given by
\begin{equation}
  \label{eq:cstnh}
\tens{T}_0 =\frac{\mu}{J}\overline{\tens{B}}_0,
\end{equation}
where $\overline{\tens{B}}_0=\overline{\tens{B}} -\frac{1}{3} \tr (\overline{\tens{B}}) \tens{I}_{\text{d}}$.

Finally, material mechanical properties are often described in terms of Young modulus $E$, Poisson ration $\nu$ and shear modulus $\mu$, which also corresponds to the second Lam{\'e} coefficient. These parameters are linked as follows:
\begin{equation}
\label{eq:Enu_relation}
\mu = \Frac{E}{2 \, (1+\nu)}.
\end{equation}
In this paper, the numerical simulations will be carried out mainly with the neo-Hookean hyperelastic constitutive law, however
we might also employ the non linear constitutive law \eqref{eq:familiy_energies} in the case $a=0$ for comparison purposes.
\subsection{Summary: Updated Lagrangian hyperelasticity for isotropic materials}
We summarize the set of partial differential equations governing the time evolution of the isotropic hyperelastic material under consideration. The conservation laws of mass, momentum and total energy read
\begin{align*}
   &\rho \mddt{\tau}-\nabla \cdot \bm{v}=0,\\
    &\rho \mddt{\bm{v}}-\nabla \cdot \tens{T}=\bm{0},\\
    &\rho \mddt{e}-\nabla \cdot (\tens{T}\bm{v})=0.
  \end{align*}
The Cauchy stress tensor is symmetric, {\it i.e.}, $\tens{T}=\tens{T}^{t}$. It is obtained deriving the free energy with respect to the left Cauchy-Green tensor $\tens{B}$. Assuming a volumetric/shear decomposition of the free energy, $\Psi=\Psi_v+\Psi_s$, the Cauchy stress tensor reads
$$\tens{T}=\rho J \left (\frac{\partial \Psi_v}{\partial J}\right)_{\theta}\tens{I}_{\text{d}} +2\rho \left [\left (\frac{\partial \Psi_s}{\partial \overline{I}_{1}}\right)_{\theta}\overline{\tens{B}}_{0}-\left (\frac{\partial \Psi_s}{\partial \overline{I}_{2}}\right)_{\theta}(\overline{\tens{B}}^{-1})_{0} \right].$$
Here, $\overline{\tens{B}}=J^{-\frac{2}{3}} \tens{B}$ denotes the isochoric part of the left Cauchy-Green tensor and $\overline{I}_{1}$, $\overline{I}_{2}$ are respectively its first and second invariants. We note also that $\Psi_v=\Psi_v(J,\theta)$ and $\Psi_s=\Psi_s(\overline{I}_{1},\overline{I}_{2},\theta)$. By construction, the foregoing constitutive law satisfies the material frame indifference principle and is thermodynamically consistent which allows to write the Gibbs identity
$$\theta \rmd \eta=-\frac{1}{2\rho}\tens{T}\tens{B}^{-1}:\rmd \tens{B}-\bm{v}\cdot \rmd \bm{v}+\rmd e.$$

This system of physical conservation laws is completed by the geometrical conservation law expressing the time rate of change of the left Cauchy-Green tensor
$$\mddt{\tens{B}}-\tens{L}\tens{B}-\tens{B}\tens{L}^{t}=0,$$
where $\tens{L}=\nabla \bm{v}$ is the Eulerian velocity gradient tensor.

It is remarkable to note that updated Lagrangian isotropic hyperelasticity requires only the knowledge of the left Cauchy-Green tensor. 
\begin{remark}[Physical admissibility]
\label{rem:PAD}
The physical admissibility property is defined by a set of so-called admissible states such that the material vector determines a valid state according to the conservation and constitutive laws.  
If the vector of variables is $\Q=(\tau,\bm{v},e, \tens{B})$ supplemented with its relationships with derived variables, $\varepsilon$, $\tens{L}$, etc. in the hyper-elastic model considered in this work the physically admissible set $\mathcal{A}$ is
\begin{equation}
\label{eqn:admissible_set}  
\mathcal{A} = \left\{ \Q \;\text{s.t.} \;\tau>0 \; \text{and}  \; \varepsilon=e-\frac12 \bm{v}^2 >0  \; \text{and}  \; \theta>0\; \text{and}\; \rho \theta \mddt{\eta}\geq 0 \right\} ,
\end{equation}
\end{remark}

%% file: Discretization.tex
	\section{Finite volume discretization}
	\label{sec.numethod} 
	%---------------------------------------------%
	
Here, $\Omega(t) \subset \mathbb{R}^d$ denotes the time-dependent polygonal/polyhedral volume in current configuration  in $d\in[2,3]$ space dimensions and $\partial \Omega(t)$ its surface defined by the outward pointing unit normal vector $\vec{n}$.

\subsection{Mesh and notation}
\label{ssec:mesh}
The computational domain $\Omega(t)$ is discretized at time $t$
by a set of non-overlapping control volumes (polygonal/polyhedral cells),
%in 2D/3D respectively), 
each denoted by $\omega(t)$. 
$N_E$ denotes the total number of elements/cells in the domain and a cell is referred to with index $c$, that is $\omega_c(t)$.
We also refer to a vertex/point with index $p$. Moreover the set of points of a cell is denoted by $\mathcal{P}(c)$ and the set of cells sharing a giving point $p$ is $\mathcal{C}(p)$. 
Next the set of the faces of a cell is $\mathcal{F}(c)$
and the set of faces sharing a node $p$ is $\mathcal{F}(p)$. 
Likewise the sets of edges of a cell is $\mathcal{E}(c)$, and impinging at a common point is denoted by $\mathcal{E}(p)$. 

\begin{figure}[!htbp]
  \begin{center}
    \begin{tabular}{ccc} 
      %\hspace{-1.7cm}
      %\includegraphics[width=0.37\textwidth]{Te_2D}  &           
      \hspace{-1.5cm}
      \includegraphics[width=0.37\textwidth]{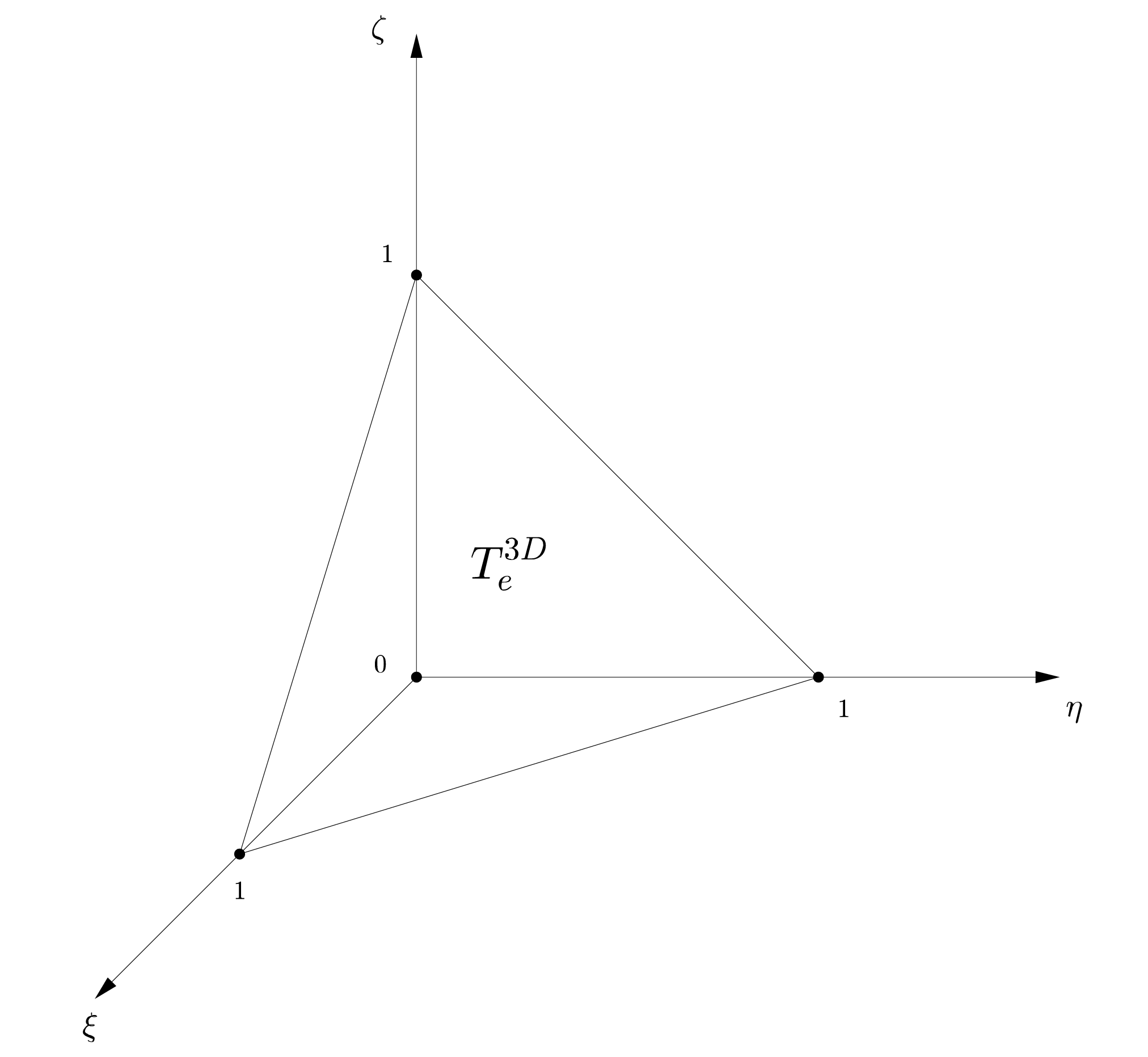}  &           
      \includegraphics[width=0.34\textwidth]{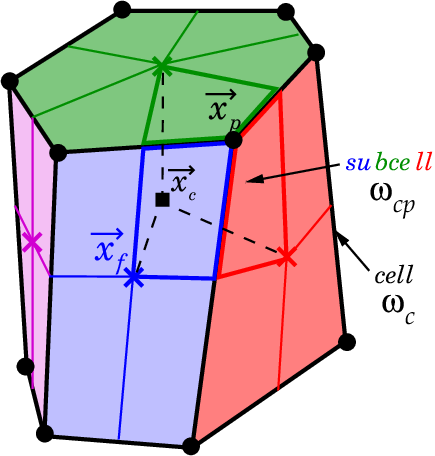}  \\  
    \end{tabular} 
    \caption{Left: Reference simplicial element $\omega_e$ in coordinates ${\vec{xi}} = (\xi, \eta, \zeta)$ for $d=3$ ---
      Right: Polyhedral cell $\omega_c$, subcell $\omega_{cp}$ and  geometrical face/cell/point centers. } 
    \label{fig:refSystem}
  \end{center}
\end{figure}

For any discrete time $t^n$, $n\in \mathbb{N}$, the union of all elements $\omega_c^n:=\omega_c(t^n)$ paving $\Omega(t^n)$ is called the \textit{current mesh configuration} $\mathcal{T}^n_{\Omega}$ of the domain 
\begin{equation}
\mathcal{T}^n_{\Omega} = \bigcup \limits_{c=1}^{N_E}{\omega^n_c}. 
\label{eqn:meshdef}
\end{equation}
Each control volume defined in the \textit{physical} space $\vec{x}=(x,y,z)$ can be mapped onto a reference element $T_e^{3D}$ in the reference coordinate system ${\vec{\xi}}=(\xi,\eta,\zeta)$ in 3D, see figure~\ref{fig:refSystem}.
In 2D the third components of $\vec{x}$ and ${\vec{\xi}}$ are maintained constant. 
%, see Figure \ref{fig.refSystem}. The reference element is the unit triangle in 2d  or the unit tetrahedron in 3d for simplex elements\footnote{In the case of polygonal or polyhedral cell, it is split into simplex elements.}.
%The spatial mapping in 3d reads
%\begin{equation} 
%\mathbf{x} = \mathbf{X}^n_{1,c} + 
%\left( \mathbf{X}^n_{2,c} - \mathbf{X}^n_{1,c} \right) \vec{x}i + 
%\left( \mathbf{X}^n_{3,c} - \mathbf{X}^n_{1,c} \right) \eta + 
%\left( \mathbf{X}^n_{4,c} - \mathbf{X}^n_{1,c} \right) \zeta,
%\label{eqn.xi} 
%\end{equation}
%where $\mathbf{X}^n_{k,c} = (X^n_{k,c},Y^n_{k,c},Z^n_{k,c})$ represents the vector of physical spatial coordinates of the $k$-th vertex of element $\omega^n_c$ for $k=1,2,3,4$.
% ---------- FIG ---------------------
% ---------- FIG ---------------------

%---------------------------------------------%
\subsubsection{Geometrical entities}
%---------------------------------------------%
The center of the cell is its centroid $\vec{x}_c$ and the center of a face $f$ is the iso-barycenter of the points defining the cell: $\vec{x}_f=\frac{1}{|\mathcal{P}(f)| }\Sum_{p\in\mathcal{P}(f)} \vec{x}_p$, where $|\mathcal{S}|$ denoted the cardinal of any set $\mathcal{S}$. 

Given a cell $c$ and a point $p$ we define a unique object called subcell, referred to with double index $cp$ which is the unique geometrical object linking a cell center $\vec{x}_c$, one of its point $\vec{x}_p$ and the face centers $\vec{x}_f$ for all face $f\in \mathcal{F}(c) \cap \mathcal{F}(p)$. 
In 3D the subcell is a hexaedron with possibly non-planar faces, in 2D it is a quadrangle. Further denoted by $\omega_{cp}$, its volume is referred to as $|\omega_{cp}|$, see figure~\ref{fig:refSystem}.
Consequently a cell $\omega_c$ is a collection of subcells: $\omega_c = \bigcup_{p\in \mathcal{P}(c)} \omega_{cp}$, each being considered as Lagrangian objects.
A dual cell $\omega_p$ is the collection of subcells sharing $\vec{x}_p$ as a node: $\omega_p = \bigcup_{c\in \mathcal{C}(p)} \omega_{cp}$. \\
In a Lagrangian framework the mass of a subcell and cell, $m_{cp}$, $m_c$ respectively, are constant in time and equal to
\bea
m_{cp} = \Int_{\omega_{cp}(t)} \rho^0(\vec{x}) \, \d v,\quad
m_c = \Int_{\omega_c(t)} \rho^0(\vec{x}) \, \d v = \Sum_{p\in \mathcal{P}(c)}m_{cp} , 
\eea
where $\rho^0(\vec{x})\equiv \rho(\vec{x},t=0)$ is the initial density distribution, and $\d v$ refers to the integral measure over volume. The mass of a dual cell, $m_p$, is the sum of the subcell masses in the dual cell.

An important geometrical object is the the so-called corner vector 
$\lncp$ which formal definition is given by
\bea \label{eq:def_lcpncp}
\lncp = \pdd{|\omega_c|}{\vec{x}_p} .
\eea
$\lcp$ represents a $(d-1)$-measure (length in 2D, area in 3D) and $\ncp$ is a unit outward pointing vector. Algebraic manipulations of (\ref{eq:def_lcpncp}) may convince the reader that the corner vector is the sum of the face outward pointing normal vectors for all face $f\in\mathcal{F}(p)$ of the current cell $c$ impinging on node $p$.
A cell being a close contour, we have the fundamental property of the corner vector
\bea
\Sum_{p\in \mathcal{P}(c)} \lncp = \vec{0}.
\eea

\subsubsection{Conservative and constitutive discrete variables}
The time dependent conserved or constitutive variables are the cell-centered approximate mass-averaged values gathered into vector $\Q_c(t)=(\tau_c(t),\vec{v}_c(t), e_c(t),\tens{B}_c(t))$.
%that is $\Q_c = \frac{1}{m_c} \Int_{\omega_{c}(t)} \Q(\vec{x},t) \d v$.
For a vector or a tensor the previous equation should be understood as component-wise. We also use in this work a point-wise velocity field $\vec{v}_p$ which represents the velocity of point $p$ and also the mean velocity in the dual cell $\omega_p(t)$: $\vec{v}_p(t) = \vec{v}( \vec{x}_p,t) = \frac{1}{m_p} \Int_{\omega_{p}(t)} \rho(\vec{x},t) \vec{v}(\vec{x},t) \d v$.
At last the density or specific volume could also be subcell centered representing \textit{de facto} the mean value over $\omega_{cp}(t)$: 
$\rho_{cp}(t) = \frac{1}{|\omega_{cp}(t)|} \Int_{\omega_{cp}(t)} \rho(\vec{x},t) \d v$. \\
For now one we implicitly assume the dependence on time and to lighten the notation we omit it.

%---------------------------------------------%
\subsection{Discrete divergence and gradient operators}
%---------------------------------------------%
Considering the discrete point-wise vector field $\vec{v}_p$ we define the cell-centered discrete divergence and adjoint gradient operators as
\bea \label{eq:div_grad_v}
(\DIV{\vec{v}})_c = \frac{1}{|\omega_c|}\Sum_{p\in \mathcal{P}(c)} \lncp \cdot \vec{v}_p, \qquad
\tens{L}_c := (\nabla{\vec{v}})_c = \frac{1}{|\omega_c|}\Sum_{p\in \mathcal{P}(c)} \lncp \otimes \vec{v}_p .
\eea
%hence we have $\tens{L}_c =(\nabla{\vec{v}})_c$. 
The discrete gradient of a scalar quantity like the cell-centered pressure $p_c$ is given by
\bea \label{eq:grad_p}
	\GRAD{p}_c = \frac{1}{|\omega_c|}\Sum_{p\in \mathcal{P}(c)} p_c \lncp .
\eea
These operators are nowadays classical in cell-center Lagrangian scheme community, see for instance \cite{Review_Handbook_16,LAM2018}.

%---------------------------------------------%
\subsection{Semi-discretization in space}
%---------------------------------------------%

\subsubsection{Conservation laws - GCL, momentum and total energy}
	
%\paragraph{Geometrical Conservation Law (GCL)}
The geometrical conservation law (GCL) is a fundamental consistency property in Lagrangian framework. Indeed it states that the discrete motion of all the points $p$ of a given cell $\omega_c$ with the trajectory equations 
\bea 
\label{eq:trajectory}
\mddt{\vec{x}_p} = \vec{v}_p,
\eea
is consistent with the volume conservation law \eqref{eq:cl1}.
Since $m_c \, \tau_c = |\omega_c|$ and taking into account the definition of corner vectors and discrete divergence, it is classical to infer the discrete version of the volume conservation law which is compatible with the GCL

%\paragraph{Momentum and total energy via subcell force discretization}
Moreover if we introduce the so-called subcell force $\vec{f}_{cp}$, which is the traction force attached to subcell $\omega_{cp}$, we can write the discrete version of the conservation laws as \cite{Mai11_subcell,LAM2018}:
\beqa 
\label{eq:tau_c}
m_c \mddt{\tau_c} - \Sum_{p\in \mathcal{P}(c)} \lncp \cdot \vec{v}_p &=& 0 , \\
\label{eq:v_c}
m_c \mddt{\vec{v}_c} - \Sum_{p\in \mathcal{P}(c)} \vec{f}_{cp} &=& \vec{0} , \\
\label{eq:e_c}
m_c \mddt{e_c} - \Sum_{p\in \mathcal{P}(c)} \vec{f}_{cp}\cdot \vec{v}_p &=& 0 .
%\label{eq:B_c}
%  \mddt{\tens{B}_c} - \tens{L}_c\tens{B}_c - \tens{B}_c\tens{L}^t_c &=& 0 .
\eeqa
Moreover the discrete version of \eqref{eq:trcb} is given by
\bea \label{eq:B_discr}
\mddt{\tens{B}_c} - \tens{L}_c \tens{B}_c - \tens{B}_c \tens{L}_c^t = 0.
\eea
One remarks that the subcell force is the last unknown of the previous discretization, our goal is to provide a compatible and consistent definition of it according to the conservation and constitutive laws. 
We refer to \cite{Mai11_subcell,LAM2018} for some details of the consequences of such a discretization, in particular the conservation properties when the hydrodynamics system of conservation law is solely considered.

\subsubsection{Semi-discrete entropy analysis - Subcell force}

The constitutive law leads to the definition of the following discrete Cauchy stress tensor: $\tens{T}_c = 2 \rho_c \pdd{\Psi}{\tens{B}}(\tens{B}_c,\theta_c) \tens{B}_c$.
Starting from the Gibbs identity \eqref{eq:gibbs} let us compute the time evolution of the entropy 
\begin{eqnarray} \label{eq:theta_deta_discr}
	m_c \theta_c \mddt{\eta_c} =
	 	-\frac{1}{2} |\omega_c| \tens{T}_c \tens{B}_c^{-1} : \mddt{\tens{B}_c} - m_c \vec{v}_c\cdot \mddt{\vec{v}_c} + m_c \mddt{e_c}.
\end{eqnarray}
Each term of the right hand side can be replaced by a more appropriate form for our analysis using (\ref{eq:v_c}), (\ref{eq:e_c}) and
\bea
\nn -\frac{1}{2}\tens{T}_c \tens{B}_c^{-1} : \mddt{\tens{B}_c}  =   -\tens{T}_c : \tens{L}_c = - \tens{T}_c : \Sum_{p\in \mathcal{P}(c)} \lcp \vec{v}_{p} \otimes \vec{n}_{pc}, %\\
%	\nn - m_c \vec{v}_c\cdot \mddt{\vec{v}_c} &=& - \Sum_{p\in \mathcal{P}(c)} \vec{f}_{cp}, \\
%	\nn m_c \mddt{e_c} &=& \Sum_{p\in \mathcal{P}(c)} \vec{f}_{cp} \cdot \vec{v}_{p}. 
\eea
which after substitution yields
\beqa 
\nn \hspace{-0.5cm}
	 m_c \theta_c \mddt{\eta_c} 
	&=&
		- \tens{T}_c : \Sum_{p\in \mathcal{P}(c)} \lcp \vec{v}_{p} \otimes \vec{n}_{pc}
		+\Sum_{p\in \mathcal{P}(c)} \vec{f}_{cp} \cdot (\vec{v}_{p} - \vec{v}_c) \\
\nn	 \hspace{-0.5cm}
   &=&
		-\Sum_{p\in \mathcal{P}(c)} (\vec{v}_{p}-\vec{v}_c) \cdot \tens{T}_c \vec{n}_{pc} + \Sum_{p\in \mathcal{P}(c)} \vec{f}_{cp} \cdot (\vec{v}_{p} - \vec{v}_c) 
	=
		\Sum_{p\in \mathcal{P}(c)} \left( -\lcp \tens{T}_c \vec{n}_{pc} + \vec{f}_{pc} \right)\cdot (\vec{v}_{p}-\vec{v}_c).
\eeqa
Therefore in order to ensure a proper entropy dissipation we propose to design
\bea \label{eq:subforce}
	\vec{f}_{pc} = \lcp \tens{T}_c \vec{n}_{pc} + \Mcp(\vec{v}_{p}-\vec{v}_c),
\eea
where the subcell matrix $\Mcp$ is symmetric positive definite.
And we easily verify that
\bea \label{eq:theta_deta_discr2}
	 m_c \theta_c \mddt{\eta_c} = 
	 	\Sum_{p\in \mathcal{P}(c)} \Mcp(\vec{v}_{p}-\vec{v}_c) \cdot (\vec{v}_{p}-\vec{v}_c) \geq 0,
\eea
which satisfies the second law of thermodynamics. 
Now it remains to determine the subcell matrix $\Mcp$, which genuinely characterizes the numerical scheme. 
Several possibilities have already been explored by different authors in \cite{Despres2005,Maire2007,Maire2009,Despres2009,ShashkovCellCentered} among others.

\subsubsection{Nodal solver - Subcell matrix}
Since the seminal works of Despres \textit{et al} \cite{Despres2005} and Maire \textit{et al} \cite{Maire2007}, a so-called nodal solver has become a classical tool for many cell-centered Lagrangian numerical schemes. 
A nodal solver could be interpreted as a local approximate multidimensional Riemann solver at a given node of the mesh. 
Our first-order discretization strictly follows the nodal solver of the Eucclhyd scheme proposed in \cite{Maire2007}. %,Maire2009,phmbn09,Maire2011}. 
It computes the nodal velocity $\vec{v}_p$ given the physical states in the surrounding cells by means of 1D half-Riemann problems invoking the conservation of momentum (or total energy). This, along with the definition of the subcell force, imply that for any point $p$ neglecting the boundary conditions
\bea \label{eq:sum_fcp}
\Sum_{c\in  \mathcal{C}(p)} \vec{f}_{pc} = \vec{0},
\eea
 yielding after substitution into (\ref{eq:subforce})
\beqa \nn
\Sum_{c\in  \mathcal{C}(p)}\lcp\tens{T}_c \vec{n}_{pc} + \Mcp(\vec{v}_{p}-\vec{v}_c) 	
&=& 
\Sum_{c\in  \mathcal{C}(p)}\lcp\tens{T}_c \vec{n}_{pc} + \left(\Sum_{c\in  \mathcal{C}(p)}\Mcp\right) \vec{v}_{p} 
- \Sum_{c\in  \mathcal{C}(p)}\Mcp \vec{v}_c 
= \vec{0}.
\eeqa
As a consequence we can compute the nodal velocity as the solution of the following linear system
\bea\label{eq:nodal_solver}
	\Mp \vec{v}_{p} =  \Sum_{c\in  \mathcal{C}(p)}\Mcp \vec{v}_c - \Sum_{c\in  \mathcal{C}(p)}\lcp\tens{T}_c \vec{n}_{pc}, \qquad \text{and} \qquad \Mp =\Sum_{c\in  \mathcal{C}(p)}\Mcp.
\eea
Notice that $\Mp$ is symmetric positive definite and, thus, invertible.
The subcell matrix in this work is given by 
\bea\label{eq:Mcp}
	\Mcp = \Sum_{c \in \mathcal{C}(p)}  z_{cp} \, \lcp \, \vec{n}_{cp} \otimes \vec{n}_{cp}, 
\eea
where we remind that $\lcp$ is the surface of the face $f$ of the three neighbor cells of $c$ sharing point $p$. $\vec{n}_{cp}$ is its outward unit normal and $z_{cp}=z_c$ is an approximation of the swept mass flux.
%\raphael{@Walter are you using $z_c$ or Dukowicz approximation involving $v_p, v_c$ in a NL way?}
%given following Dukowicz \cite{Duk85} as 
%\bea \label{eq:zfcp}
%	z_{cpf} &=& \rho_c \left[a_c + \Gamma_c  | (\vec{v}_c-\vec{v}_p).\vec{n}_{cpf} |\right]
%\eea
Once the velocity is determined thanks to (\ref{eq:nodal_solver}) then the trajectory equation can be invoked to compute the new point position.

%---------------------------------------------%
\subsection{Space-Time discretization --- ADER methodology} \label{ssec:time_discr}
%---------------------------------------------%
The time interval $[0,T] $  is discretized into time-steps  such that $t \in [t^n,t^{n+1}]$,
\begin{equation}
  t = t^n + \alpha \dt, \qquad  \alpha \in [0,1], 
\label{eqn:time}
\end{equation} 
where $t^n$ and $\dt$ represent the current time and time-step respectively. 
For evaluating the magnitude of $\dt$ we use a classical CFL condition and a criterion to avoid a too large increase of cell volume in a single time-step \cite{Maire2007,Maire2009}. \\
The time discretization simply consists in evaluating (\ref{eq:tau_c}-\ref{eq:e_c}) from %the geometry and
the state vectors given at $t^* \in [t^n,t^{n+1}]$, that is 
\beqa 
\label{eq:tau_c^n}
  \tau_c^{n+1} &=& \tau_c^n + \Frac{\dt}{m_c} \Sum_{p\in \mathcal{P}(c)} \lcp^n \vec{n}_{cp}^n \cdot \vec{v}_p^*  , \\
\label{eq:v_c^n}
 \vec{v}_c^{n+1} &=& \vec{v}_c^{n} + \Frac{\dt}{m_c}\Sum_{p\in \mathcal{P}(c)} \vec{f}_{cp}^* , \\
\label{eq:e_c^n}
 e_c^{n+1} &=& e_c^{n} +\Frac{\dt}{m_c} \Sum_{p\in \mathcal{P}(c)} \vec{f}_{cp}^* \cdot \vec{v}_p^* , 
\eeqa
and the trajectory equation as
\bea
\vec{x}_p^{n+1} = \vec{x}_p^{n+1} + \dt \, \vec{v}_p^*,
\eea
where $\vec{v}_p^*$ is obtained from the nodal solver 
\bea\label{eq:nodal_solver^n}
\Mp \vec{v}_{p}^* =   \Sum_{c\in  \mathcal{C}(p)}\Mcp \vec{v}_c^* -  \lcp^n \tens{T}_c^* \vec{n}_{cp}^n ,
\eea
thanks to the discrete subcell and nodal matrices $\Mcp$, $\Mp$, 
\bea\label{eq:Mcp^n}
\Mcp = \Sum_{c \in \mathcal{C}(p)}  z_{cp}^* \, \lcp^n \, \vec{n}_{cp}^n \otimes \vec{n}_{cp}^n,  \qquad
\Mp =  \Sum_{c\in  \mathcal{C}(p)}\Mcp,
\eea
and the subcell force (\ref{eq:subforce})
\bea \label{eq:subforce^n}
\vec{f}_{cp}^* = \lcp^n \tens{T}_c^* \vec{n}_{cp}^n + \Mcp (\vec{v}_{p}^*  -\vec{v}_c^* ).
\eea
The first-order time discretization simply considers $t^*=t^n$ and the cell-centered values of the state vector $\Q_c^n=(\tau_c,\vec{v}_{c},e_c)^n$. To obtain second order of accuracy in space a piece-wise linear reconstruction of the numerical solution $\Q_c$ must be carried out, thus obtaining higher order polynomials $\w_h^n(\vec{x})$ \cite{LAM2018,Maire2009}. Second-order time stepping demands that $t^* = t^{n+1/2} = \Frac12 (t^n + t^{n+1})$, which corresponds to the use of a midpoint rule to perform the time integration.
Classically a predictor-corrector \cite{Caramana-Burton-Shashkov-Whalen-98} or a Generalized-Riemann-Problem (GRP) scheme \cite{Maire2009} are used for this matter.
Contrarily, in this work, the second-order time discretization relies on the concept of the ADER (Arbitrary high order schemes using DERivatives) methodology following \cite{LAM2018}. 

The ADER procedure aims at computing high order space-time polynomials $\q_h(\vec{x},t)$ starting from the spatial reconstructed solution $\w_h^n(\vec{x})$ and performing a local time evolution of the governing equations \eqref{eq:cl}, that is
\bea 
\begin{aligned}
& \int \limits_{t^n}^{t^{n+1}}\rho\mddt{\q} - \DIV{\f(\q)} = 0,  \qquad \q=(\tau,\vec{v},e), \qquad \f(\q)=(\vec{v},\tens{T},\tens{T} \vec{v}), \\
& \int \limits_{t^n}^{t^{n+1}}\mddt{\vec{x}} = \vec{v}.           
\end{aligned}
 \label{eqn.NCPDE2}
\eea
The trajectory equation is coupled with the evolution of the governing PDE, thus the above nonlinear system \eqref{eqn.NCPDE2} is solved iteratively up to convergence for both the numerical solution $\q_h$ and the local geometry configuration $\vec{x}_h$. The space-time polynomials $\q_h$ coincide by construction with the high order spatial polynomials $\w_h^n$ at time $t^n$, i.e. $\q_h(\vec{x},t^n)=\w_h^n$, and all the details for the computation of a second order ADER predictor can be found in \cite{LAM2018}. Once the predictor is available, the subcell forces and the node values in \eqref{eq:tau_c^n}-\eqref{eq:e_c^n} are simply fed with the high order extrapolated values of the predictor, hence for any variable it holds $\q^*(\vec{x})=\q_h(\vec{x},t^*)$ for any $\vec{x}$.

%with a space-time piecewise linear reconstruction referred to as $\tens{T}_h(\vec{x},t)$ and $\vec{v}_h(\vec{x},t)$ evaluated at point $\vec{x}_p$ at time $t^{n+1/2}$, that is:
%\bea\label{eq:nodal_solver2^n}
%	\Mp \vec{v}_{p} =   \Sum_{c\in  \mathcal{C}(p)}\Mcp \vec{v}_h\left(\vec{x}_p,t^{n+1/2}\right) -  \lcp^n\tens{T}_h\left(\vec{x}_p,t^{n+1/2}\right) \vec{n}_{pc}^n , \\
%	\label{eq:subcell2^n}
%	\vec{f}_{pc} = \lcp^n \tens{T}_h\left(\vec{x}_p,t^{n+1/2}\right) \vec{n}_{pc}^n + \Mcp \left(\vec{v}_{p}^n-\vec{v}_h\left(\vec{x}_p,t^{n+1/2}\right)\right).
%\eea
%More precisely we average the $t^n$ and $t^{n+1}$ values of these polynomials to obtain the values at $t^{n+1/2}$ as 
%\bea
%\tens{T}_h\left(\vec{x}_p,t^{n+1/2}\right) = \Frac{ \tens{T}_h\left(\vec{x}_p,t^{n}\right)+\tens{T}_h\left(\vec{x}_p,t^{n+1}\right)}{2},
%\quad
%\vec{v}_h\left(\vec{x}_p,t^{n+1/2}\right) = \Frac{ \vec{v}_h\left(\vec{x}_p,t^{n}\right)+\vec{v}_h\left(\vec{x}_p,t^{n+1}\right)}{2} ,
%\eea
%while the geometrical vector and subcell matrices are frozen at $t^n$, see \cite{Maire_elasto}.
%In the ADER phraseology these space-time polynomial reconstructions are referred to as the \textit{predictors}, see \cite{LAM2018} for more details. \\
%In the context of hyper-elasticity equation~(\ref{eq:B_c^n}) can be updated once the predictor $\vec{v}_h$ is known using the definition of tensors $\tens{L}_c$ as a function of $\vec{v}_h$ and knowing $\tens{B}_c^n$.

%\paragraph*{Time evolution of the left Cauchy-Green tensor.}
The governing PDE system includes also the constitutive law \eqref{eq:trcb}, which describes the time evolution of $\tens{B}$. A semi-discrete form writes
\bea
\tens{B}^{n+1} = \tens{B}^{n} + \dt \, \int \limits_{t^n}^{t^{n+1}} \left( \tens{L}\tens{B} + \tens{B}\tens{L}^t \right) \, dt, \qquad \tens{L} = \nabla \vec{v}.
\label{eqn.semidiscrB}
\eea
The first order scheme is simply given by the Euler method in time and no reconstruction in space, thus it reads
\bea
\tens{B}_c^{n+1} = \tens{B}_c^{n} + \dt \left( \tens{L}_c\tens{B}_c + \tens{B}_c\tens{L}^t_c \right)^n,
\label{eqn.firstB}
\eea
with the spatial discretization of $\tens{L}_c$ given by \eqref{eq:div_grad_v}. A second order update of $\tens{B}$ is obtained by applying a Crank-Nicolson method to solve the integral ODE \eqref{eqn.semidiscrB}, hence one has
\beqa
\tens{B}_c^{n+1} &=& \tens{B}_c^{n} + \frac{\dt}{2} \left[ \left( \tens{L}_c\tens{B}_c + \tens{B}_c\tens{L}^t_c \right)^n + \left( \tens{L}_c\tens{B}_c + \tens{B}_c\tens{L}^t_c \right)^{n+1} \right], \nonumber \\
\tens{B}_c^{n+1} - \left( \tens{L}_c\tens{B}_c + \tens{B}_c\tens{L}^t_c \right)^{n+1} &=& \tens{B}_c^{n} + \frac{\dt}{2} \left( \tens{L}_c\tens{B}_c + \tens{B}_c\tens{L}^t_c \right)^n.
\label{eqn.CN_B}
\eeqa
The knowledge of $\vec{v}^{n+1}$ is required for the computation of $\tens{L}_c^{n+1}$ in the left hand side of \eqref{eqn.CN_B}. The second order nodal solver \eqref{eq:nodal_solver^n} provides the velocity at time level $t^{n+\frac{1}{2}}$, while the velocity at the current time level $t^{n}$ is known. To obtain a compatible velocity at the new time level and therefore be able to compute $\tens{L}_c^{n+1}$, let consider the equivalence of the midpoint and the trapezoidal rule for solving the trajectory equation \eqref{eqn.trajODE} with second order of accuracy:
\begin{equation}
\left. \begin{array}{lll}
\vec{x}^{n+1} &=& \vec{x}^{n} + \dt \, \vec{v}^{n+\frac{1}{2}} \\
\vec{x}^{n+1} &=& \vec{x}^{n} + \frac{\dt}{2} \, \left( \vec{v}^{n+1} + \vec{v}^{n} \right)
\end{array} \right\} \qquad \Rightarrow \qquad \vec{v}^{n+1} = 2 \vec{v}^{n+\frac{1}{2}} - \vec{v}^{n}.
\end{equation}
Once $\tens{L}_c^{n+1}$ is evaluated, equation \eqref{eqn.CN_B} constitutes a linear system for the unknown $\tens{B}_c^{n+1}$ that can be analytically inverted and solved.

%---------------------------------------------%
\subsection{Limiting: \aposteriori MOOD loop} \label{ssec:MOOD}
%---------------------------------------------%
%
%\textbf{MOOD loop, cascade and parachute scheme.}
While in the original ADER schemes the limiting relies on \apriori limited WENO reconstructions for all variables \cite{DumbserKaeser06b,ADERNC},
here we adopt an \aposteriori MOOD paradigm, see \cite{CDL1,LAM2018}.
Indeed the MOOD method is based on an \emph{a posteriori} evaluation of
the numerical solution, that is at $t^{n+1}$, to determine if some dissipation is needed.\\
The technique is \aposteriori in the sense that we compute a solution
at time $t^{n+1}$, and, then, determine if this candidate solution is acceptable, or not.
The candidate solution is first computed with a second-order accurate unlimited scheme using a centered reconstruction stencil.
Then a detection procedure determines the problematic cells, i.e.
the cells where the approximation does not respect some user-given criteria.
For those cells the solution is locally recomputed with a lower-order but more robust scheme.
In this work we consider three schemes forming a \emph{cascade}, each of them chosen to comply with one specific objective:
\begin{enumerate}
\item \textit{Accuracy}   gained with the unlimited piece-wise-linear polynomial reconstruction: maximal second-order of accuracy, possibly oscillating;
\item \textit{Robustness} gained with the piece-wise-linear polynomial reconstruction supplemented with Barth-Jespersen (BJ) \cite{BarthJespersen} slope limiter: between first- and second-order of accuracy, essentially-non-oscillatory, may not be positivity-preserving; 
\item \textit{Fail-safe}  gained without any polynomial reconstruction: first-order of accuracy, positivity preserving, hyper-robust and dissipative.
\end{enumerate}
A cell which does not satisfy all detection criteria is recomputed with the next scheme in the cascade.
This procedure, called the MOOD loop, is repeated until each cell satisfies all detection criteria or if the latest scheme of the cascade is selected. 
In this case, the robust positivity preserving first order finite volume scheme is employed.
The role of this so-called \emph{parachute} scheme is to always produce a meaningful physical solution at the price of
an excessive numerical dissipation. Notice that in practice it is almost never used, and, the BJ slope limiter can be substituted by any other reasonalbe one.
The process of dropping in the cascade is called \emph{decrementing} and a numerical solution not yet valid is referred to as being a
\textit{candidate solution}. \\
%\textbf{Detection/Admissibility criteria.}
The efficiency of the \textit{a posteriori} MOOD paradigm is brought by the fact that usually few cells need decrementing.
As such the extra-work needed to recompute only few problematic cells is usually low.
In the present implementation, the MOOD loop simply embraces the main evolution routines of the
ADER method and iterates to recompute those cells with invalid values, detected by the admissibility criteria.
In the worst case scenario all cells in the domain are updated with the parachute scheme, leading to the true first-order
accurate and robust numerical solution.
On the other hand, in the best case scenario, all cells are admissible at first MOOD iterate, that is
with the first scheme of the cascade leading to a truly second-order accurate numerical solution --- no limiting whatsoever.
In any other case, the MOOD loop always converges and produced an acceptable numerical solution,
assuming that the parachute scheme does so.\\
In the case of hyper-elasticity the detection/admissible criteria are based on the discrete version of $\mathcal{A}$,
see remark~\ref{rem:PAD}, that is, a candidate solution $\Q_h^{n+1}$ is physically admissible if it belongs to 
\begin{equation} 
\mathcal{A}_h^n = \left\{ \Q_c^n=(\tau_c^n,\vec{v}_c^n, e_c^n,\tens{B}_c^n) \text{ s.t. } 
\tau_c^n>0, \;
\varepsilon_c^n >0, \;
%=e_c^n-\frac12 (\vec{v}^n_c)^2 >0, \;
\theta_c^n>0,    \right\} .
\label{eq:PAD}
\end{equation}
Notice that we do not really use the entropy production in each cell, i.e see in (\ref{eqn:admissible_set}), 
%$\rho_c^n \theta_c^n \Frac{\eta^{n+1}-\eta_c^n}{\dt}\geq 0$,
because it produces excessive dissipative numerical solutions without any apparent gain. \\
Moreover to avoid spurious oscillations we also demand that the candidate solution fulfills a Relaxed Discrete Maximum Principle
(RDMP) that is
\bea \label{eq:RDMP}
-\delta_c^n + m_c^n \leq \rho_c^{n+1,*} \leq M_c^n + \delta_c^n, \quad \text{with} \quad
\left\{\begin{array}{l}
\delta_c^n = \max( \delta_0, \delta_1|M_c^n-m_c^n| ), \\
m_c^n = \min_{c \in \mathcal{V}_c} ( \rho_c^n ), \\
M_c^n = \max_{c \in \mathcal{V}_c} ( \rho_c^n ),
\end{array} \right.
\eea
$\mathcal{V}_c$ is the von Neumann neighborhood of cell $c$ used to reconstruct the piece-wise polynomials.
We fix $\delta_0=10^{-4}$ and $\delta_1=10^{-3}$.
Otherwise noticed only the density variable is tested for the RDMP.
Any cell which is not belonging to $\mathcal{A}_h^n$ or does not fulfill (\ref{eq:RDMP}) is declared troubled
and sent back to $t^n$ along with its neighbors for their re-computation using the next scheme in the cascade, see \cite{LAM2018}. \\
% --- Constitutive law
The \aposteriori detection and correction allows to monitor mathematical or model involution to ensure that the numerical errors remain at an acceptable level\footnote{Such a concern was raised in \cite{Bauer_Caramana_Loubere_06} in the  context of hydrodynamics solved by a staggered Lagrangian scheme where the cell volume can be computed either from the point coordinates or a PDE for the specific volume $\tau$. The difference between the two ``measures'' was monitored to assess the internal consistency of the scheme.}.
The fact that modern cell-centered Lagrangian schemes fulfill the GCL by construction is one kind of such involution.
For the hyper-elasticity model, the identity $\det \tens{B} = J^2 = \left(\Frac{\rho}{\rho^0}\right)^2$ should also be ensured.
For each cell, numerically, $\rho_c^{n+1}$ is not directly identified as: $\rho_c^{n+1} = \rho_c^{0} \sqrt{\det \tens{B}_c^{n+1}}$ but
deduced from the new point positions $\vec{x}_p^{n+1}$ which further yield the cell volume $V_c^{n+1}$ and the density as $\rho_c^{n+1}=\Frac{m_c}{V_c^{n+1}}$.
Therefore no process in the numerical scheme ensures such equality to hold.
We therefore monitor their difference as a goodness criteria as
\bea \label{eq:test_detB_rho}
\big| \sqrt{\det \tens{B}_c^{n+1}} - \Frac{\rho_c^{n+1}}{\rho^0_c} \big| < L_c^3,
\eea
where $L_c$ is a cell characteristics length, computed in this work as the smallest diameter of the in-spheres.

%---------------------------------------------%
\subsection{Time-step monitoring} \label{ssec:time_discr}
%---------------------------------------------%
The time-step is restricted by the classical CFL condition in our Lagrangian context \cite{Maire2009}
\bea \label{eq:CFL}
	\dt = \min \left( \dt_{\text{volume}}, \; \dt_{\text{acoustic}}, \; \dt_\text{increase} \right),
\eea
where we have used a criterion to avoid a too large increase of cell volume in a single time-step
\bea \label{eq:CFL2}
	\dt_{\text{vol.}} = C_v \min_c \left(\Frac{V_c^n}{V_c^{n+1}-V_c^n} \right), \,
	\dt_{\text{acoust.}} = C_{\text{CFL}} \min_c \left( \Frac{L_c}{a_c} \right), \,
        \dt_{\text{incr.}} = C_i ( t^{n}-t^{n-1} ),
\eea
where $L_c$, $a_c$ are a cell characteristics length and sound-speed respectively and $\left\{ C_v,C_{\text{CFL}},C_i \right\}\in [0,1]^3$.
The last constrain is designed to avoid a loo large increase of $\dt$.
Notice that the \aposteriori detection allows to ensure the positivity of the cell volume and the internal
energy provided that the parachute first-order scheme does.
As such the time-step control must be suited for the parachute scheme.
In our simulations we take $C_v=0.2$, $C_i=0.1$ and $C_{\text{CFL}}=0.25$ otherwise noticed.  \\
Notice that the \aposteriori MOOD loop may also be used to try to exceed the time-step restrictions (\ref{eq:CFL}) at the price of creating more troubled cells, for instance by setting $C_{\text{CFL}}$ closer to one.

%---------------------------------------------%
\subsection{Boundary condition treatments}
\label{sec:BCs}
%---------------------------------------------%

The Boundary Conditions (BCs) play a crucial role in the time evolution of the numerical solution.
In the context of an hyper-elasticity model solved by the Lagrangian numerical scheme we consider several types of BCs, such as free traction, restricted normal/tangential displacement and contact/symmetry plans.
These classical BCs are described in appendix~\ref{app:BCs} in the context of hyper-elastic materials, and all are applied through the nodal solver, differently from other face-based FV schemes. \\ 
To enlarge even further the ability of the code to handle complex situations, we have added the possibility
for a BCs to change its type during the simulation, for instance transitioning from free-traction to null normal velocity.
Generally such BC type evolution is driven by the nullification of a cost or distance function $\mathcal{D}$.
For instance an elastic material balistically flying, impacting onto a wall, spreading and 
 eventually detaching, demands such type of evolving BCs, see for instance the test case 'Rebound of a hollow bar' in section~\ref{ssec.BarRebound}.\\
The transition from BC type $A$ ($\BC_A$) to $B$ ($\BC_B$) can be imposed in two different ways:
\begin{itemize}
\item at a prescribed instant $t=t_{BC}$ the type of BCs changes, hence $\BC_A$ $\rightarrow$ $\BC_B$;
\item when the moving medium approaches a prescribed target located, $\vec{x}_T$, \textit{i.e} the distance function $\mathcal{D}= |\vec{x}_p -\vec{x}_T|<\epsilon_\mathcal{D}$, where $\epsilon_\mathcal{D}$ is a user-given threshold value, and, the velocity vector points in the direction of the target, then $\BC_A$ $\rightarrow$ $\BC_B$. Later, if the medium happens to detach from the target, then the distance function becomes again greater than the threshold value and the original BC is restored, that is $\BC_B$ $\rightarrow$ $\BC_A$.
\end{itemize}
Finally, from a practical point of view a hierarchy between the type of BCs must be imposed.
For instance when two faces sharing the same node must fulfill two different types of BCs,
then they must be applied in a hierarchical manner, taking into account the most important one first, possibly relaxing the fulfillment of the other ones.
Also at a material corner, a wall type BC must prevail compared to free traction BC, in such a way it avoids the boundary node to penetrate into the wall line/plane. Our hierarchy is as follows:
1- wall BC (restricted normal/tangential displacement),
2- space-dependent BC on velocity or pressure,
3- symmetry BC,
and 4- free-traction BC. \\
Although it seems at first glance to be ``only'' implementation issues, the treatment of BCs is of utmost importance for 3D mesh-moving numerical scheme like ours.

%---------------------------------------------%
%\clearpage \newpage
%---------------------------------------------%
\section{Implementation considerations}
\label{sec:computer_science}
%---------------------------------------------%

\subsection{Algorithm} \label{ssec:algorithm}
\begin{figure}[!htbp]
  \begin{center}    
    \includegraphics[width=0.9\textwidth]{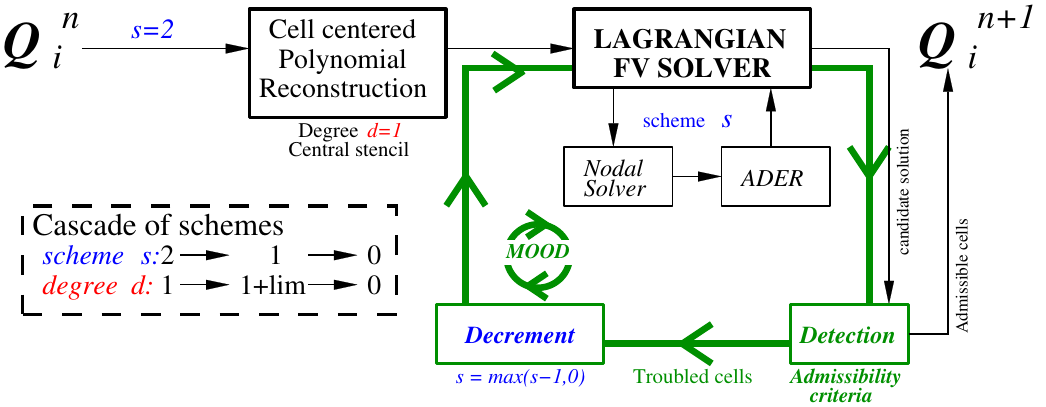}
    \caption{Sketch of the current Lagrangian numerical method and its MOOD loop.} 
    \label{fig:sketch}
  \end{center}
\end{figure}

In this section we recall the main steps of the MOOD loop applied to this cell-centered Lagrangian scheme sketched in figure~\ref{fig:sketch}.
First of all, cell centered unlimited polynomials of degree $d=1$ are reconstructed for any cell $i$ starting from data at $t^n$, $\Q_i^{n+1}$. Then a nodal solver and the ADER methodology allows to compute a candidate solution at $t^{n+1}$ with this 1st order accurate scheme labeled with $s=2$.
This candidate solution in cell $i$ can be either acceptable or numerically/physically wrong.
This is the purpose of the 'Detection' box to determine which cells are troubled, and, on the contrary to accept the admissible ones.
For those troubled cells, we pick the next scheme in the 'cascade' labeled $s=\max( s-1, 0)$, the scheme employs a piecewise-limited reconstruction (BJ limiter), or, no reconstruction at all, \textit{i.e}  the parachute scheme, in the latter, the first order Godunov scheme is used. 
Those troubled cells and their Voronoi neighbors are solely sent back for re-computation with this more robust scheme. This is the purpose of the 'Decrement' box.
This part of the solution which has been recomputed is re-tested against the detection criteria. New admissible cells are accepted, while troubled ones are again sent for re-computation with a more robust scheme.
Notice that this MOOD loop converges in a finite number of steps because the number of schemes in the cascade is fixed as well as the number of cells. \\
Once the slope limiter is chosen, the only parameters of the numerical method are the $\delta$'s ($\delta_0$ and $\delta_1$) in (\ref{eq:RDMP}) and the time-step control parameters (\ref{eq:CFL2}).
% ---- FIG ---------
% ---- FIG ---------

\subsection{Meshing and parallelization} \label{ssec:parallelisation}
The 3D Lagrangian simulation code is fully coded in Fortran and relies on MPI protocol
for the parallelization and the free graph partitioning software METIS \cite{metis}.
More precisely the computational domain $\Omega \in \mathbb{R}^3$ is first meshed with
a genuinely Coarse mesh made of large tetrahedra, say $N_C$, using any classical 3D mesher.
$N_C$ is chosen small enough for the resulting coarse mesh to be handled by one processor
without any difficulty.
This primary mesh is then partitioned among the total number of threads $N_\CPU$, see figure~\ref{fig:tets}-right for $N_\CPU=4$ in 2D and the coarse mesh in black.
Each MPI rank \emph{locally} refines its portion of the primary mesh by an arbitrary refinement factor $\ell>0$.
$N_\CPU$ and $\ell>0$ are given by the user.
A local structured recursive refinement is further applied.
The $\ell=0$th level corresponds to one of the primary tetrahedron, that is $N_\ell=1$ cell.
The $\ell=1$st level consists of its division into eight sub-tetrahedra, see remark~\ref{rem:split_tets}, to get $N_\ell=8$ sub-tetrahedra. 
The $\ell$th level consists of the division of all sub-tetrahedra obtained at level $\ell-1$, leading to $N_\ell=8\ell$ sub-tetrahedra.
In 2D the subdivision of one triangle is made into $4$ sub-triangles.
Each thread possesses only a portion of the full mesh and writes also its own output files.
As such the full mesh is never really assembled on a single thread leading to a reduction of memory storage. 
\begin{remark} \label{rem:split_tets}
To split one single tetrahedron we insert new vertices at the midpoints of each edge and connect the vertices together to form four new sub-tetrahedra associated to the vertices.
When removed from the parent tetrahedron, it leaves one central octahedron which can further be split into four more sub-tetrahedra by arbitrarily choosing an octahedron diagonal, see figure~\ref{fig:tets}-left.
% ---- FIG ---------
\begin{figure}[!htbp]
  \begin{center}
      \includegraphics[width=0.475\textwidth]{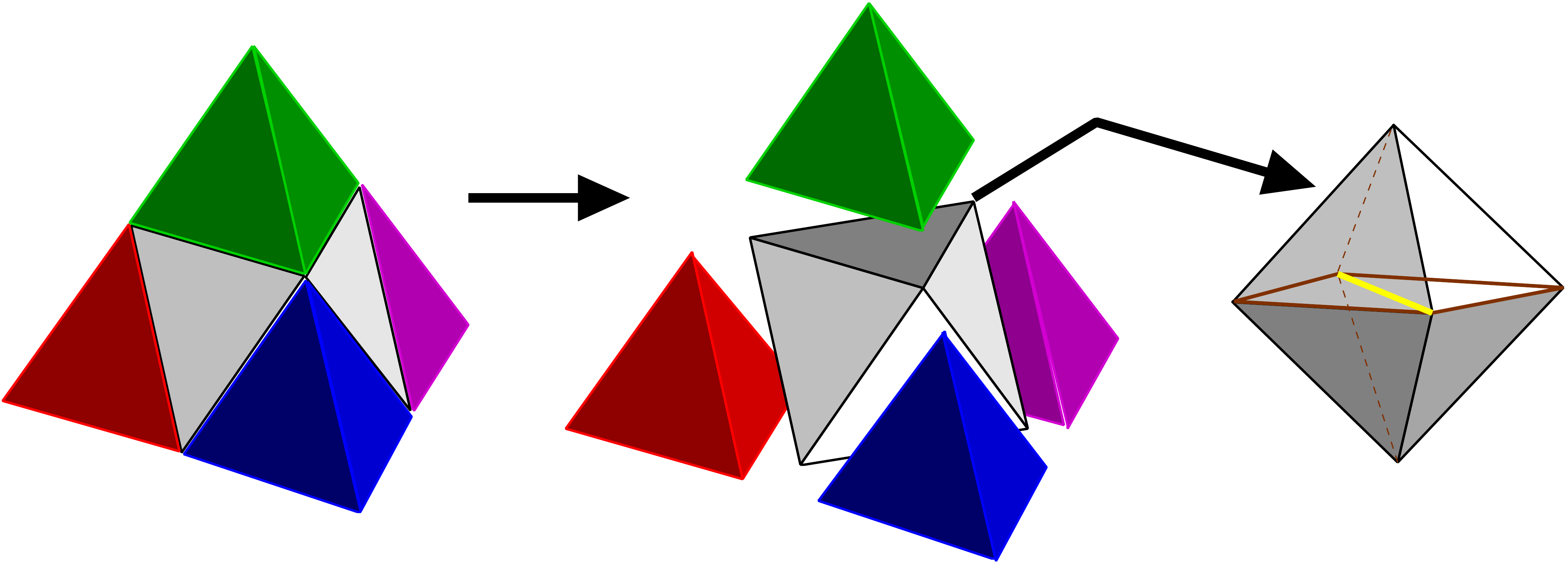}
      \includegraphics[width=0.5\textwidth]{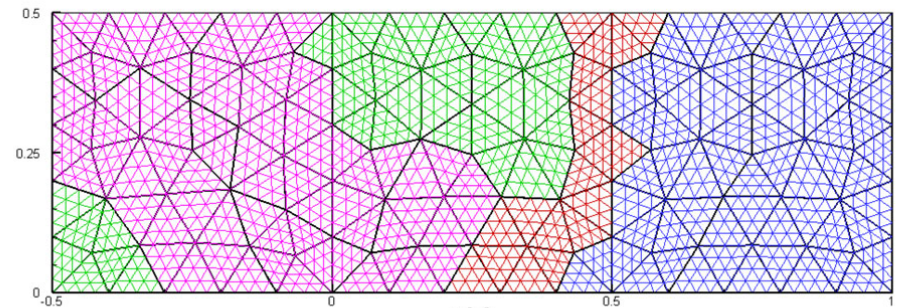}
      \caption{
        Left: Split of a tetrahedron into sub-tetrahedra by inserting six new midpoint edge vertices
        to get four corner sub-tetrahedra (colored ones). After choosing a diagonal (yellow line)
        to split the remaining central octahedron into fours more sub-tetrahedra, it yields a total
        of eight sub-tetrahedra ---
        Right: example of 2D partitioning on $N_\CPU=4$ threads (colors), the refinement is performed locally to each thread, only the coarse partition of large black triangles is actually built across the threads.} 
    \label{fig:tets}
  \end{center}
\end{figure}
% ---- FIG ---------
\end{remark}

%% file: Numerics.tex
\section{2D and 3D test problems}
\label{sec.validation}  
%---------------------------------------------%
In the following we present the results for a set of 2D and 3D benchmark test cases.
For each test problem the CFL stability coefficient is assumed to be $0.4$ in 2D and $0.25$ in 3D.
The time-dependent computational domain is addressed with $\Omega(t)$, while $\Q(\vec{x},t=0)\equiv\Q_0(\vec{x})=(1/\rho_0,\vec{v}_0,p_0,\tens{B}_0)$
denotes the vector of initial primitive variables typically used to setup the test problems.
$\tens{B}_0$ is set to the identity matrix as we only consider initially unloaded materials.
% Mesh construction
The unstructured tetrahedral meshes are obtained by meshing softwares, such as GMSH \cite{GMSH}
which takes a characteristics target length $h$ as input parameter.
%Refined meshes are obtained by reducing the objective characteristics length parameter
%by an appropriate factor, see section~\ref{ssec:parallelisation}.

In order to highlight the advantages of adding a second order limited scheme in the cascade compared to a first order discretization, according to \cite{Haider_2018}, the numerical dissipation $\boldsymbol{\delta}_h$ is monitored and here evaluated as
\begin{equation}
\boldsymbol{\delta}_h = \frac{\Psi+k-E_0}{E_0},
\label{eqn.numdiss}
\end{equation}
with the kinetic and total energy at the initial time $t_0$ defined by
\begin{equation*}
k_0 = \frac{1}{2}\vec{x}^2, \qquad E_0=\Psi_0 + k_0.
\end{equation*} 
Finally, if not stated otherwise, the simplified neo-Hookean equation of state \eqref{eq:EOS_neoHook} is adopted, while in the last test the stiffened gas EOS \eqref{eq:EOS_stiffened_gas} is used.

%
% TEST # 1 :  2D Swinging test
%
\subsection{2D swinging plate}
\label{ssec.swinging}
%\walter{setting by Scovazzi \cite{scovazzi3}}
The 2D swinging plate test problem, see \cite{Gil2D_2014,scovazzi3}, is employed to evaluate the
numerical order of convergence. The computational domain is $\Omega=[0,2]^2$
and the analytical solution for the displacement is given by
\bea
%\vec{u}^{ex}( \vec{x}, t ) = U_0 \sin ( \omega t)
%\left( \begin{array}{l}
%  -\sin \left( \Frac\pi2 x \right) \cos \left( \Frac\pi2 y \right) \\
%   \cos \left( \Frac\pi2 x \right) \sin \left( \Frac\pi2 y \right)
%\end{array} \right), \; \; \;
%U_0 = 10^{-2}~\text{m}.
\vec{v}^{ex}( \vec{x}, t ) = \omega U_0 \cos ( \omega t)
\left( \begin{array}{l}
  -\sin \left( \Frac\pi2 x \right) \cos \left( \Frac\pi2 y \right) \\
   \cos \left( \Frac\pi2 x \right) \sin \left( \Frac\pi2 y \right)
\end{array} \right), \qquad \omega=\frac{\pi}{2} \sqrt{\frac{2\mu}{\rho^0}},
\label{eqn.SwingIni}
\eea
with $U_0=5\cdot 10^{-4}~\text{m}$. The material under consideration is characterized by $\rho^0=1100~\text{kg}.\text{m}^{-3}$ with Young's modulus $E=1.7 \cdot 10^7~\text{Pa}$ and Poisson ratio $\nu=0.45$. 
%and the bulk modulus $\kappa$ by
%\bea
%\lambda=\Frac{R \nu}{(1+\nu)(1-2\nu)}, \qquad
%\mu = \Frac{E}{2(1+\nu)}, \qquad
%\kappa = \lambda + \Frac{2\mu}{3}.
%\eea
%For nearly incompressible material $\nu \to \frac12$.
%Here we consider $E=1.7 \times 10^6~\text{Pa}$ and $\nu=0.45$ \raphael{@Walter, could you please double check?}.
The velocity and displacement fields are divergence-free, leading to the exact pressure $p^{ex}=0$. Space-time dependent boundary conditions are prescribed for the normal velocities, according to the exact solution \eqref{eqn.SwingIni}.
Notice that the exact solution is a smooth one and the final time is set to $t_\text{final}=\pi/\omega$, so that
$\cos (\omega t_\text{final})= 1 $ and the final displacement corresponds to the initial one.
% Numerical results
In table~\ref{tab:convRates} we report the $L_2$ errors $\epsilon$ at the final time for the horizontal velocity $u$,
the first component of the left Cauchy-Green tensor $\tens{B}_{11}$ and of the Cauchy stress tensor $\tens{T}_{11}$.
The unstructured mesh made of triangles is successively refined and the final characteristics length $L_c(\Omega(t_\text{final}))$ is measured and further used to compute the numerical order of convergence $\mathcal{O}$ as reported in table~\ref{tab:convRates}, where one can notice that the numerical scheme is able to retrieve the second-order of convergence on this regular solution.
% ---- TAB ---------
\begin{table}[!htbp] 
  \begin{center}
    \numerikNine
    \begin{tabular}{|c|cc|cc|cc|}
      \hline
      $L_c(\Omega(t_\text{final}))$ & $\epsilon(u)$ & $\mathcal{O}(u)$ & $\epsilon(\tens{B}_{11})$ & $\mathcal{O}(\tens{B}_{11})$ & $\epsilon(\tens{T}_{11})$ & $\mathcal{O}(\tens{T}_{11})$ \\ 
      \hline
      \hline
      7.81E-02 & 2.144E-03 & ---  & 1.581E-04 & ---  & 9.681E+02 & --- \\
      5.21E-02 & 8.206E-04 & 2.37 & 7.072E-05 & 1.98 & 4.258E+02 & 2.03 \\
      3.91E-02 & 4.650E-04 & 1.97 & 3.914E-05 & 2.06 & 2.343E+02 & 2.08 \\
      3.13E-02 & 3.085E-04 & 1.84 & 2.473E-05 & 2.06 & 1.477E+02 & 2.07 \\ 
      2.60E-02 & 2.212E-04 & 1.82 & 1.699E-05 & 2.06 & 1.015E+02 & 2.06 \\
      \hline
      \multicolumn{2}{|c}{\text{Expected orders} $\rightarrow$}  & 2 & & 2 & & 2 \\
      \hline
    \end{tabular}
  \caption{
    \label{tab:convRates}
    Numerical errors in $L_2$ norm and convergence rates for the 2D swinging plate test computed with second order of accuracy Lagrange ADER scheme at time $t_\text{final}=\pi/\omega$. The error norms refer to the variables $u$ (horizontal velocity), $\tens{B}_{11}$ (first component of the left Cauchy-Green tensor $\tens{B}$) and $\tens{T}_{11}$ (first component of the Cauchy stress tensor $\tens{T}$).} 
  \end{center}
\end{table}
% ---- TAB ---------

%
% TEST # 2 :  Beryllium plate
%
\subsection{Elastic vibration of a Beryllium plate}  \label{ssec.BePlate}

\begin{figure}[!htbp]
  \begin{center}
    \includegraphics[width=0.8\textwidth]{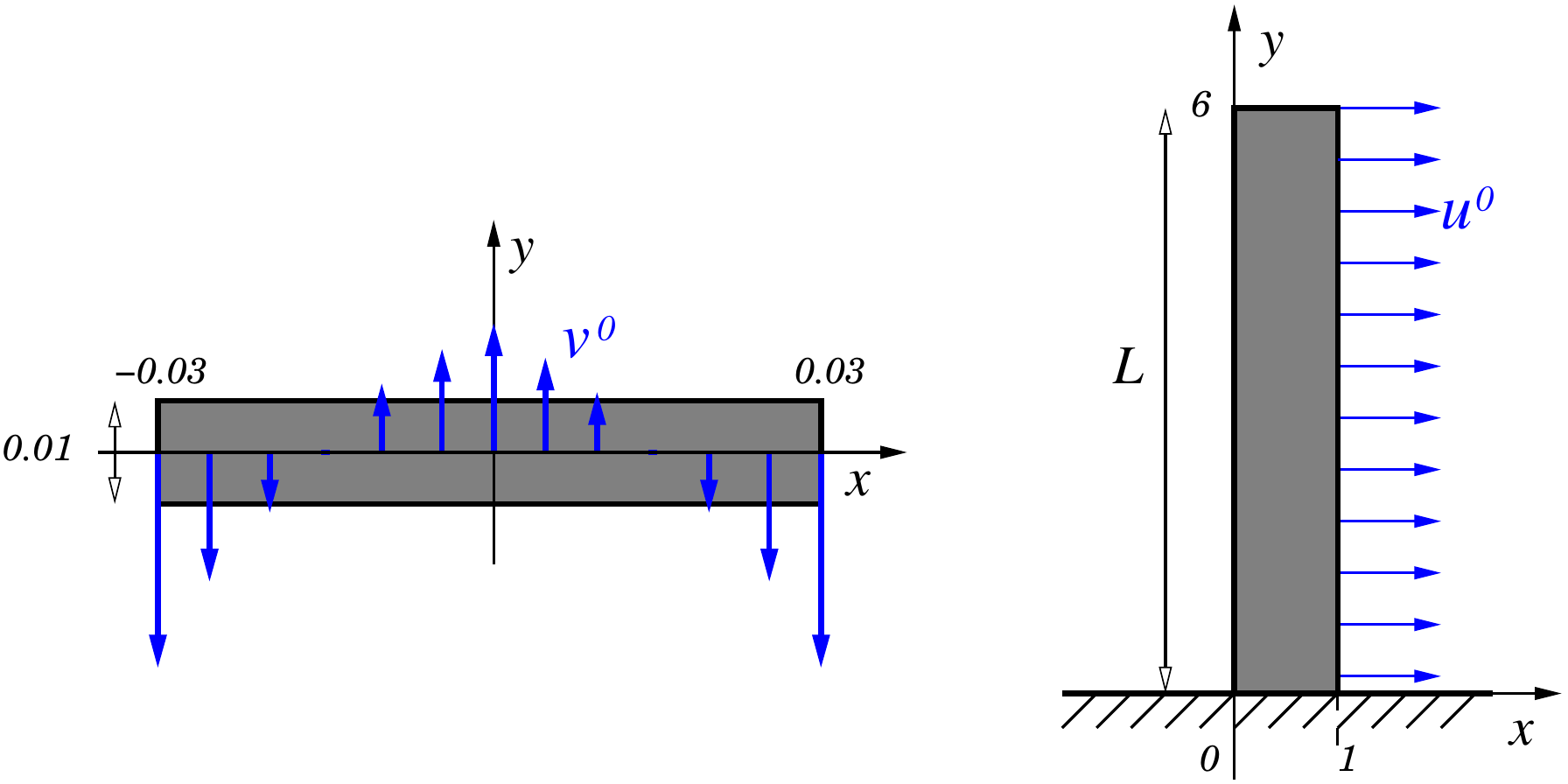}
    \caption{Sketch for the elastic vibration of a beryllium plate in section~\ref{ssec.BePlate} (left) and the finite deformation of a cantilever thick beam in section~\ref{ssec.BendCol} (right).} 
    \label{fig.sketch_Be_bar}
  \end{center}
\end{figure}

\begin{figure}[!htbp]
  \begin{center}
    \begin{tabular}{cc} 
      \includegraphics[width=0.47\textwidth]{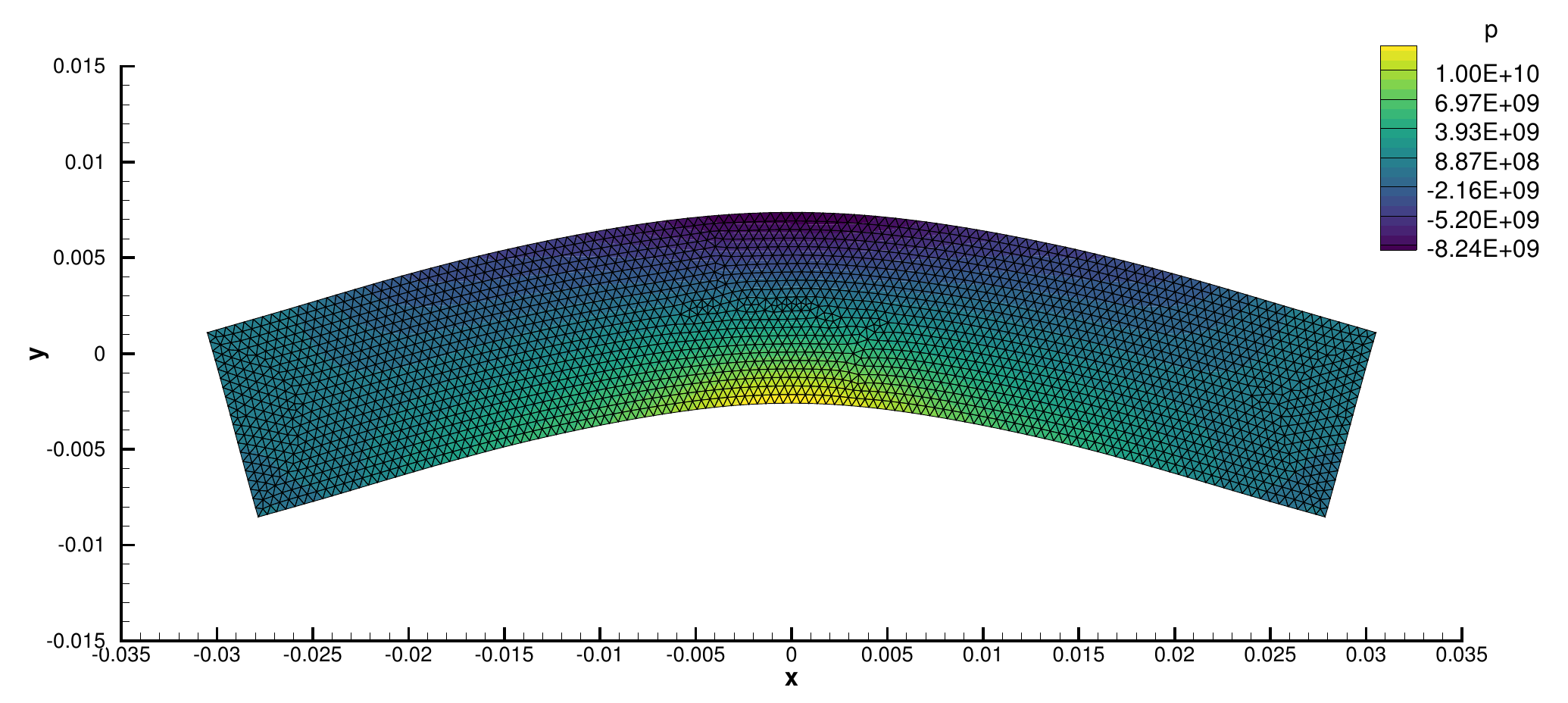} &       
      \includegraphics[width=0.47\textwidth]{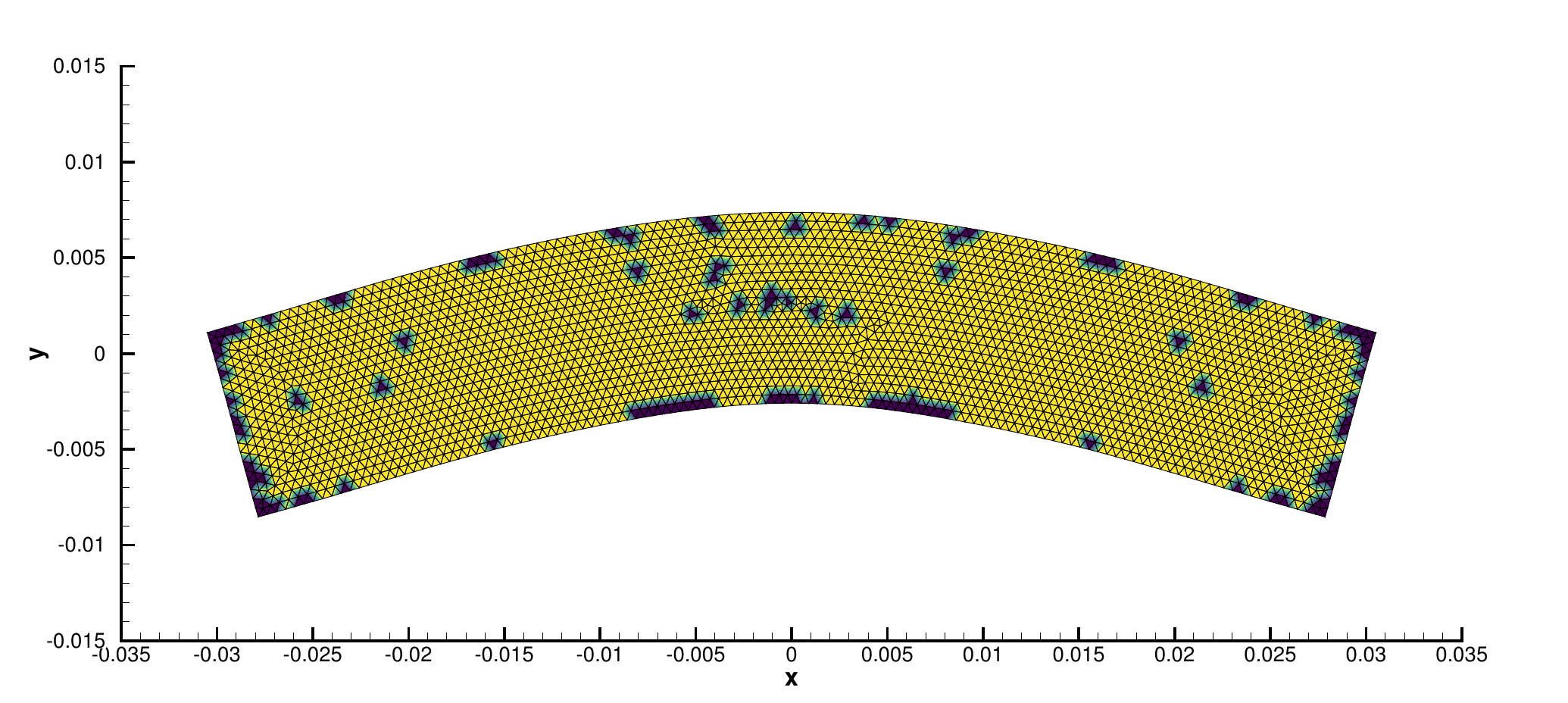} \\     
      \includegraphics[width=0.47\textwidth]{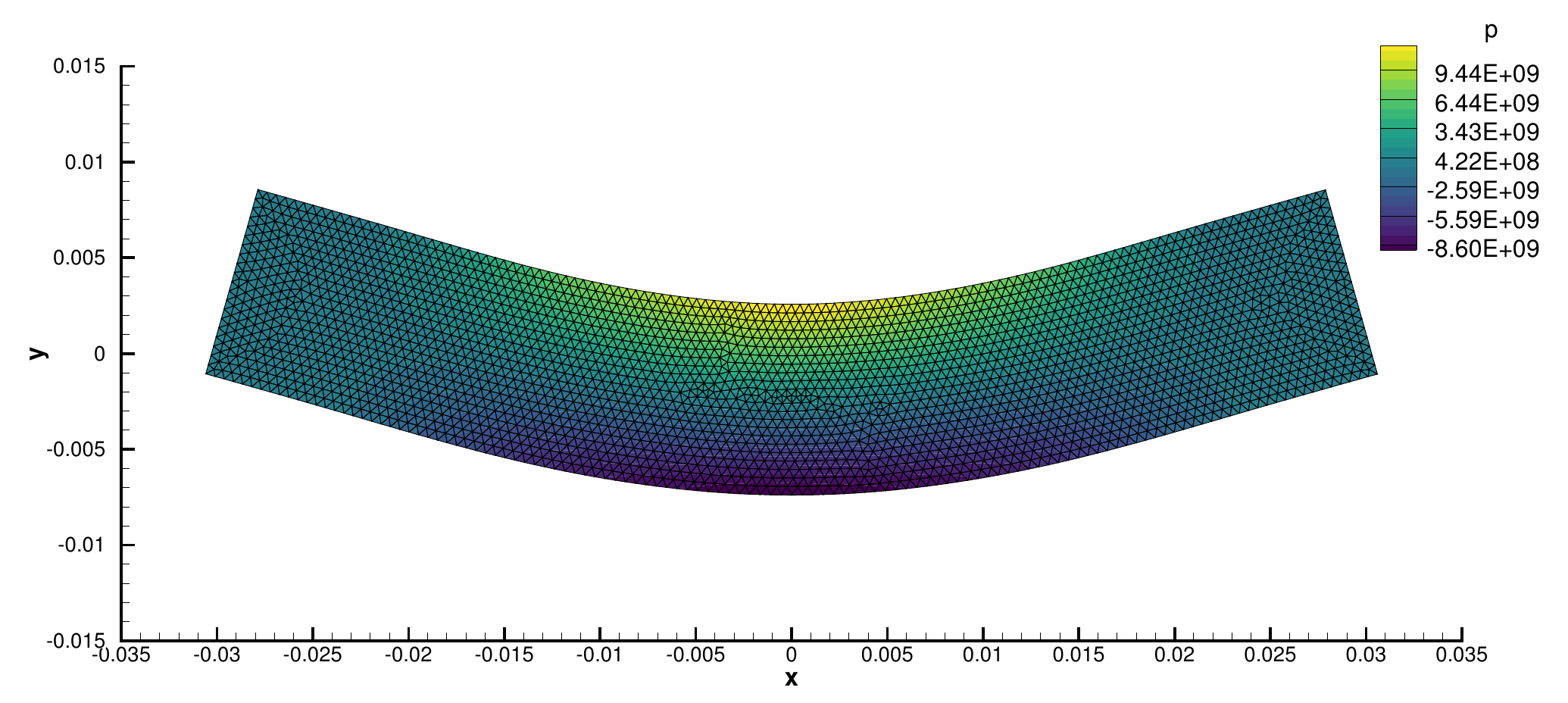} &
      \includegraphics[width=0.47\textwidth]{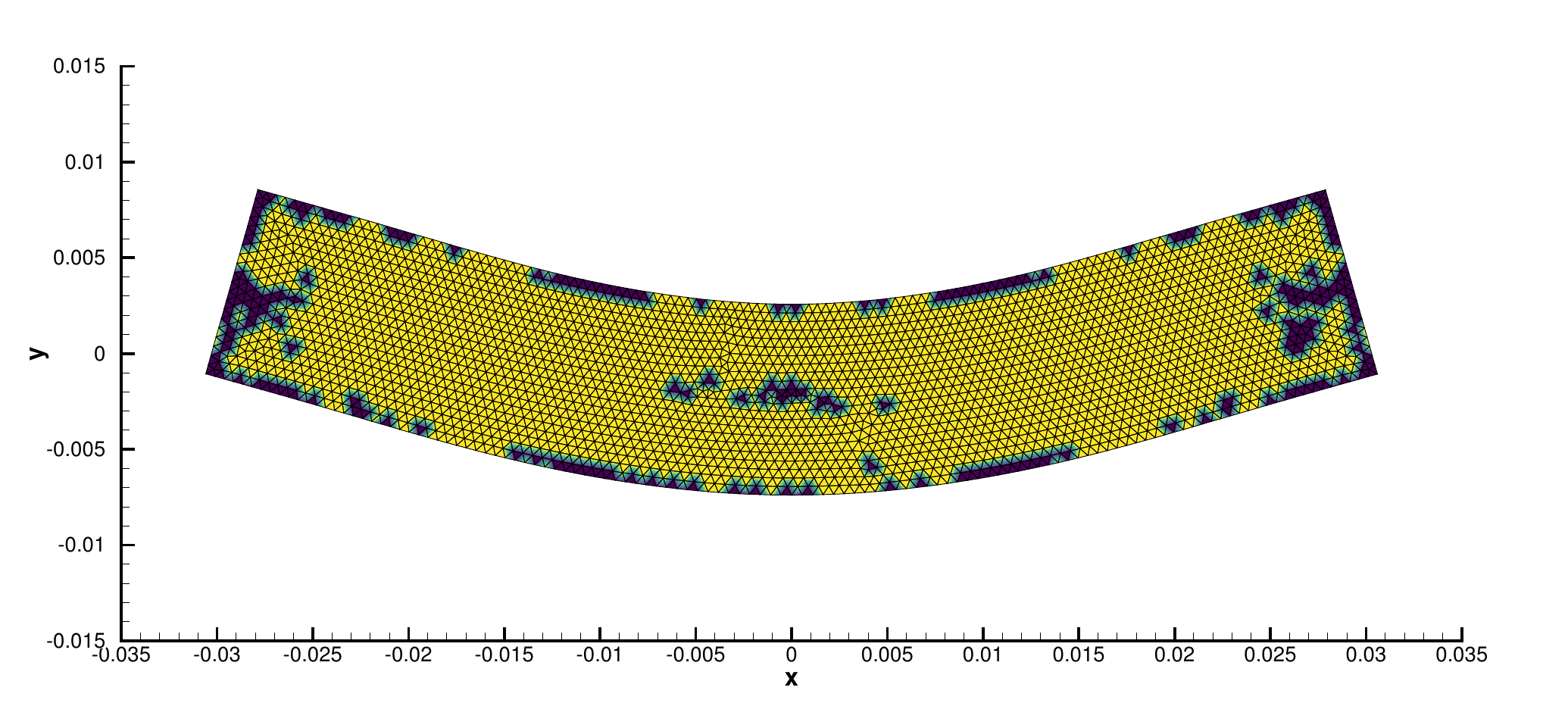} \\
      \includegraphics[width=0.47\textwidth]{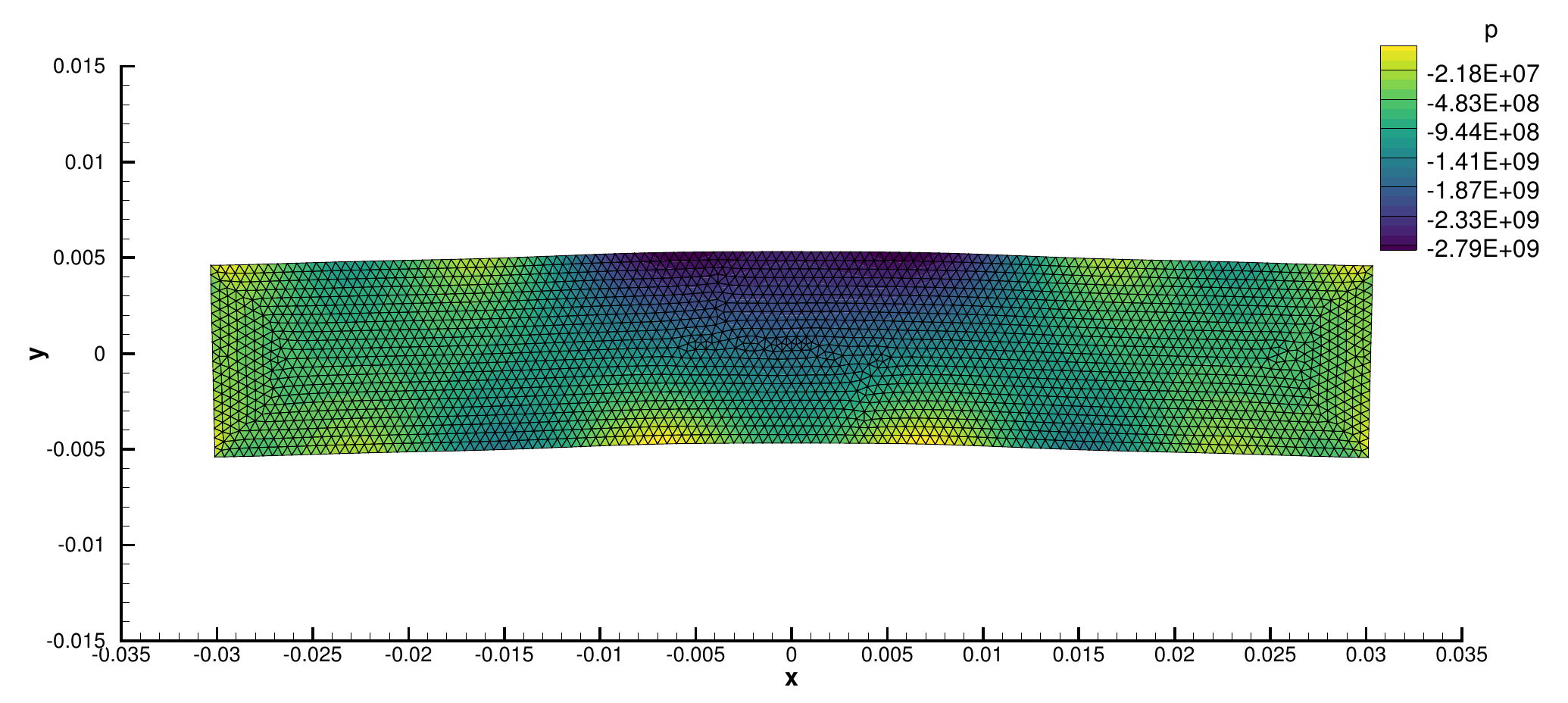}  &
      \includegraphics[width=0.47\textwidth]{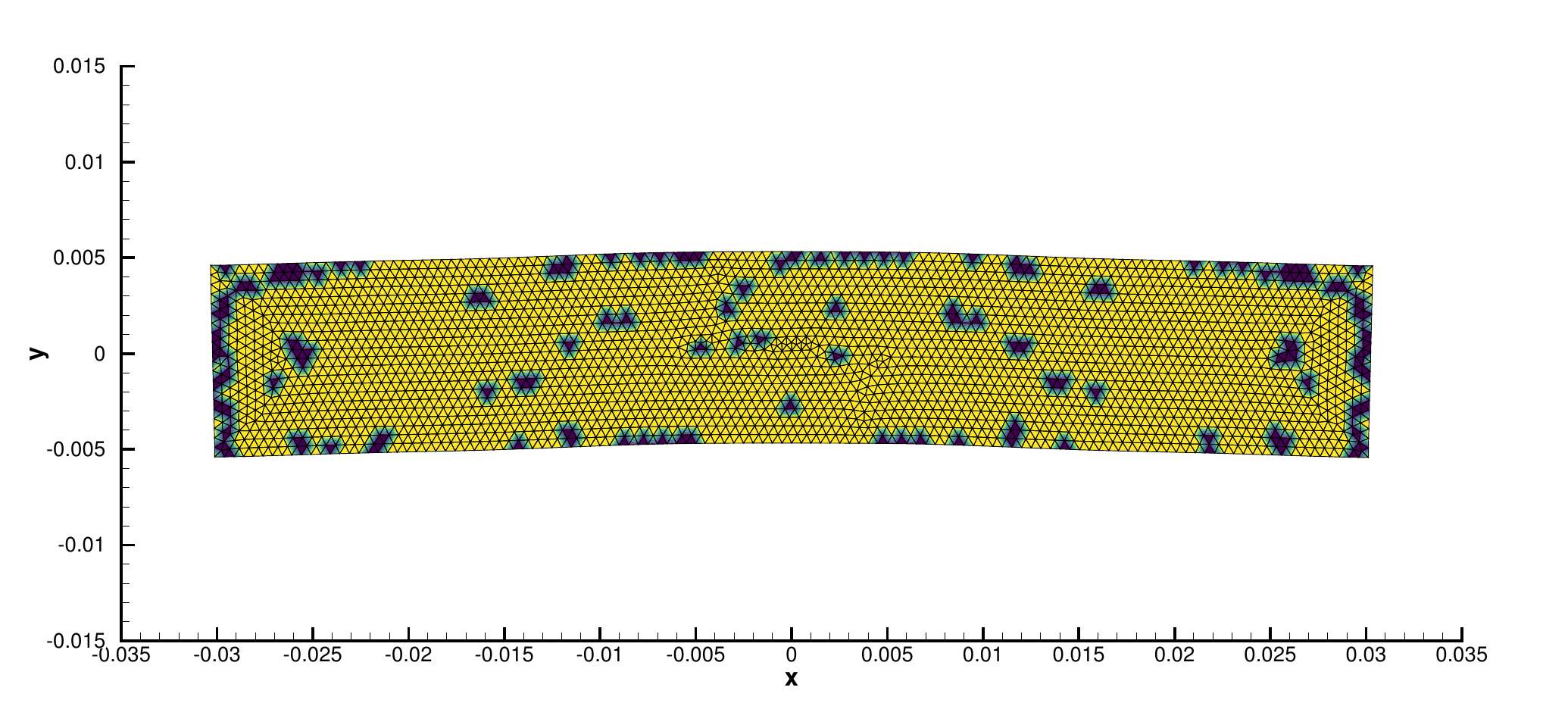}  \\
    \end{tabular} 
    \caption{
      Elastic vibration of a beryllium plate ---
      Numerical results at output times $t=10^{-5}$ (top), $t=2\cdot 10^{-5}$ (middle) and $t=3\cdot 10^{-5}$ (bottom) for pressure (left) and cell order map (right), the cells in yellow are at unlimited order 2, while the blue ones are the BJ limited ones. No first-order updated cell is observed.} 
    \label{fig.BePlate2D}
  \end{center}
\end{figure}

% Description
This test case describes the elastic vibration of a beryllium plate or bar, see \cite{Peshkov_Boscheri_Loub_hyper_hypo19,CCL2020} for instance.
Here we consider the 2D version, that is the vibration of a plate.
The computational domain is $\Omega(t=0)=[-0.03,0.03]\times[-0.005,0.005]$ of length $L=0.06~\text{m}$.
The material under investigation is characterized by:
$\rho^0=1845~\text{kg}.\text{m}^{-3}$, $E=3.1827 \cdot 10^{11}~\text{Pa}$ and $\nu=0.0539$. 
The initial material is loaded via a perturbed initial velocity field $\vec{v}^0=(0,v^0(x))$ of the form
\bea
v^0(x) = A \omega \left[  a_1(\sinh(x')+\sin(x')) - a_2(\cosh(x')+\cos(x')) \right],
\eea
where $x'=\alpha(x+L/2)$, $\alpha=78.834~\text{m}^{-1}$, $A=4.3369\times 10^{-5}~\text{m}$,
$\omega=2.3597\times 10^5~\text{s}^{-1}$, $a_1=56.6368$ and $a_2=57.6455$.
The final time is $t_\text{final}=3\cdot 10^{-5}~\text{s}$, see figure~\ref{fig.sketch_Be_bar}
for a sketch.
% ---- FIG ---------
% ---- FIG ---------
% BCs, mesh and runs
Free boundary conditions are applied on the plate faces. The unstructured triangulation is constituted of $N_c=5344$ cells.
% Results
In figure~\ref{fig.BePlate2D} we present the numerical results obtained at different output times for the pressure (left panels) and cell orders (right panels).
The pressure field is coherent with results from the litterature.
On the right panel we plot the cell order, which is equivalent to record which scheme from the cascade is actually employed. Yellow cells are dealt with the unlimited second order scheme (maximal order, prine to oscillation), while the blue ones employ a piecewise reconstruction limited by BJ slope limiter, via the \aposteriori MOOD loop. For this relatively mild problem, no cell is updated with the parachute scheme. Moreover no spurious modes nor artificial oscillations are observed.
% ---- FIG ---------
% ---- FIG ---------

\begin{figure}[!htbp]
  \begin{center}
    \begin{tabular}{cc} 
      \includegraphics[width=0.47\textwidth]{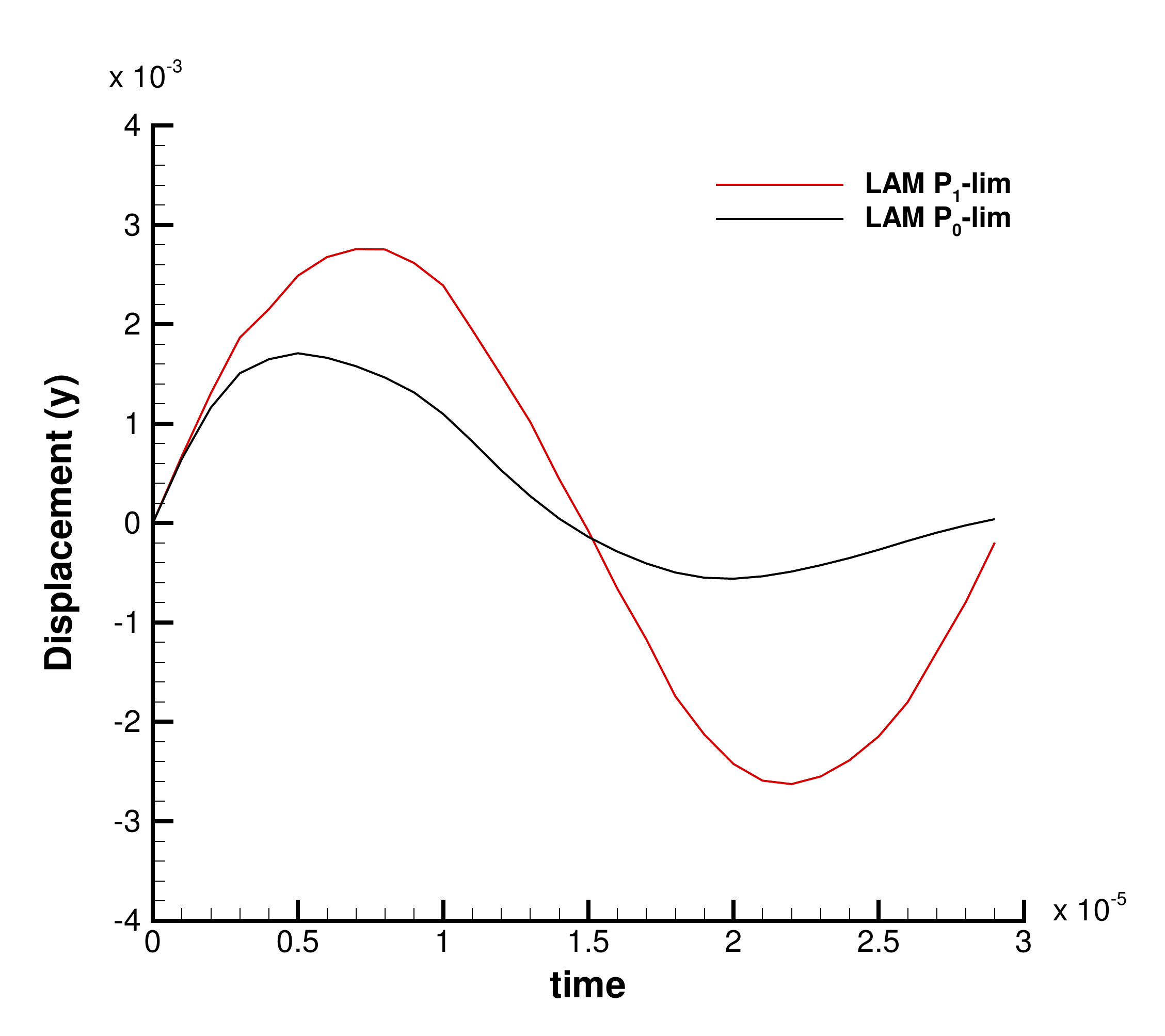} &       
      \includegraphics[width=0.47\textwidth]{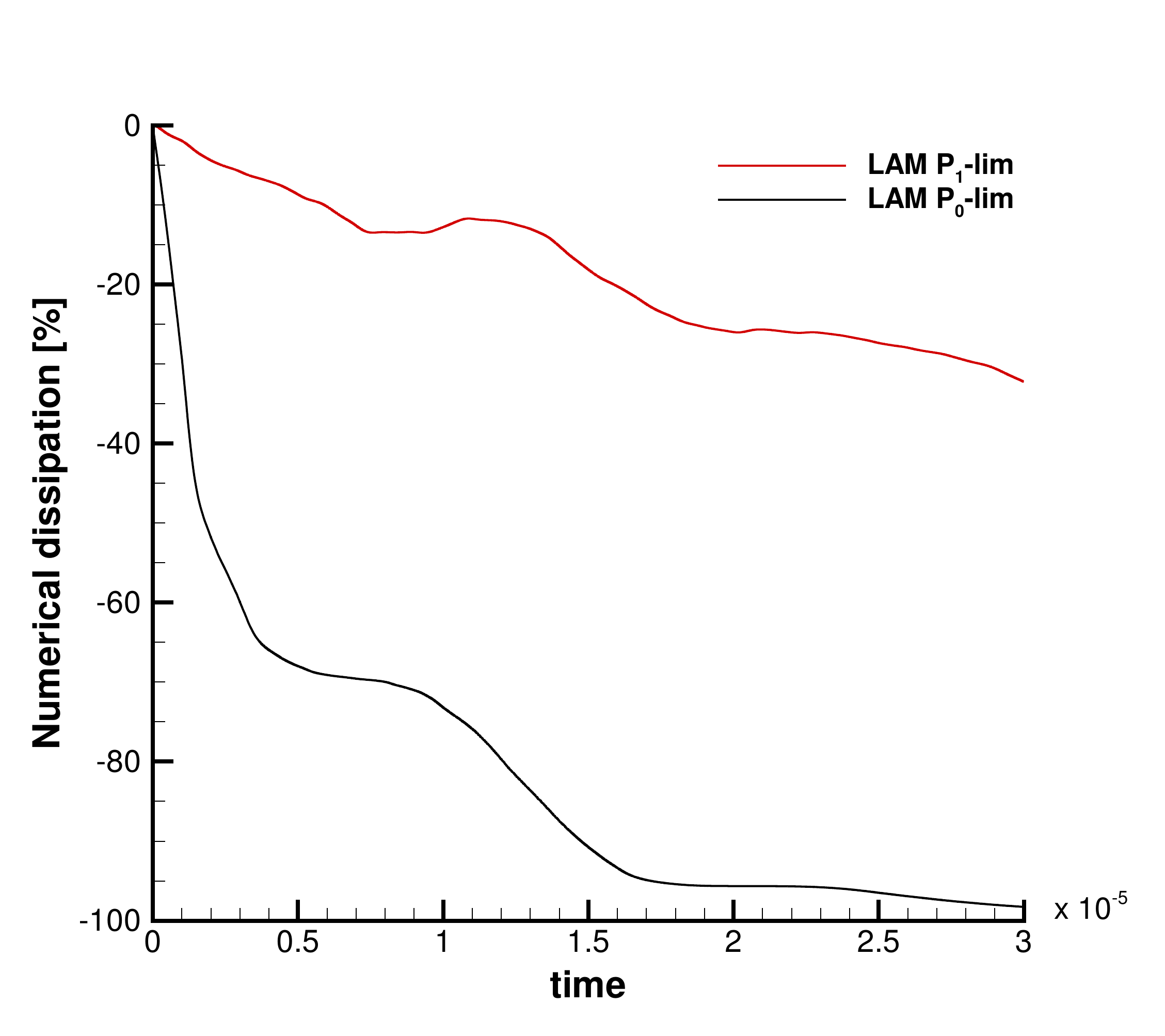} \\     
    \end{tabular} 
    \caption{Elastic vibration of a beryllium plate ---
      Comparison between the MOOD cascades:
      $\mathbb{P}_1 \rightarrow \mathbb{P}_0$ (LAM $\mathbb{P}_0$-lim) and
      $\mathbb{P}_1 \rightarrow \mathbb{P}_1^\text{lim} \rightarrow \mathbb{P}_0$ (LAM $\mathbb{P}_1$-lim) for the vertical displacement at the barycenter of the plate (left)
      and the computed numerical dissipation as a function of time (right).} 
    \label{fig.BePlate2D_comp}
  \end{center}
\end{figure}

In order to illustrate the reduction of dissipation when the cascade is not $\mathbb{P}_1 \rightarrow \mathbb{P}_0$, like in \cite{LAM2018}, but $\mathbb{P}_1 \rightarrow \mathbb{P}_1^\text{lim} \rightarrow \mathbb{P}_0$ instead, we show in figure~\ref{fig.BePlate2D_comp} two diagnostics.
First, on the left panel, the vertical displacement at the barycenter of the plate as a function of time is presented for the two cascades.
As expected the nominally second order scheme is able to follow the barycenter with lower dissipation.
On the right panel we enhance the actual numerical dissipation computed with \eqref{eqn.numdiss} which confirms that a high order scheme reduces the numerical viscosity by about $75\%$ at final time. 
% ---- FIG ---------
% ---- FIG ---------

%
% TEST # 3 :  Bending column
%
\subsection{Finite deformation of a cantilever thick beam}  \label{ssec.BendCol}
\begin{figure}[!htbp]
  \begin{center}
    \begin{tabular}{cc}
      \vspace{-0.5cm}
      \includegraphics[width=0.4\textwidth]{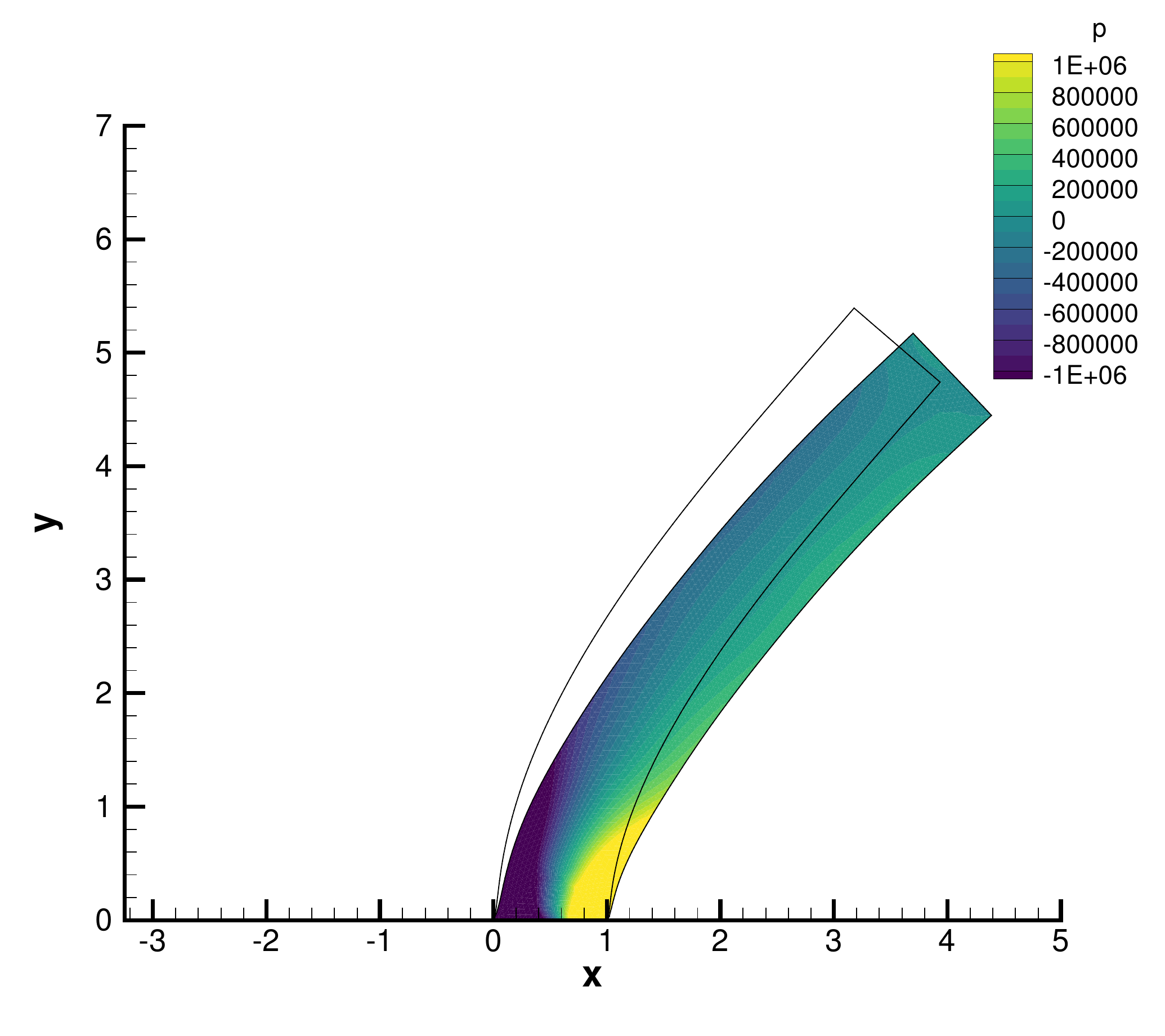} &
      \includegraphics[width=0.4\textwidth]{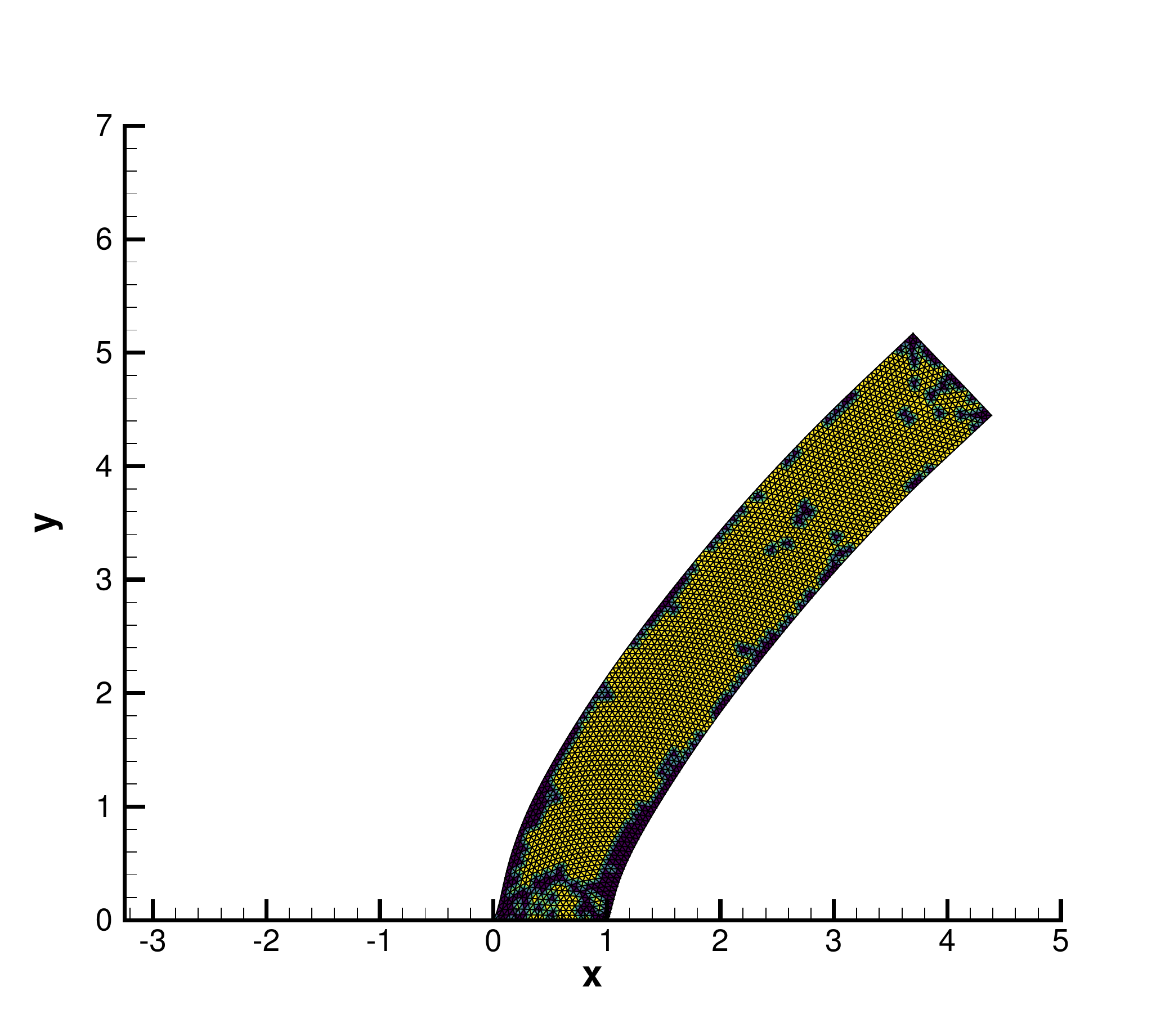} \\ 
      \vspace{-0.5cm}      
      \includegraphics[width=0.4\textwidth]{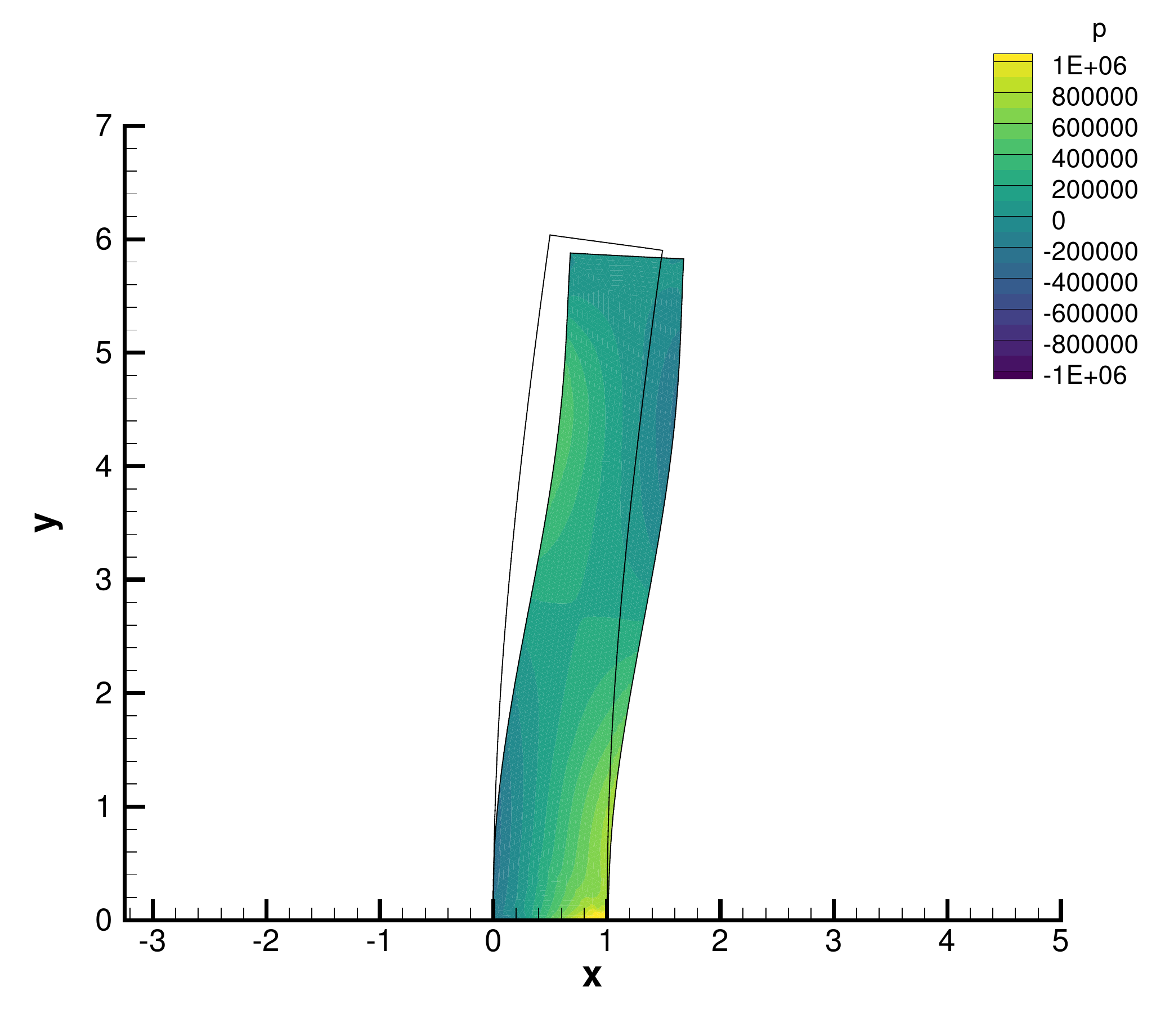} &
      \includegraphics[width=0.4\textwidth]{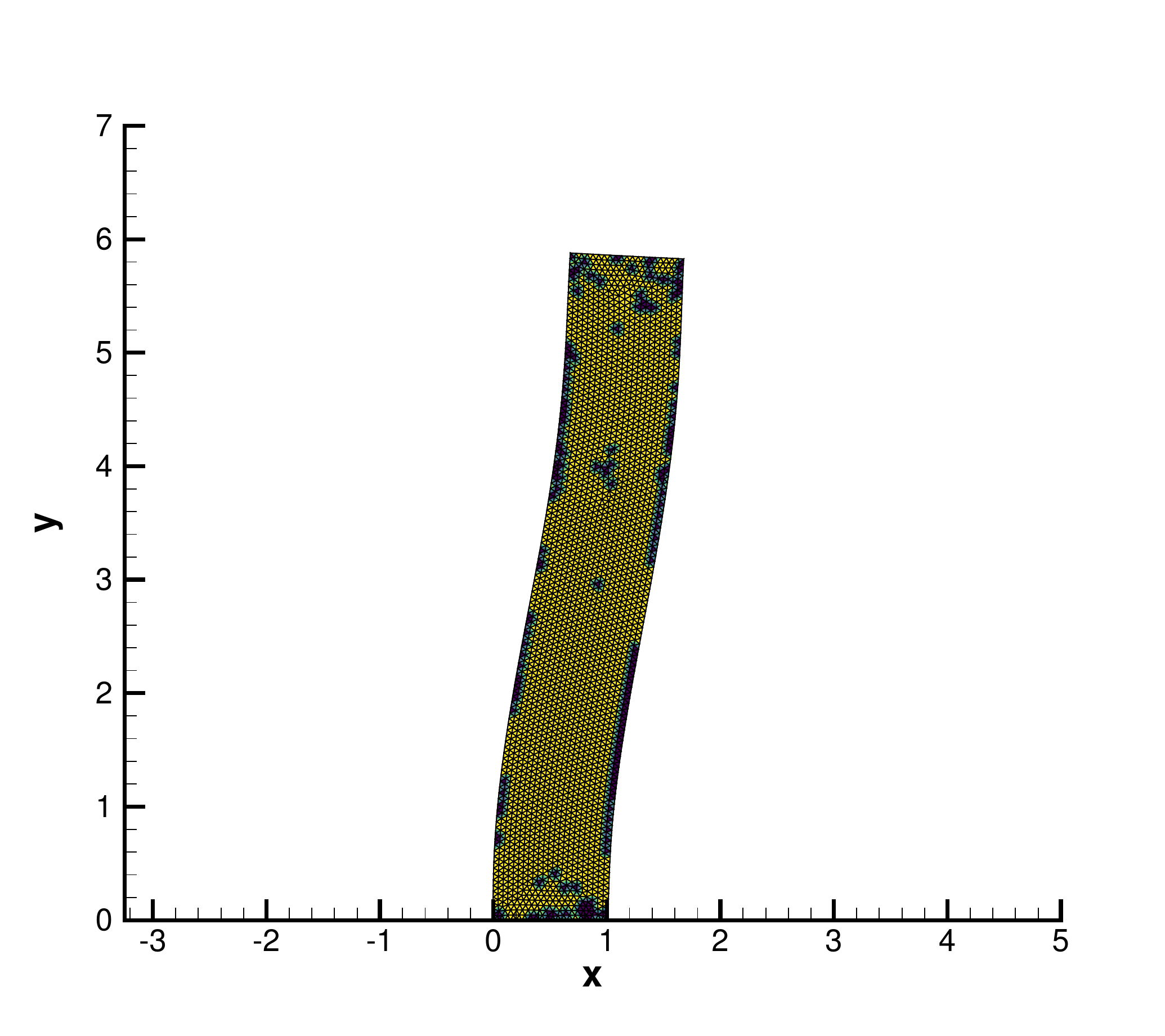} \\ 
      \vspace{-0.5cm}    
      \includegraphics[width=0.4\textwidth]{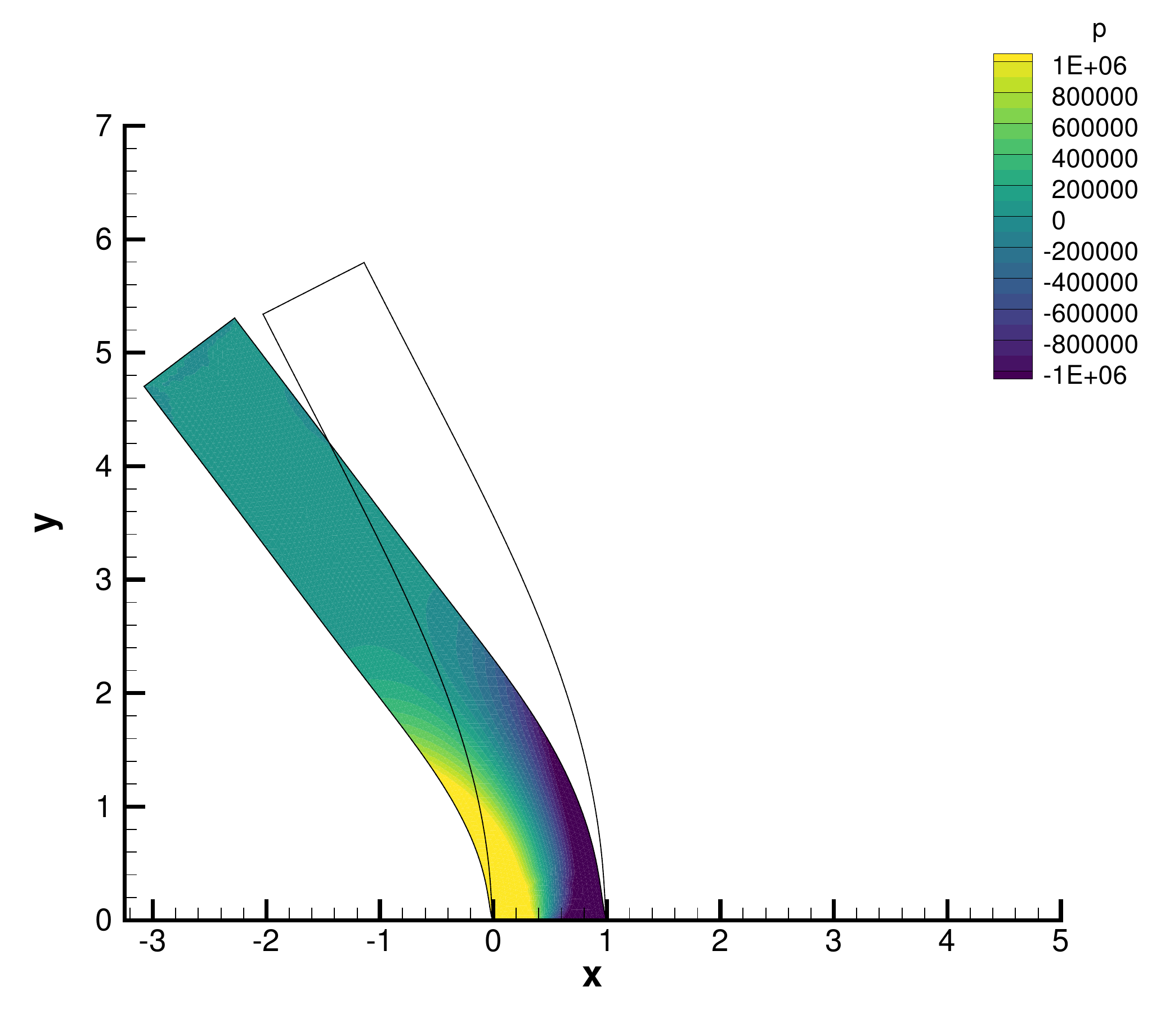} &
      \includegraphics[width=0.4\textwidth]{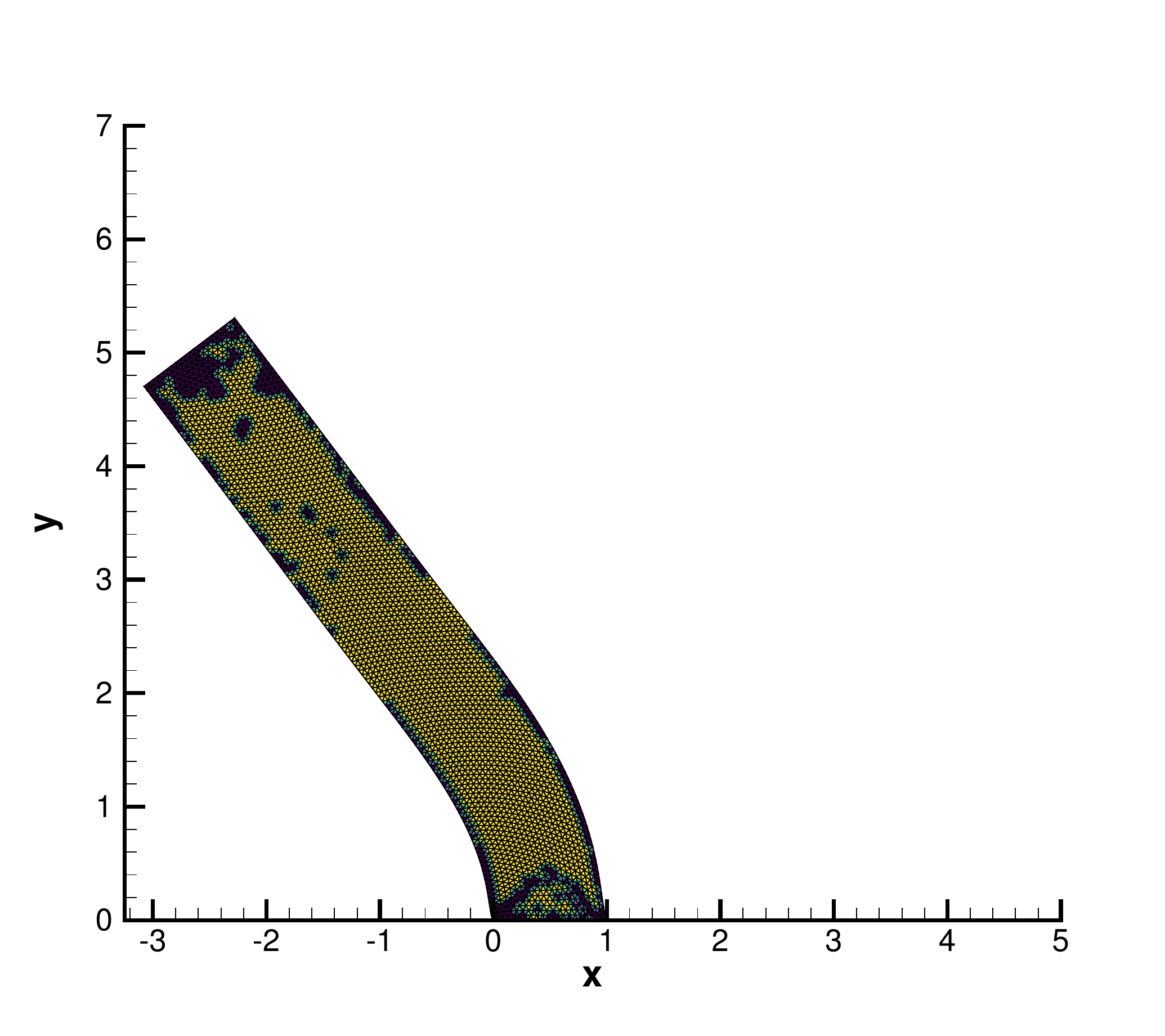}\\	
      \vspace{-0.5cm}	
      \includegraphics[width=0.4\textwidth]{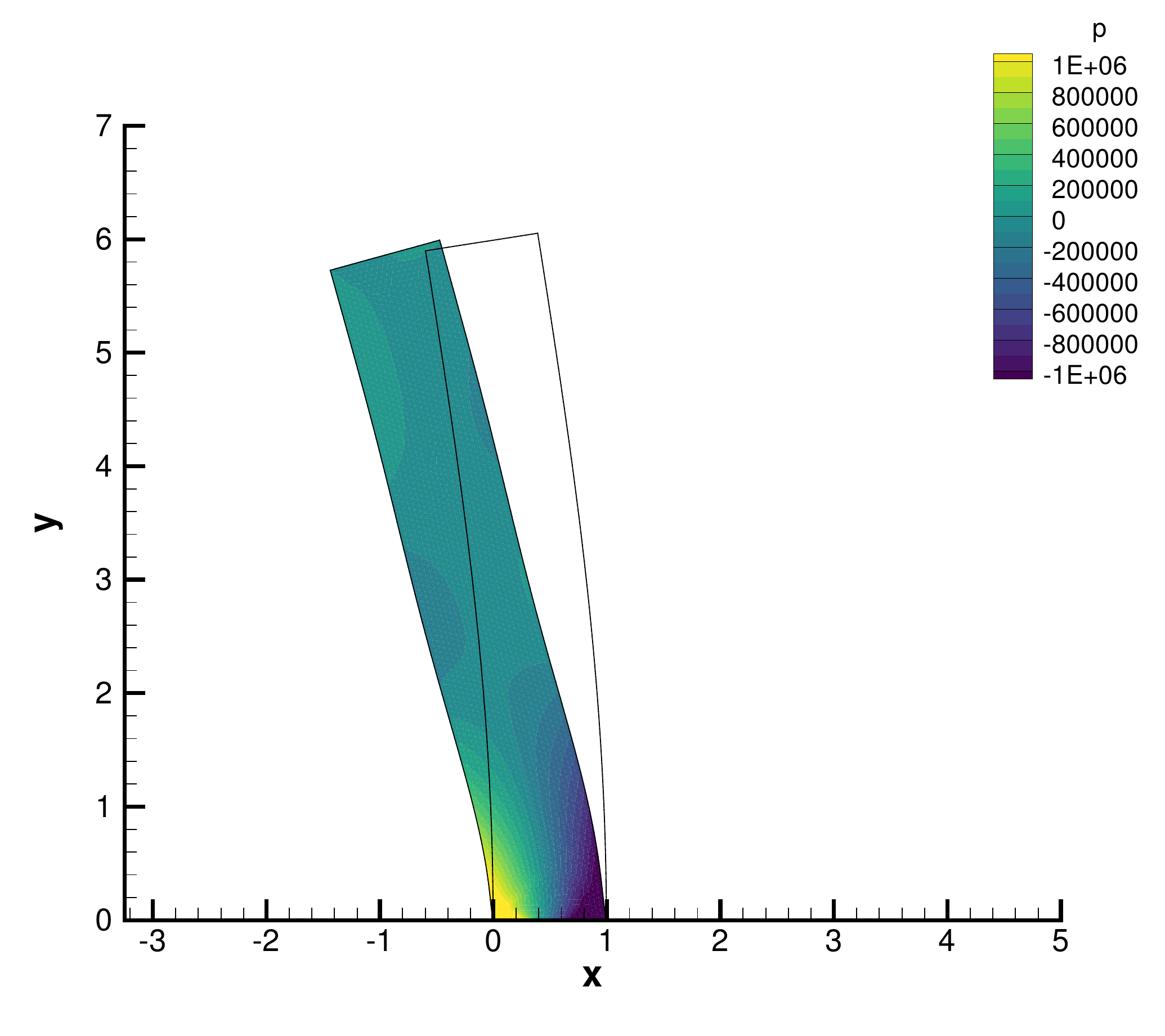} &
      \includegraphics[width=0.4\textwidth]{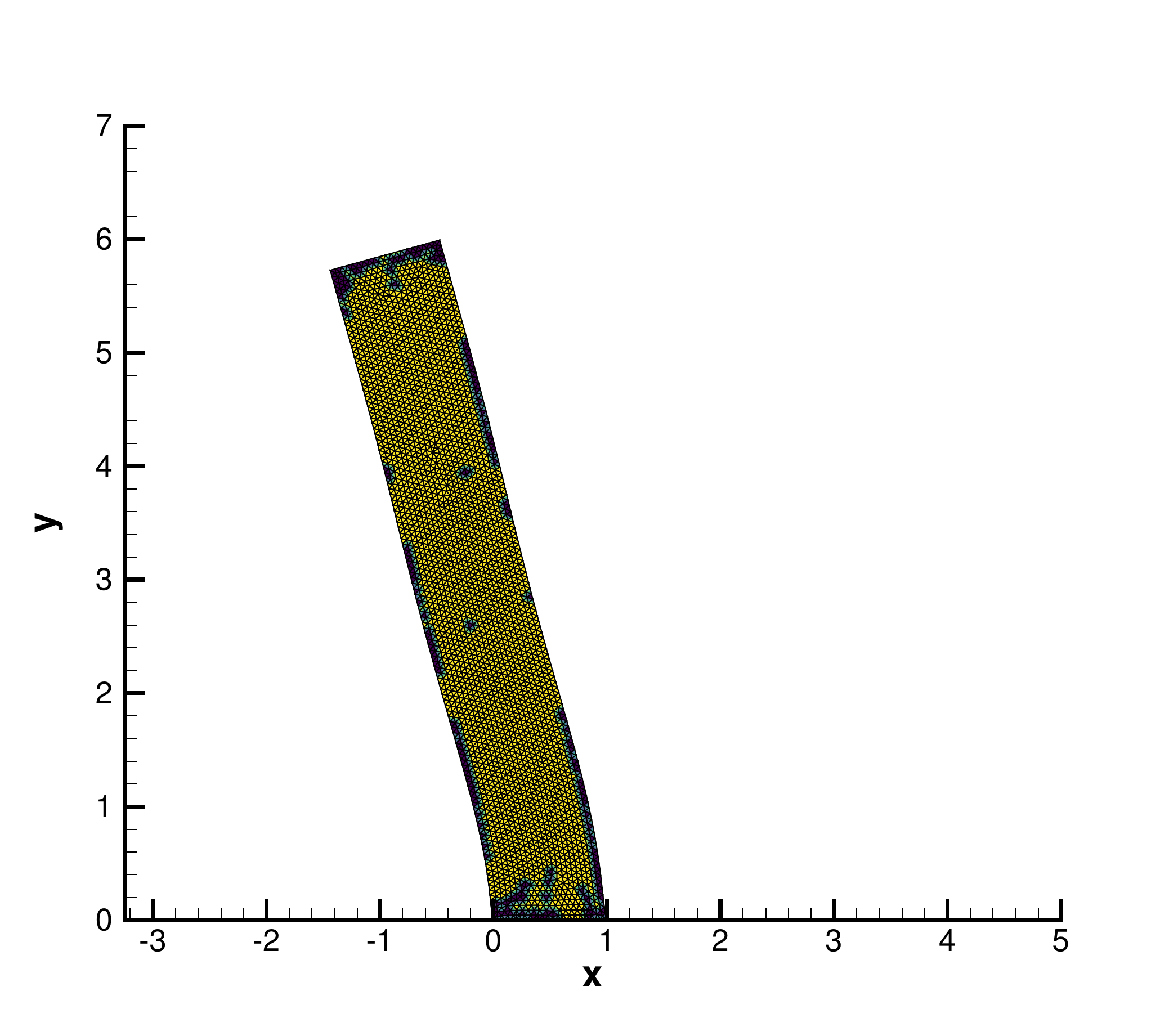}
    \end{tabular} 
    \caption{
      Cantilever thick beam test case ---
      Pressure distribution with deformed shapes (left column) and cell order map (right column)
      with the second-order \aposteriori limited Lagrangian scheme at output times
      $t=0.375$, $t=0.75$, $t=1.125$ and $t=1.5$ (from top to bottom row) ---
      Comparison of the deformed shape computed using the simpler cascade $\mathbb{P}_1 \rightarrow \mathbb{P}_0$ in black line on the left panels only.} 
    \label{fig.BendCol2D}
  \end{center}
\end{figure}

% Description
In \cite{Gil2D_2014} the authors present a test case involving
a finite deformation of a 2D cantilever vertical thick beam of length $L$ having a unit square cross section
and initially loaded by a uniform horizontal velocity $u^0 = 10~\text{m}.\text{s}^{-1}$
whilst the unit width base is maintained fixed, see figure~\ref{fig.sketch_Be_bar}
for a sketch.
We consider the initial computational domain $\Omega(t=0)=[0;1]\times[0;6]$
leading to $L=6~\text{m}$ and material characteristics $\rho^0=1100~\text{kg}.\text{m}^{-3}$, $E=1.7 \cdot 10^{7}~\text{Pa}$ and $\nu=0.45$.
% BCs, resh, run
Free boundary conditions are considered apart from the fixed-wall bottom part of the bar.
The mesh is made of $N_c=5442$ triangles.
% Results
The simulations are run with %the nominally second order Lagrangian scheme with an \aposteriori limiter using
the cascade $\mathbb{P}_1 \rightarrow \mathbb{P}_1^{\text{lim}} \rightarrow \mathbb{P}_0$.
On the left panels of figure~\ref{fig.BendCol2D}, we present the pressure distribution along with the deformed shapes at four different output times.
The results are qualitatively in adequation with the published ones from the litterature.
Moreover we observe on the right panels that the yellow cells (unlimited second-order scheme) are massively represented, while only few demand dissipation (blue cells).
For comparison purposes we also superimpose in black line the shapes obtained with the simpler cascade $\mathbb{P}_1 \rightarrow \mathbb{P}_0$ from \cite{LAM2018}.
As can be observed, this latter scheme is genuinely more dissipative, and, it numerically justifies the need for using a second order limited reconstruction within the cascade.
% ---- FIG ---------
% ---- FIG ---------

\begin{figure}[!htbp]
  \begin{center}
    \begin{tabular}{cc} 
      \includegraphics[width=0.47\textwidth]{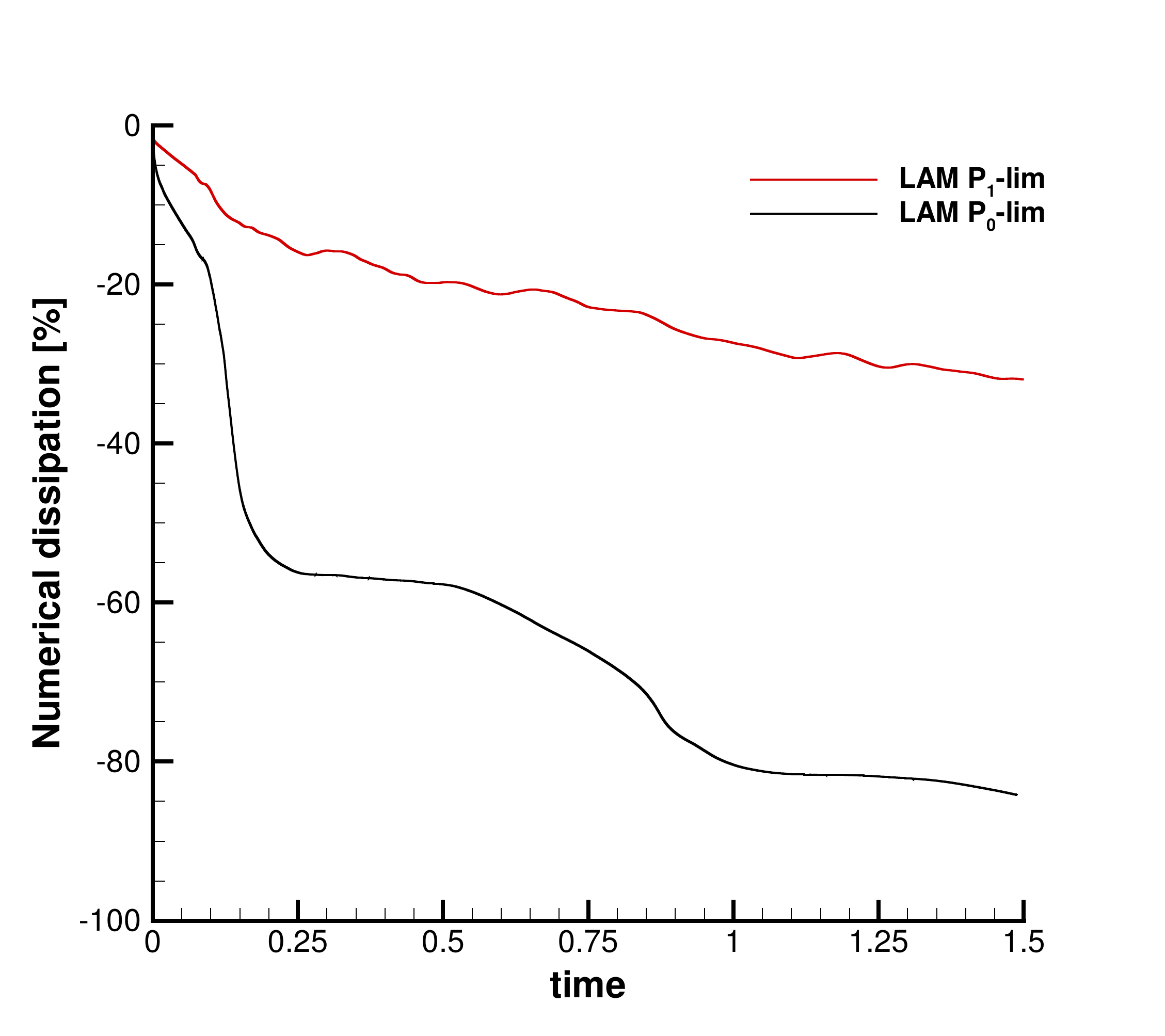} &       
      \includegraphics[width=0.47\textwidth]{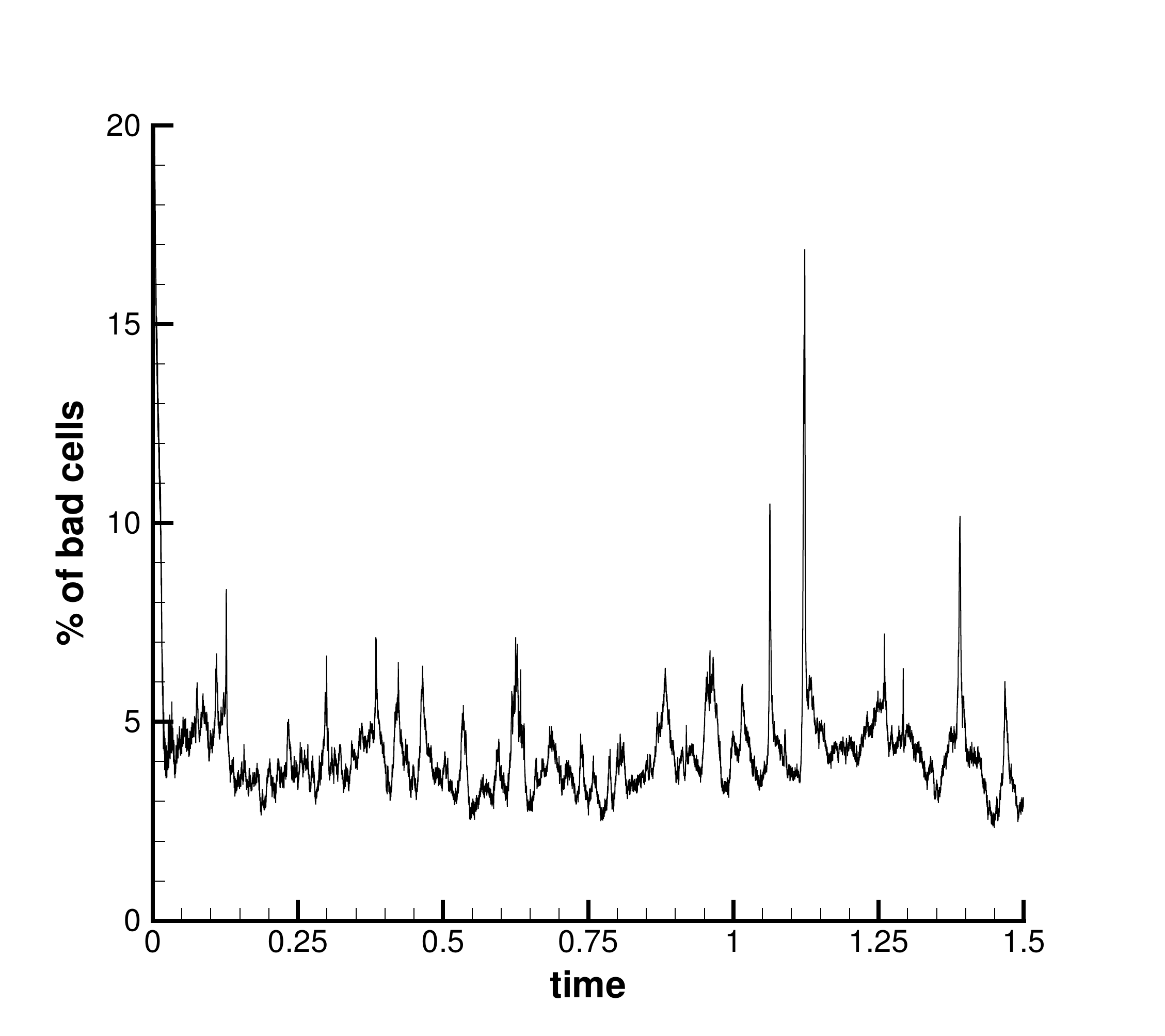} \\     
    \end{tabular} 
    \caption{Cantilever thick beam ---
      Comparison between the MOOD cascades:
      $\mathbb{P}_1 \rightarrow \mathbb{P}_0$ (LAM $\mathbb{P}_0$-lim) and
      $\mathbb{P}_1 \rightarrow \mathbb{P}_1^\text{lim} \rightarrow \mathbb{P}_0$ (LAM $\mathbb{P}_1$-lim) for the computed numerical dissipation as a function of time (left) and percentage of bad cells detected at each time step (right).}
    \label{fig.BendCol2D_diss_bad-cell}
  \end{center}
\end{figure}

Then in figure~\ref{fig.BendCol2D_diss_bad-cell} we present the computed numerical dissipation as a function of time for the two cascades, where about $60\%$ less dissipation is obtained by the current 3 scheme cascade. At last the right panel presents the percentage of troubled cells encountered as a function of time. On average about $5\%$ of cells are re-computed at each timestep.
% ---- FIG ---------
% ---- FIG ---------

%
% TEST # 4 :  Blakes
%
\subsection{Blake's problem} \label{ssec.Blake}	
\begin{figure}[!htbp]
  \begin{center}
    \begin{tabular}{c} 
      \includegraphics[width=0.9\textwidth]{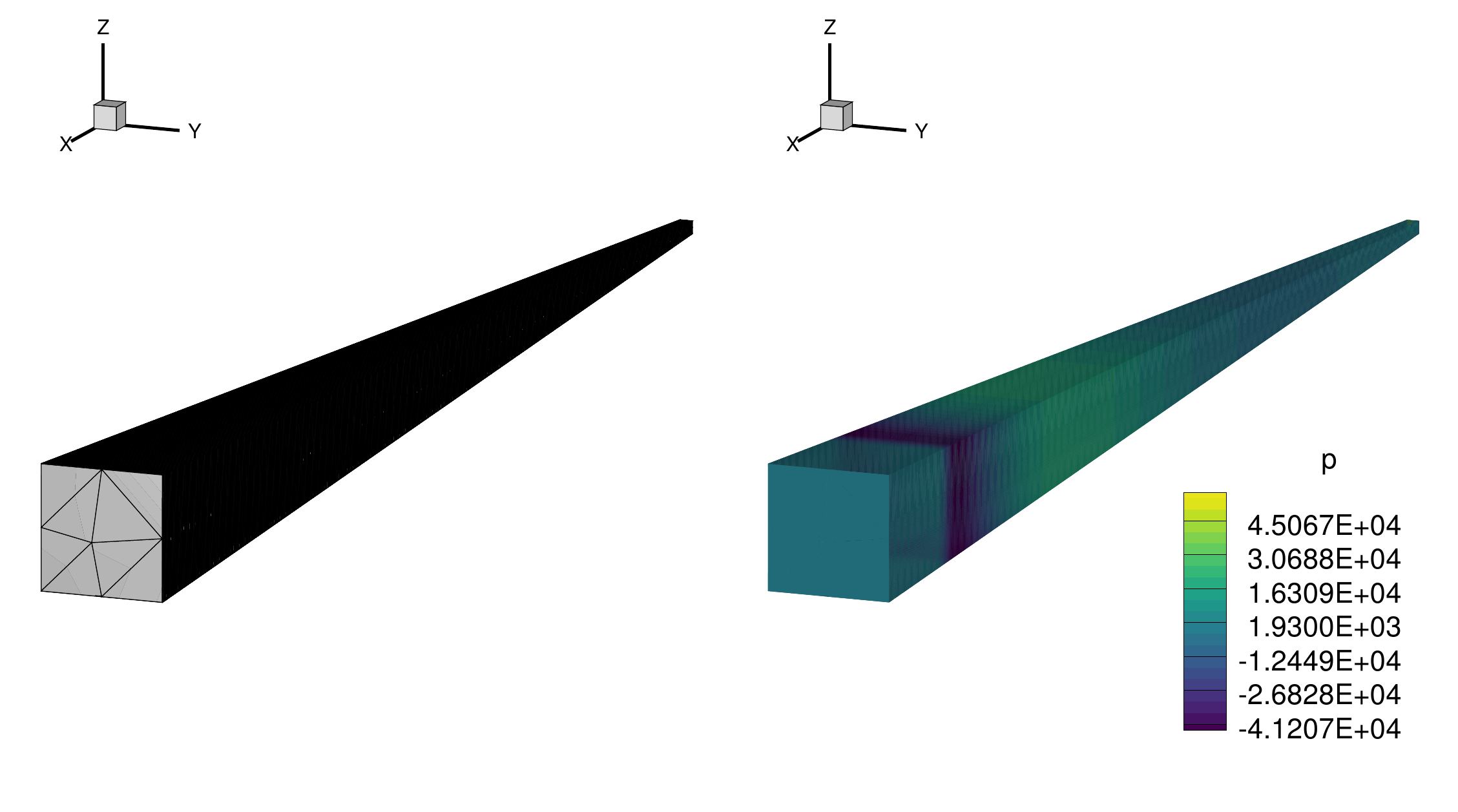}      
    \end{tabular} 
    \caption{Blake's problem --- Computational mesh of the needle domain $\Omega=[r,\theta,\phi]=[0.9,\pi/180,\pi/180]$ with $h=1/1000$ (left) and pressure distribution at the final time $t_\text{final}=1.6 \cdot 10^{-4}$ (right). }
    \label{fig.Blake3D}
  \end{center}
\end{figure}

\begin{figure}[!htbp]
  \begin{center}
    \begin{tabular}{cc} 
      \includegraphics[width=0.47\textwidth]{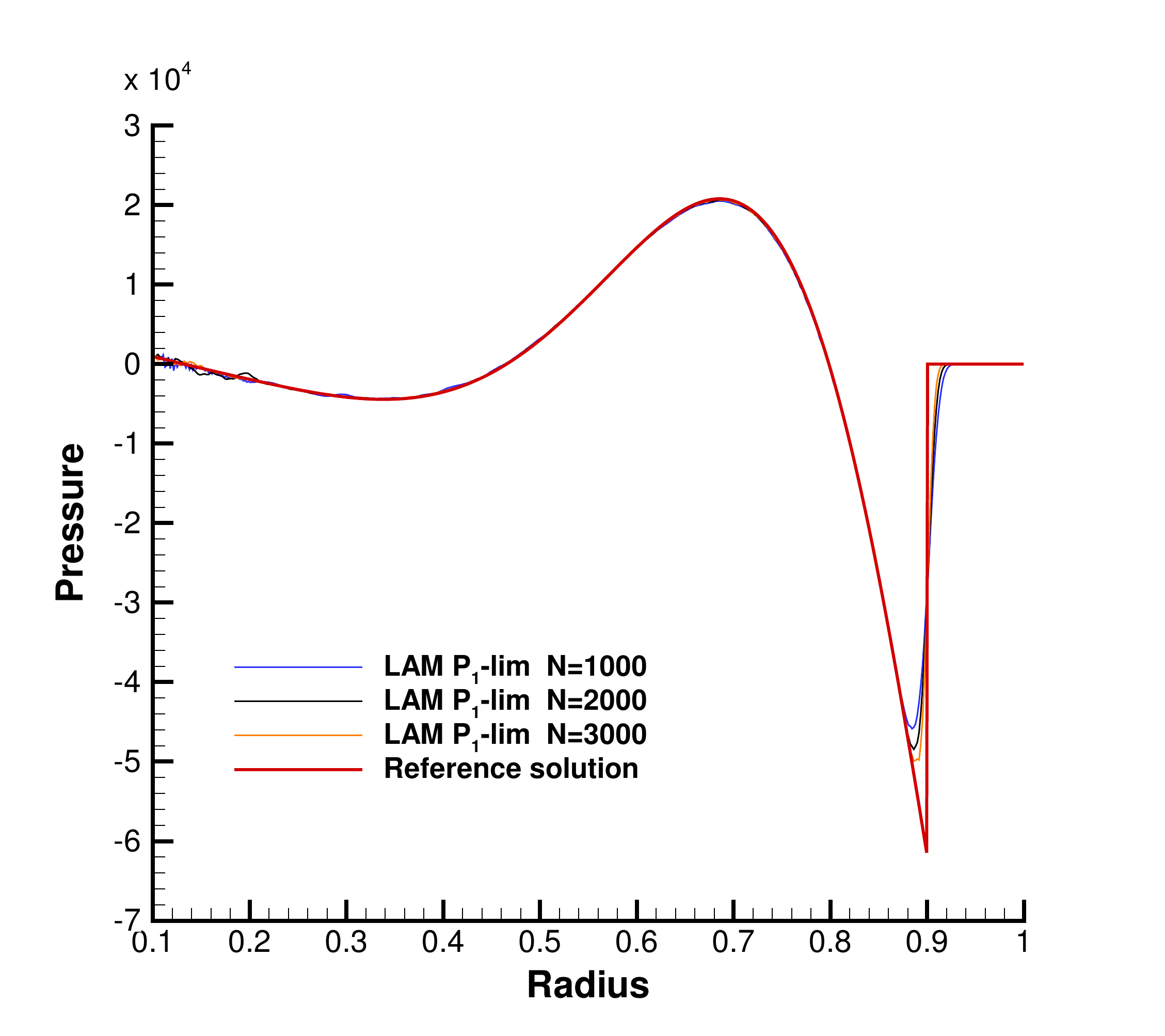} &       
      \includegraphics[width=0.47\textwidth]{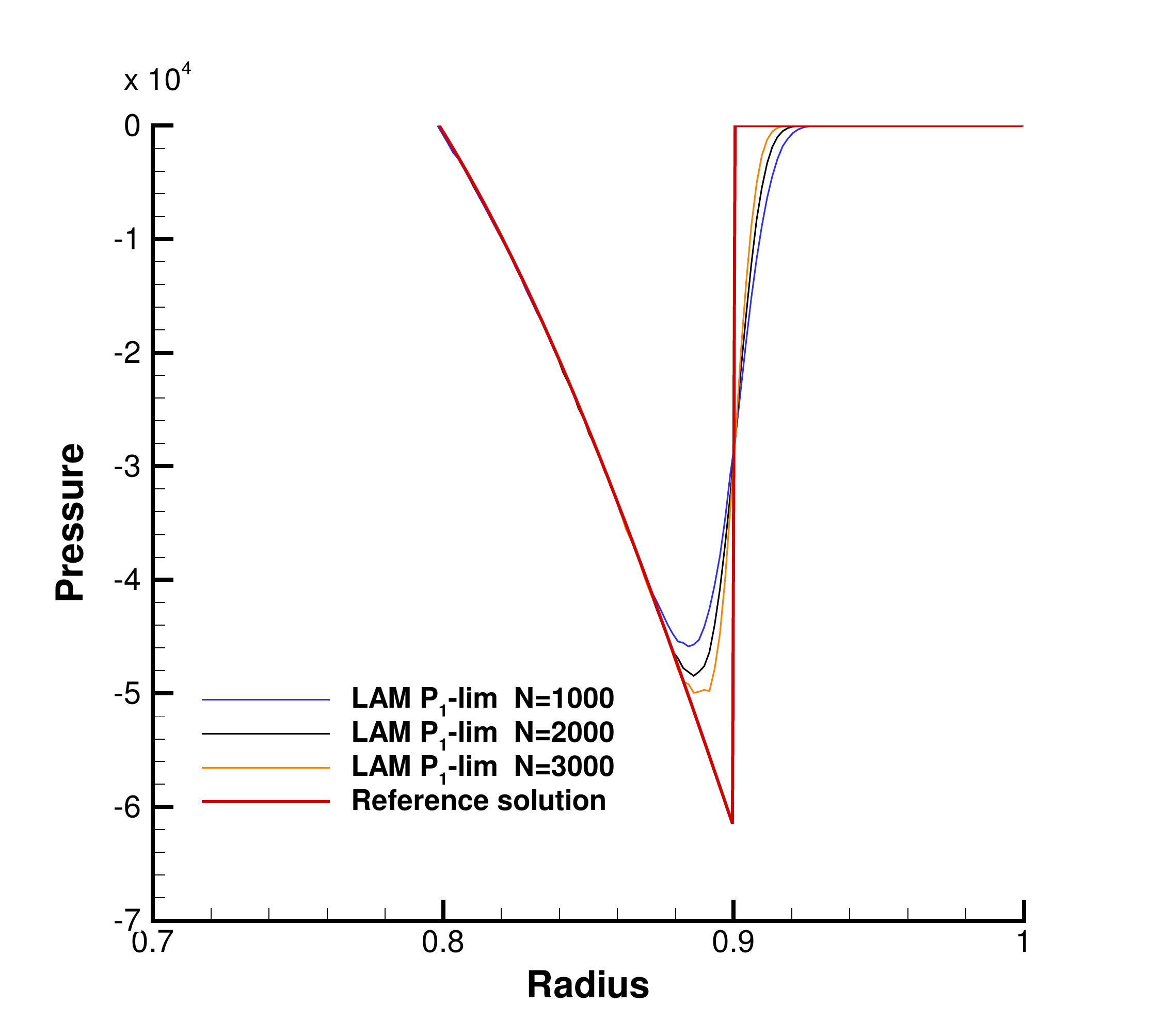} \\ 
      \includegraphics[width=0.47\textwidth]{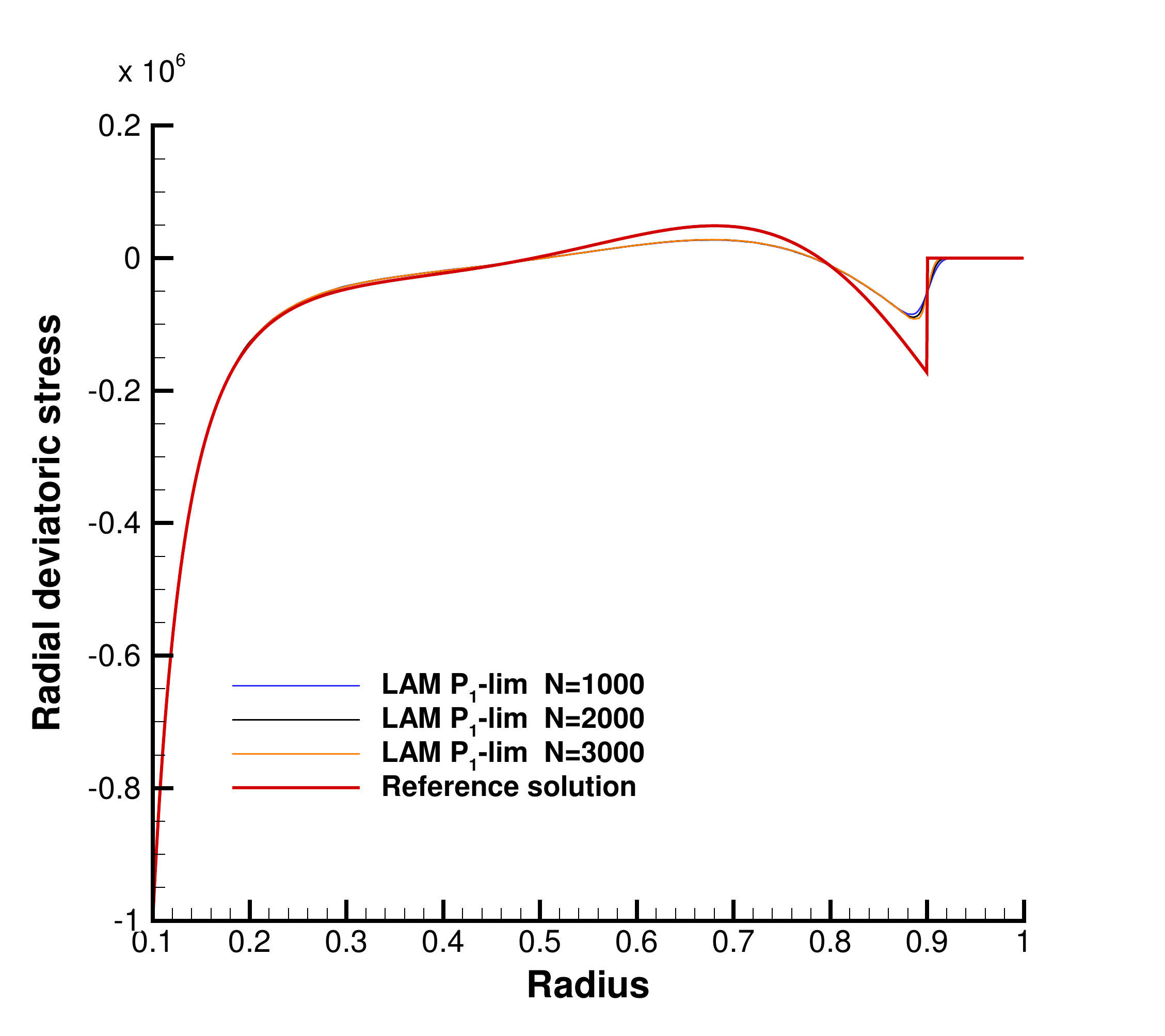} &       
      \includegraphics[width=0.47\textwidth]{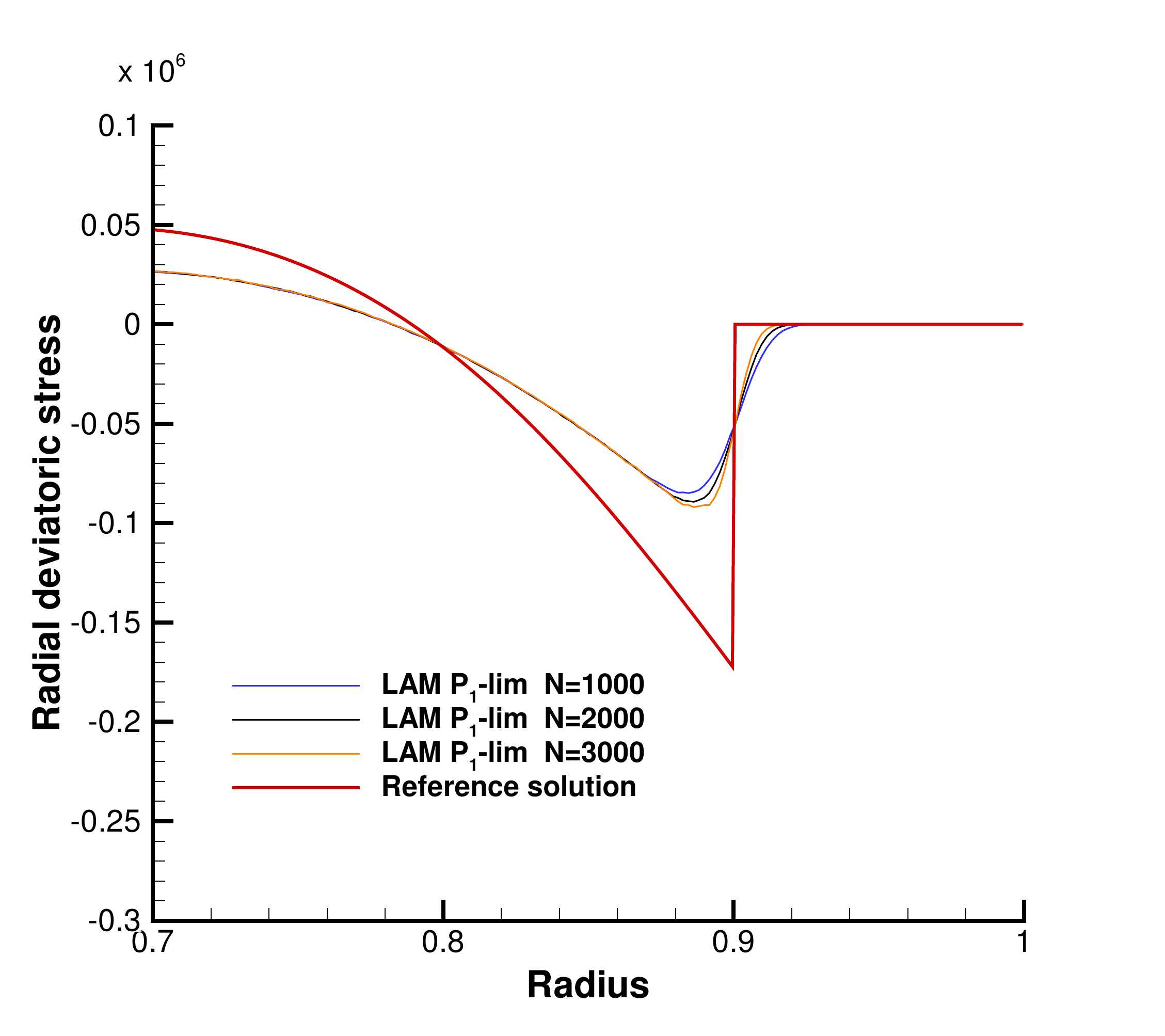} \\     
    \end{tabular} 
    \caption{Blake's problem--- Convergence of the second order solution towards the reference solution for the radial pressure (top row) and radial deviatoric stress (bottom row) at time $t_\text{final}=1.6 \cdot 10^{-4}$ (left) and zoom across the shock (right).}
    \label{fig.Blake-pressure}
  \end{center}
\end{figure}

% Problem set-up
Blake's problem is a classical spherical test derived from the small strain linear elasticity theory 
\cite{Kamm_Blake09}. The domain is a shell of inner radius $r_{in}= 0.1~\text{m}$ and outer radius $r_{out}=1~\text{m}$.
The shell material is isotropic with parameters: $\rho_0 = 3000~\text{kg.m}^{-3}$,
Young's modulus $E = 62.5\cdot 10^9~\text{Pa}$ and Poisson's ratio $\nu = 0.25$.
The inner face of the shell is driven by a pressure constrain of magnitude $10^6~\text{Pa}$
whereas the outer face is a stress free boundary condition.
The final time is $t_\text{final}=1.6 \cdot 10^{-4}$.
% Numerical setup
In practice, for computational time reasons, the domain is not a complete
shell but a needle-like domain of one degree aperture angle. All the boundary faces
introduced by this geometrical simplification are then symmetry boundary
conditions. As such the computational domain is defined by 
$\Omega=[r,\theta,\phi]=[0.9,\pi/180,\pi/180]$ and three meshes with characteristics length $h=1/N_s$ are considered ($N_s=1000 \cdot s$ cells with $s=1,2,3$). 
% Computational mesh
An additional difficulty arise in the context of three-dimensional unstructured meshes, which is related to the spatial discretization of the needle-like computational domain for the Blake problem. In order to avoid ill-conditioned reconstruction matrices due to the high difference in cell size between elements close to the origin of the needle and the ones very far from that location, the entire computational domain has to be mapped onto a reference system $[\bar{r},\bar{\theta},\bar{\psi}]$ such that all coordinates are defined within the interval $[0;1]$. This is sufficient to carry out a second order reconstruction on a more uniform tessellation of the domain with tetrahedra.
% Figure
In figure~\ref{fig.Blake3D} we present the mesh of the needle and the pressure distribution
at final time as illustration with $N_s=1000$. In order to provide more quantitative analysis, 
in figure~\ref{fig.Blake-pressure} we display the numerical results for the pressure and radial deviatoric stress (and zooms) as a function of radius
for a sequence of meshes: $N_1=1000$, $N_2=2000$ and $N_3=3000$. The solution is then compared against the reference solution. We can observe not only accuracy but also convergence even though is some perturbations are seen for small radius on pressure variables.
% ---- FIG ---------
% ---- FIG ---------
% ---- FIG ---------
% ---- FIG ---------

%
% TEST # 5 : Twisting column
%
\subsection{Twisting column} \label{ssec.TwistCol}
% ---- FIG ---------
\begin{figure}[!htbp]
  \begin{center}
    \includegraphics[width=0.6\textwidth]{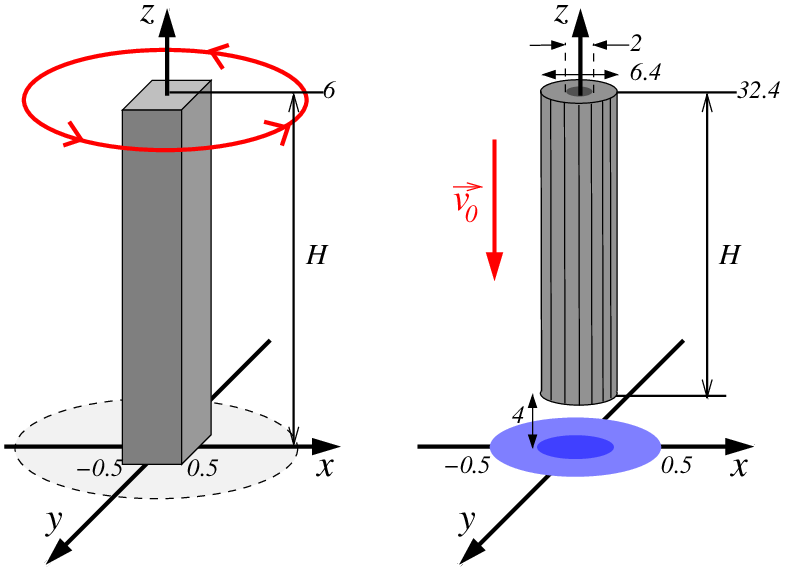}
    \caption{Sketch for the twisting column in section~\ref{ssec.TwistCol} (left)
    and the rebound of a hollow circular bar from section~\ref{ssec.BarRebound} (right).} 
    \label{fig:sketch_column}
  \end{center}
\end{figure}
\begin{figure}[!htbp]
  \begin{center}
%    \begin{tabular}{ccc} 
%      \includegraphics[width=0.53\textwidth]{TwistColumn3D_t00375} &  
%      \includegraphics[width=0.53\textwidth]{TwistColumn3D_t0075} &
%      \includegraphics[width=0.53\textwidth]{TwistColumn3D_t01125} \\
%      \includegraphics[width=0.53\textwidth]{TwistColumn3D_t015} &
%      \includegraphics[width=0.53\textwidth]{TwistColumn3D_t01875} &  
%      \includegraphics[width=0.53\textwidth]{TwistColumn3D_t0225} \\
%      \includegraphics[width=0.53\textwidth]{TwistColumn3D_t02625} &
%      \includegraphics[width=0.53\textwidth]{TwistColumn3D_t03} &   
%      \includegraphics[width=0.53\textwidth]{TwistColumn3D_t00} \\
%    \end{tabular}  
%    \begin{tabular}{cccccccc}
%      \hspace{-0.5cm}
%      \includegraphics[width=0.12\textwidth]{twc_t00375} &  
%      \hspace{-0.5cm}
%      \includegraphics[width=0.12\textwidth]{twc_t0075}  & 
%      \hspace{-0.5cm} 
%      \includegraphics[width=0.12\textwidth]{twc_t01125} & 
%      \hspace{-0.5cm}
%      \includegraphics[width=0.12\textwidth]{twc_t015}   & 
%      \hspace{-0.5cm} 
%      \includegraphics[width=0.12\textwidth]{twc_t01875} &  
%      \hspace{-0.5cm}
%      \includegraphics[width=0.12\textwidth]{twc_t02625}  & 
%      \hspace{-0.5cm} 
%      \includegraphics[width=0.12\textwidth]{twc_t03}  &  
%      \hspace{-0.5cm}
%      \includegraphics[width=0.1\textwidth]{twc_legend2}  
%    \end{tabular} 
        \begin{tabular}{cccc} 
         \includegraphics[width=0.22\textwidth]{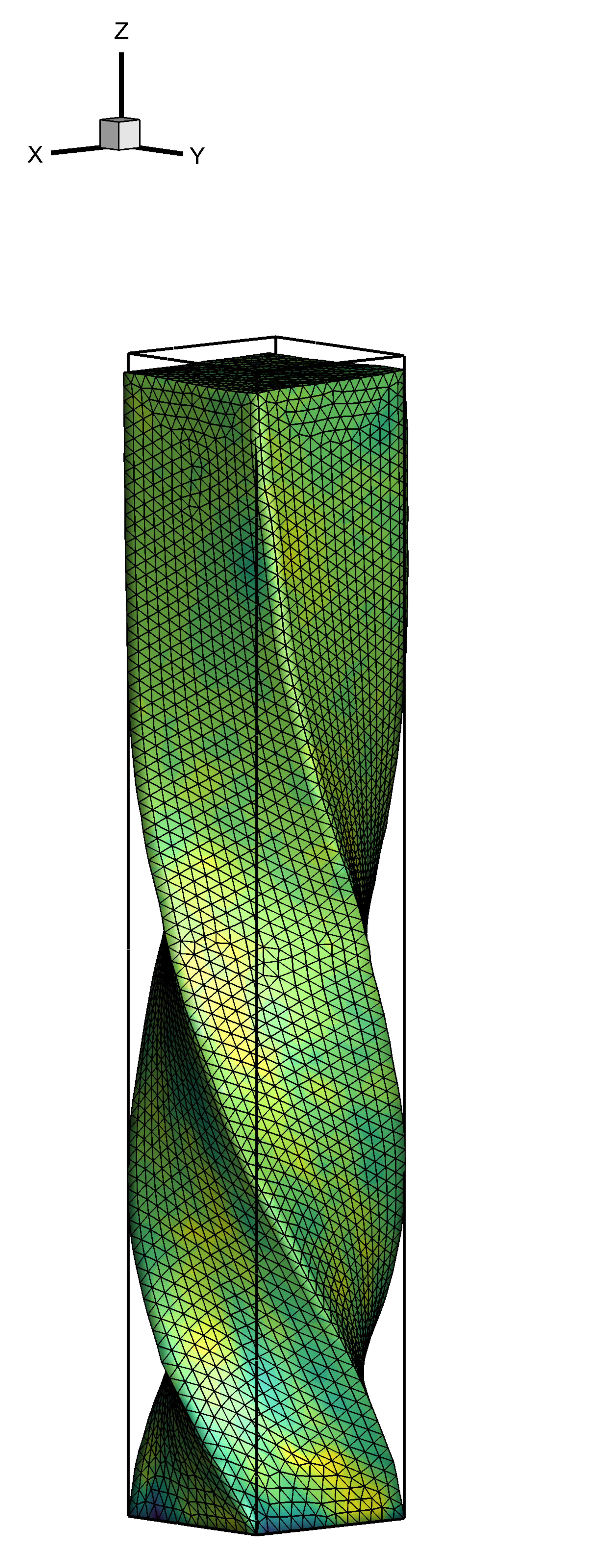} &
         \includegraphics[width=0.22\textwidth]{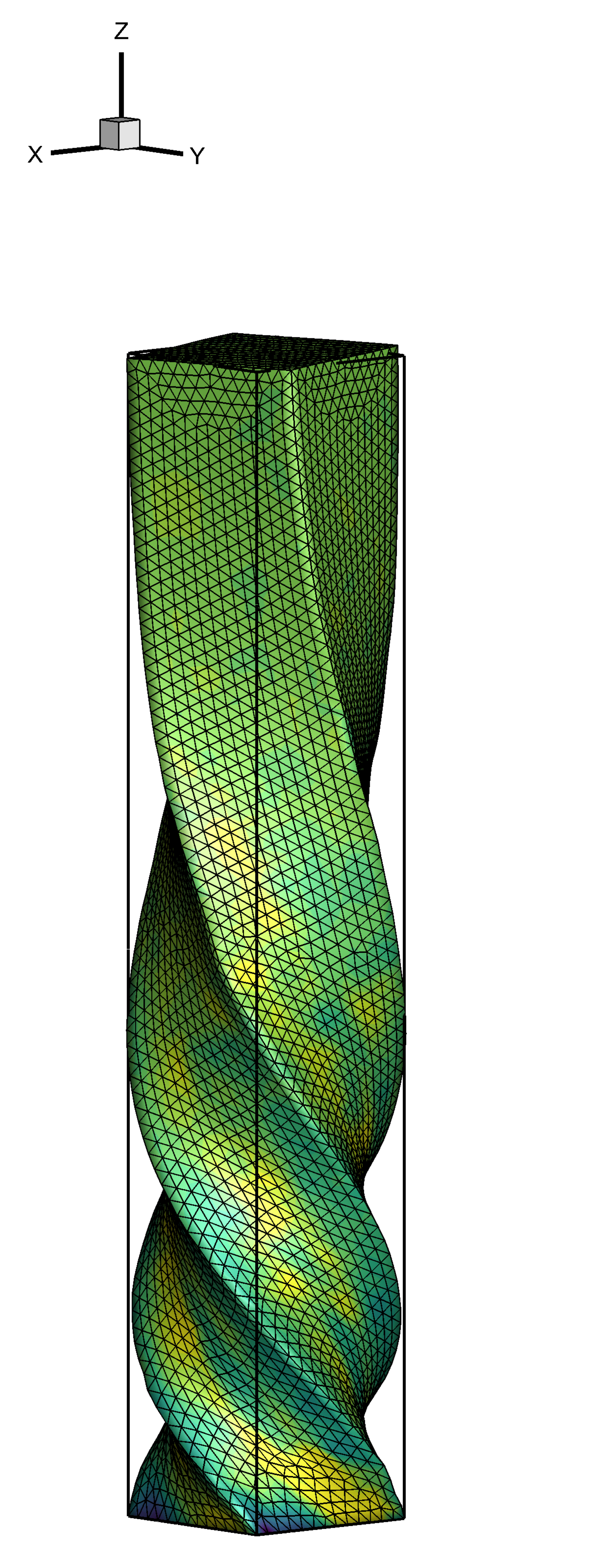} &
         \includegraphics[width=0.22\textwidth]{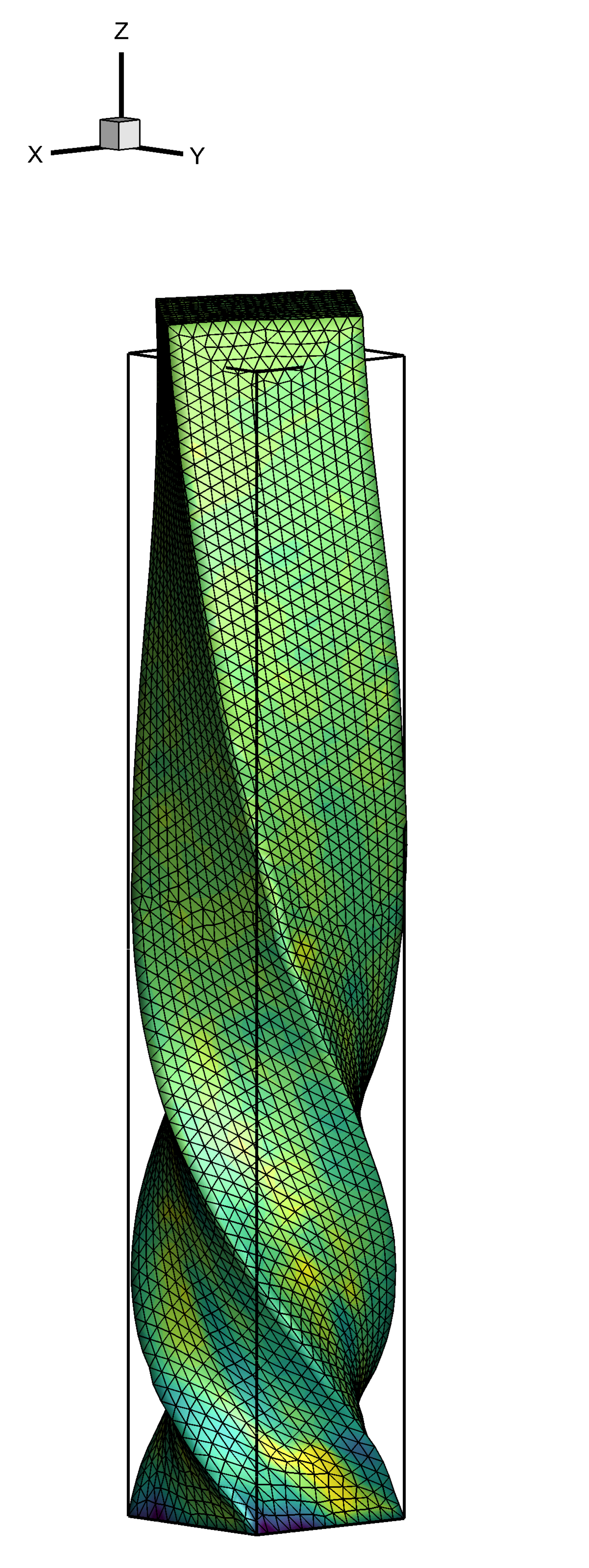} &
         \includegraphics[width=0.22\textwidth]{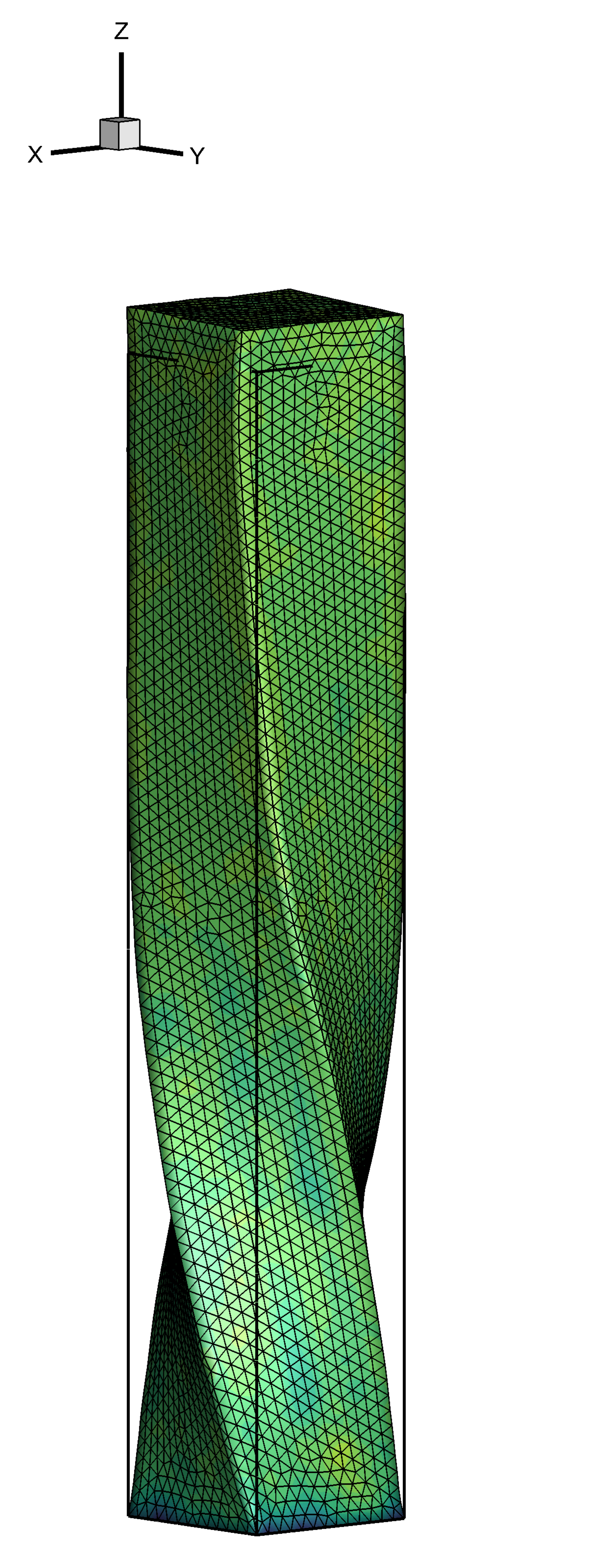} \\
         \includegraphics[width=0.22\textwidth]{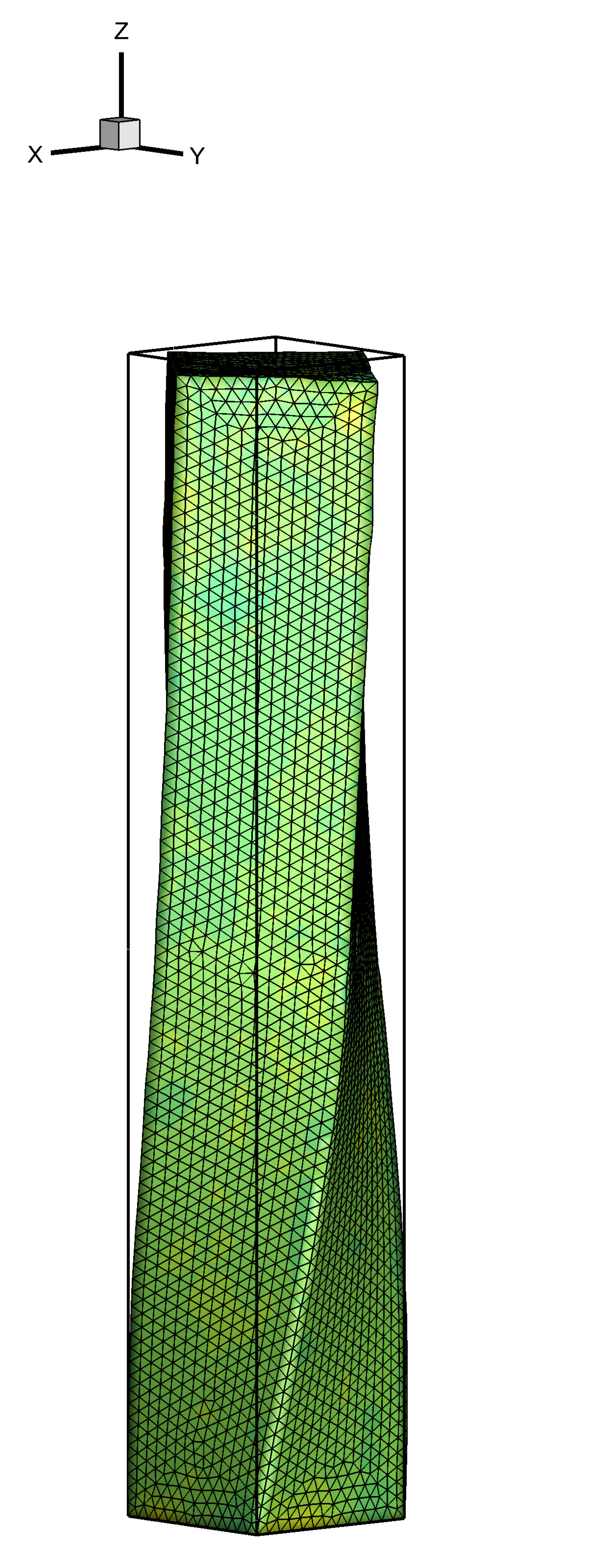} &
         \includegraphics[width=0.22\textwidth]{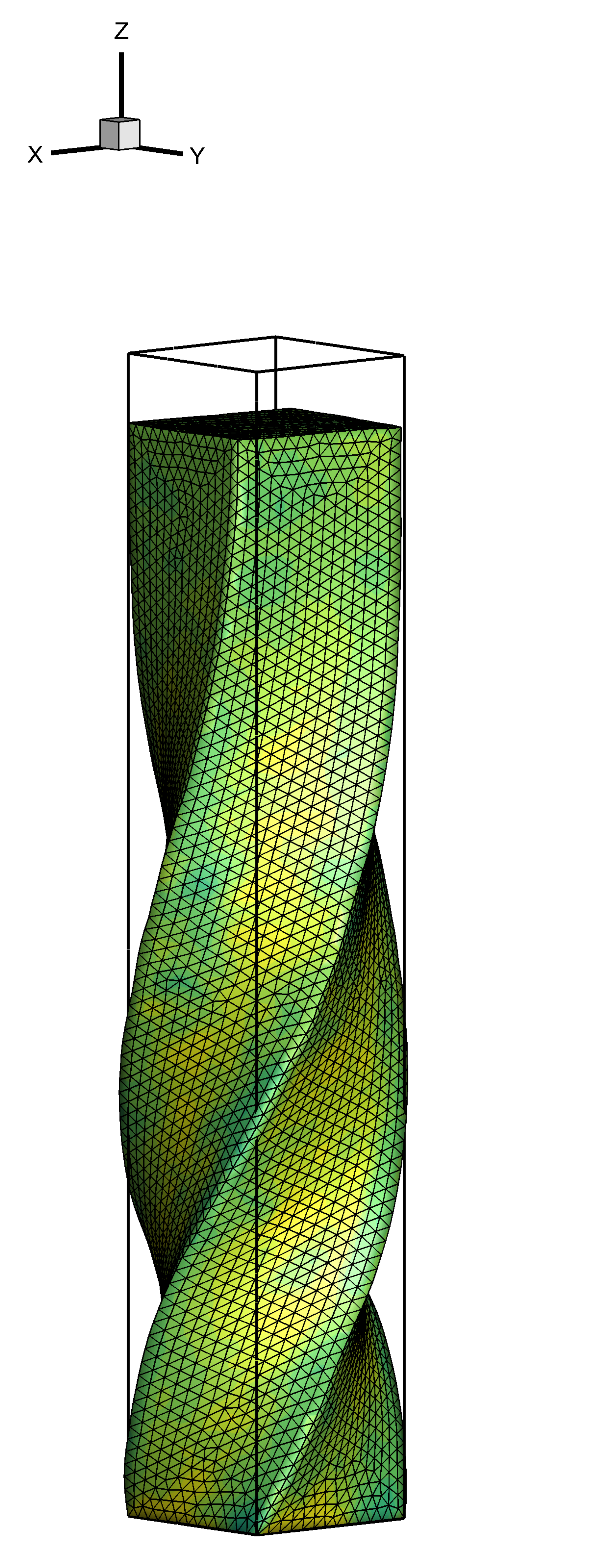} &
         \includegraphics[width=0.22\textwidth]{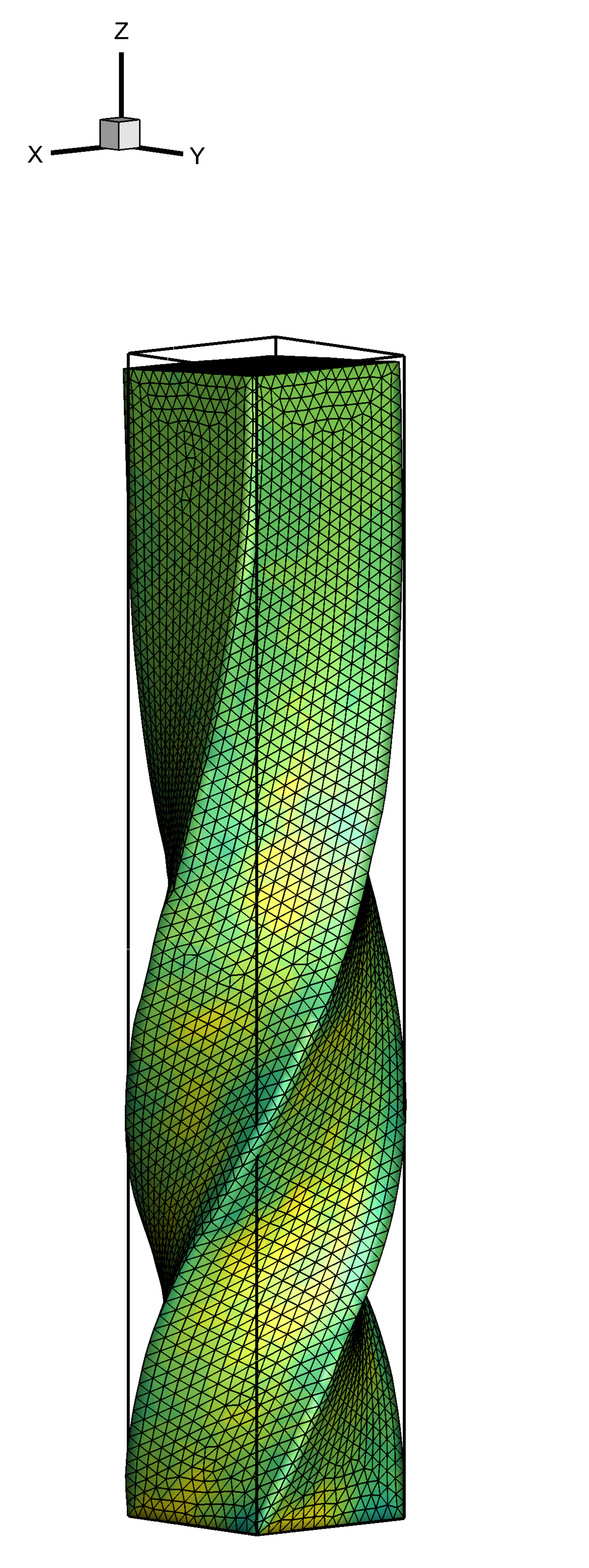} &
         \includegraphics[width=0.22\textwidth]{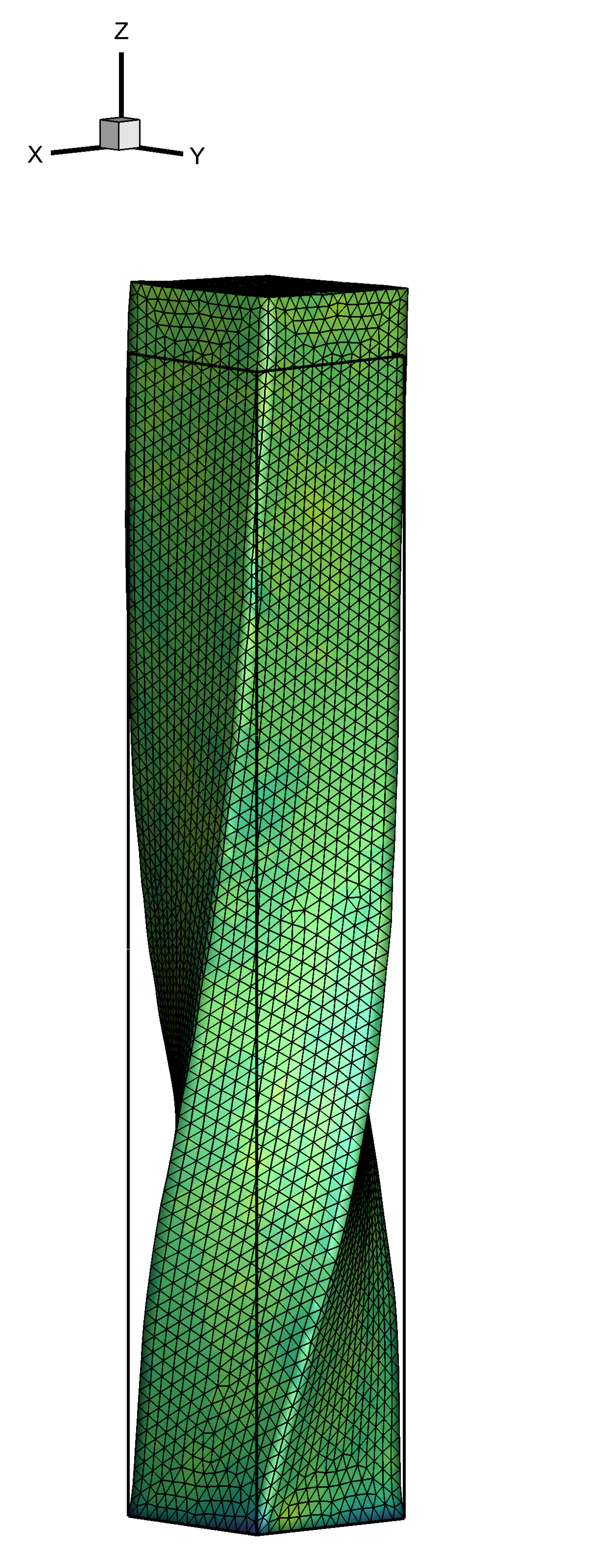} \\
         \multicolumn{4}{c}{\includegraphics[width=0.8\textwidth]{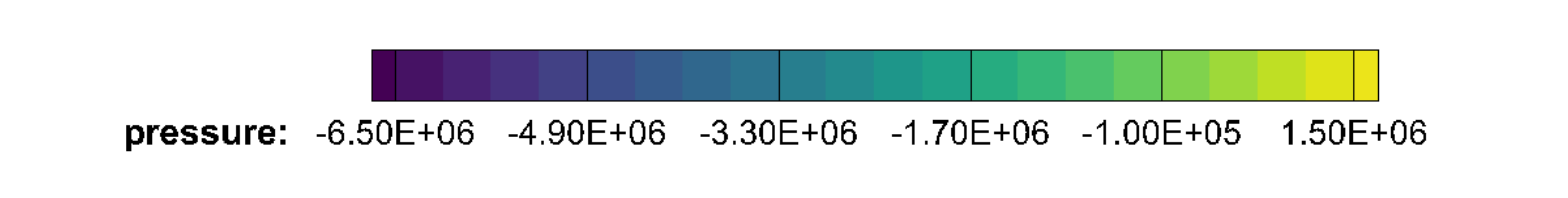}}\\
    \end{tabular} 
    \caption{Twisting column ---
      Beam shape and pressure distribution at output times $t=0.00375$, $t=0.075$, $t=0.1125$, $t=0.15$, $t=0.1875$, $t=0.225$, $t=0.2625$ and $t=0.3$ (from top left to bottom right).
      The shape is compared with respect to the initial configuration (hollow box). %shown in the bottom right panel with velocity vectors and magnitude of the velocity field at time $t=0$.
    }
    \label{fig.TwistCol3D}
  \end{center}
\end{figure}
\begin{figure}[!htbp]
  \begin{center}
    \begin{tabular}{ccc}
    \hspace{-0.95cm} 
      \includegraphics[width=0.37\textwidth]{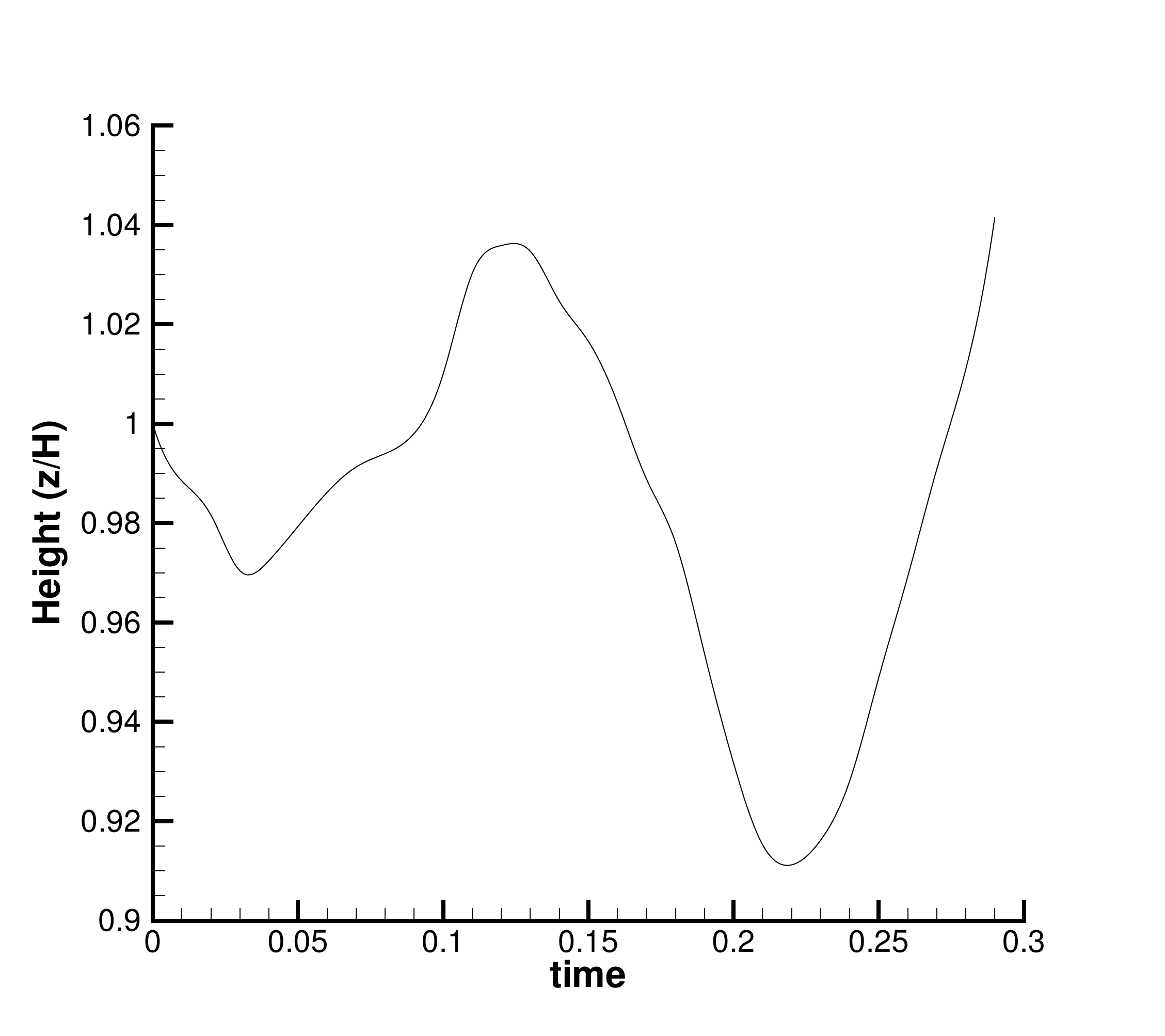} &
    \hspace{-0.95cm} 
      \includegraphics[width=0.37\textwidth]{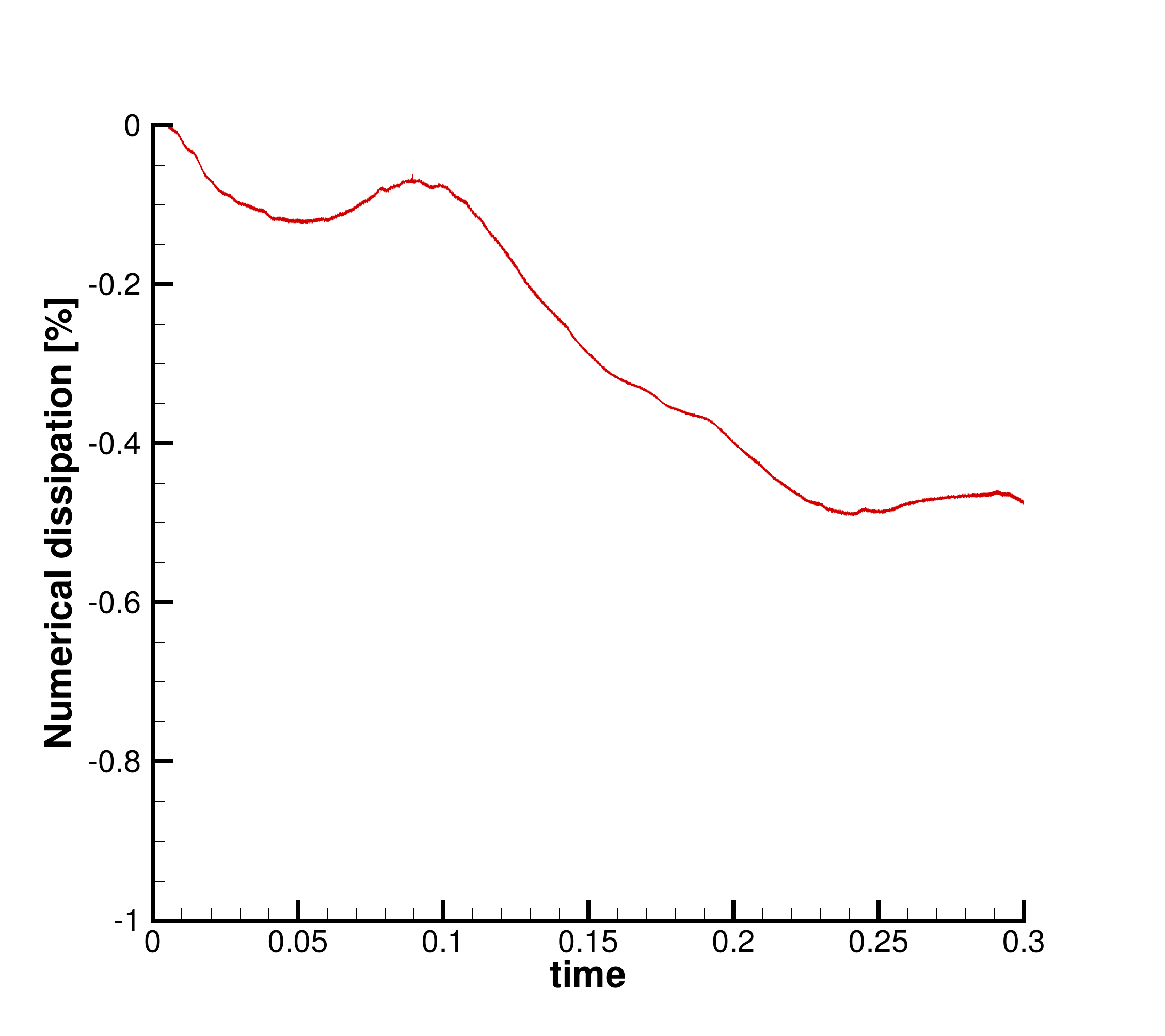} &  
    \hspace{-0.95cm}      
      \includegraphics[width=0.37\textwidth]{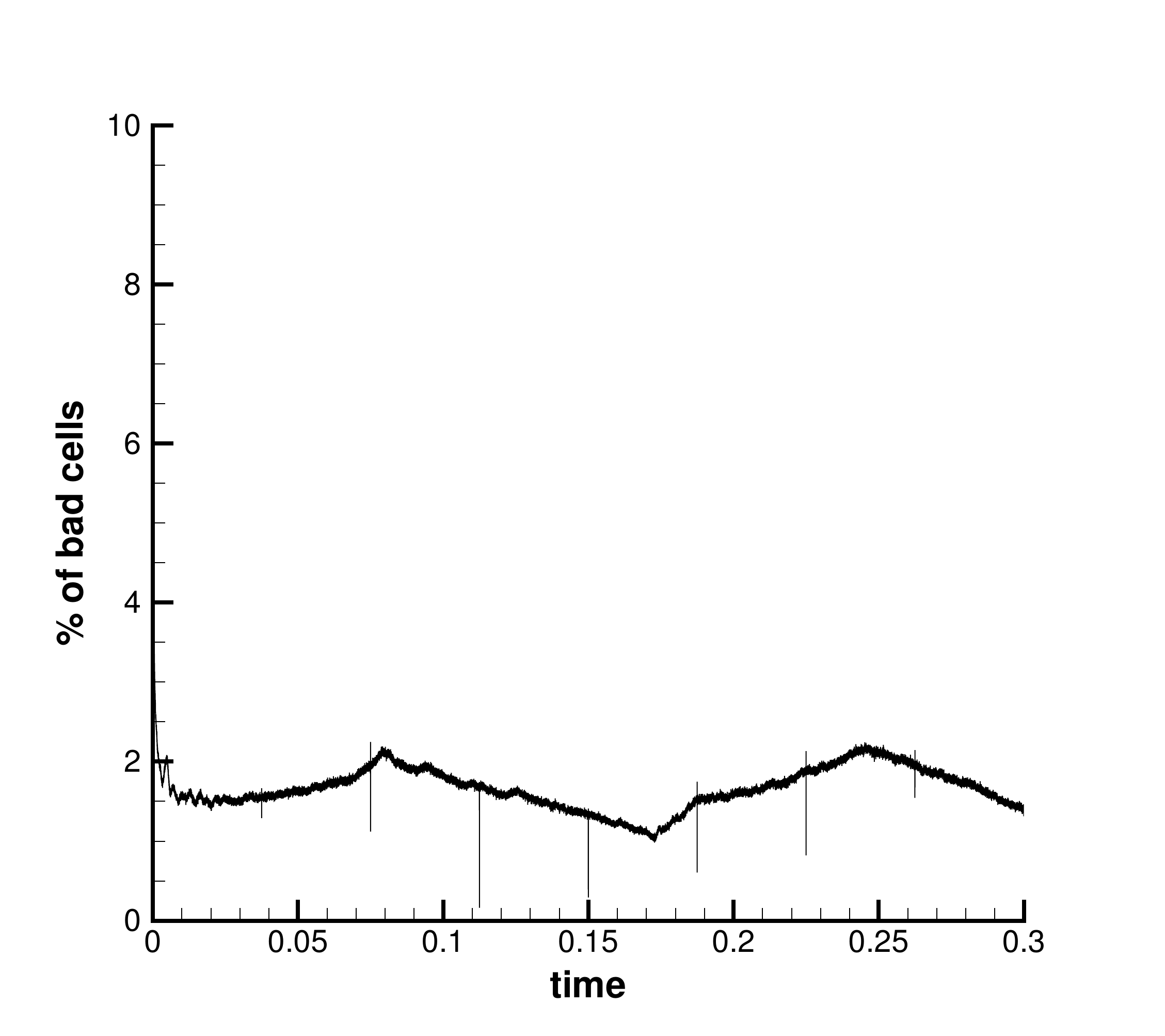}      
    \end{tabular} 
    \caption{Twisting column ---
      Time evolution of non-dimensionalised height of the column measured at initial point $\mathbf{x}_T=(0,0,6)$ (left) ---
      Numerical dissipation of second order scheme (center) ---
      Percentage of bad cells detected at each time step (right).}
    \label{fig.TwistCol3D_diss_bad-cell}
  \end{center}
\end{figure}

% ---- FIG ---------
%\walter{setting from \cite{CCL2020}. Note that there is an error in IC, thus we have $\mathbf{V}(\vec{x},0)=100 \sin \left(\frac{\pi z}{2L}\right) \, \left(0,y,.x\right)^\top$}
%
% Initialization
A twisting column test case aims at examining the effectiveness of the proposed methodology in highly nonlinear scenarios, see \cite{Haider_2018} and the reference therein.
An initial unit squared cross section column of height $H = 6$~m is considered, $\Omega=[-0.5;0.5]\times[-0.5;0.5]\times[0;6]$.
The $z=0$ face of the column is embedded into a wall.
An initial sinusoidal angular velocity field relative to the origin is given by
$\vec{v}_0 = 100\sin(\pi \frac{z}{2H}) ( y, -x, 0)^t$~rad/s, see figure~\ref{fig:sketch_column}.
The main objective of this problem is to assess the capability of the proposed methodology to still perform when approaching the limit of incompressibility.
% Numerical setup
A neo-Hookean material is used with material density $\rho_0 = 1100$~kg/m$^3$, Young's modulus $E = 1.7 \cdot 10^7$~Pa and Poisson's ratio $\nu = 0.45$.
The simulation is run till time $t_{\text{final}}=0.3$~s. Qualitatively one should observe at time $t\sim 0.1$~s a counter-clockwise rotation
and a severe twist of the column which returns to its initial position at about $t\sim 0.2$~s.
Driven by its own inertia, the bar twists clockwise until the final time.
The mesh of the column is made of  $N_c=119092$ tetrahedra with characteristic length of $1/80$.
Stress free BCs are imposed everywhere apart from the bottom face for which we impose a wall type boundary with zero displacement.
% Results
In figure~\ref{fig.TwistCol3D} we plot the shape of the column colored by the pressure distribution for different output times. The initial column is represented as a hollow bar for comparison purposes.
The main behaviors are reproduced by the numerical simulation.
Notice that there is no spurious oscillations nor suspicious pressure distribution.
% Diagnostics
In figure~\ref{fig.TwistCol3D_diss_bad-cell} we gather several diagnostics of this simulation.
First on the left panel we plot the time evolution of dimensionless height of the column measured at the point
initially located at $\mathbf{x}_T=(0,0,6)$.
Next, in the middle panel, we plot the numerical dissipation of the second-order scheme computed as the percentage of energy loss computed by means of \eqref{eqn.numdiss} as a function of time and observe that at final time only $0.5\%$ is lost.
For a numerical simulation recall that the twisting period does not only depend on the material but also on the numerical dissipation of the scheme. Usually first-order schemes are extremely dissipative and can not perform adequately, i.e the column barely twists.
At last in the right panel we present the percentage of bad cells detected at each time step by the \aposteriori limiting procedure and observe that on average only $2\%$ of the cells are recomputed due to spurious numerical issues.
% ---- FIG ---------
% ---- FIG ---------	
% ---- FIG ---------
% ---- FIG ---------
	
%
% TEST # 6 : Rebound of a hollow circular bar
%
\subsection{Rebound of a hollow circular bar} \label{ssec.BarRebound}
\begin{figure}[!htbp]
  \begin{center}
    \begin{tabular}{cccc} 
      \includegraphics[width=0.24\textwidth]{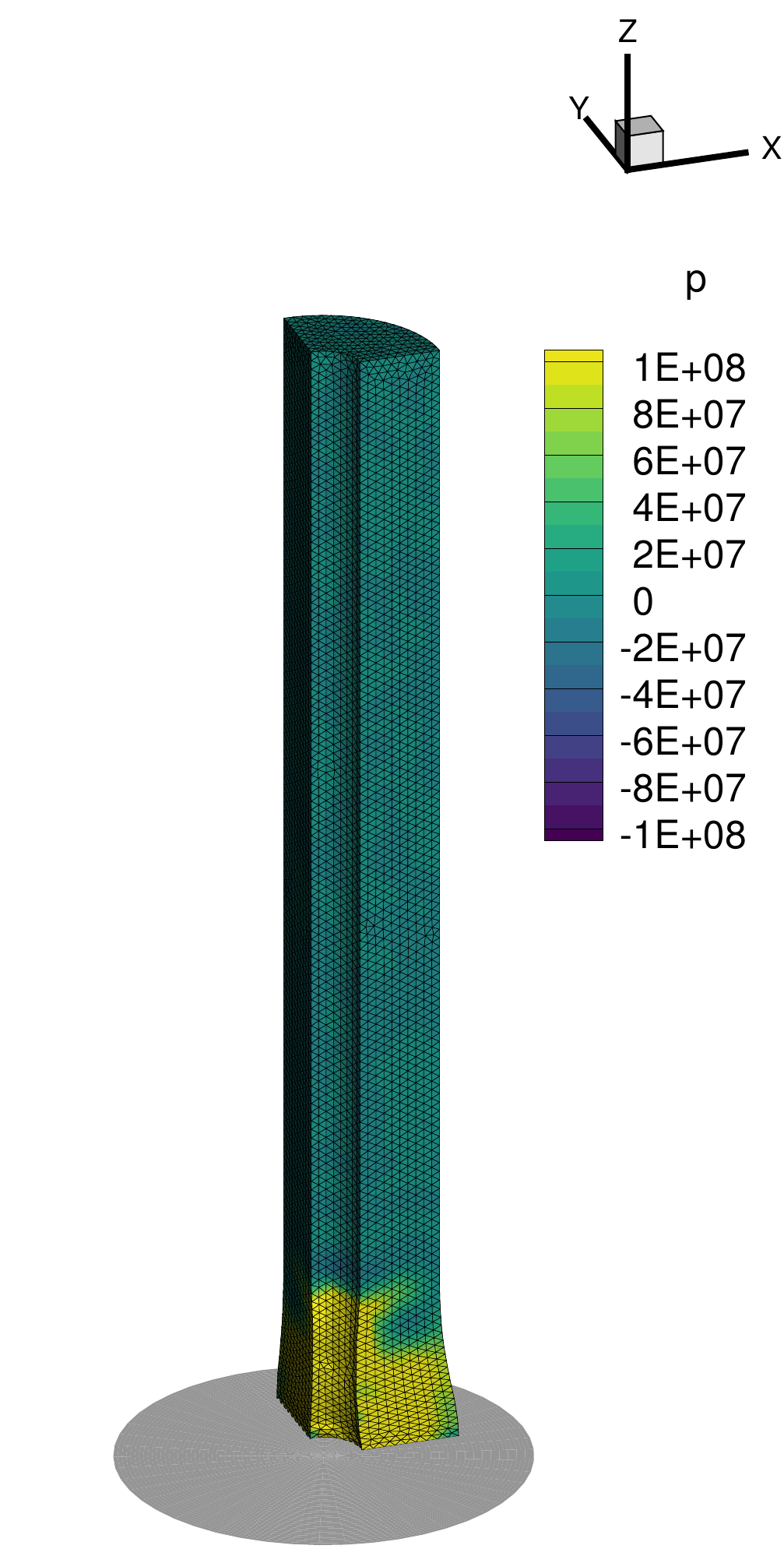} &       
      \includegraphics[width=0.24\textwidth]{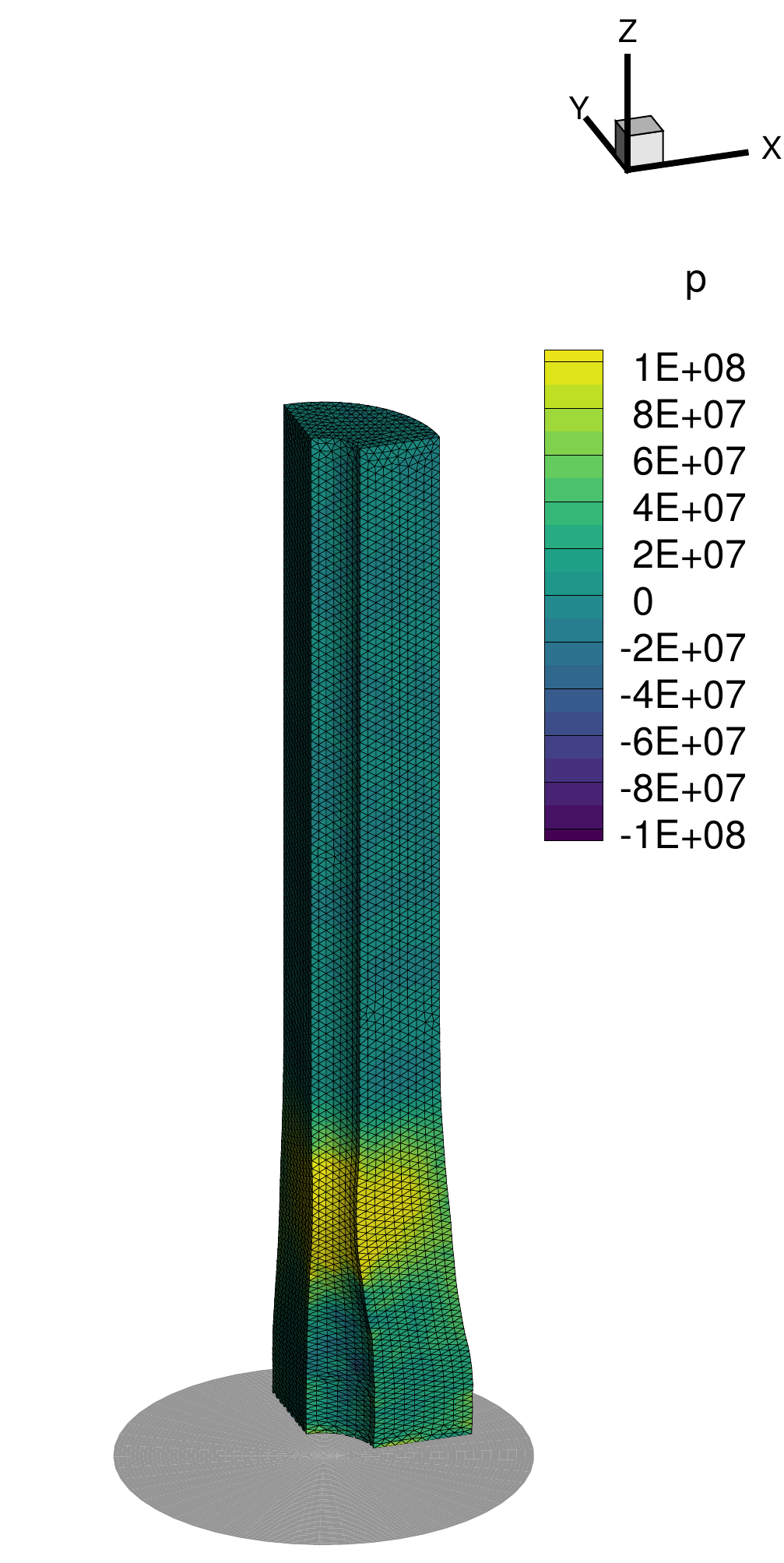} &
      \includegraphics[width=0.24\textwidth]{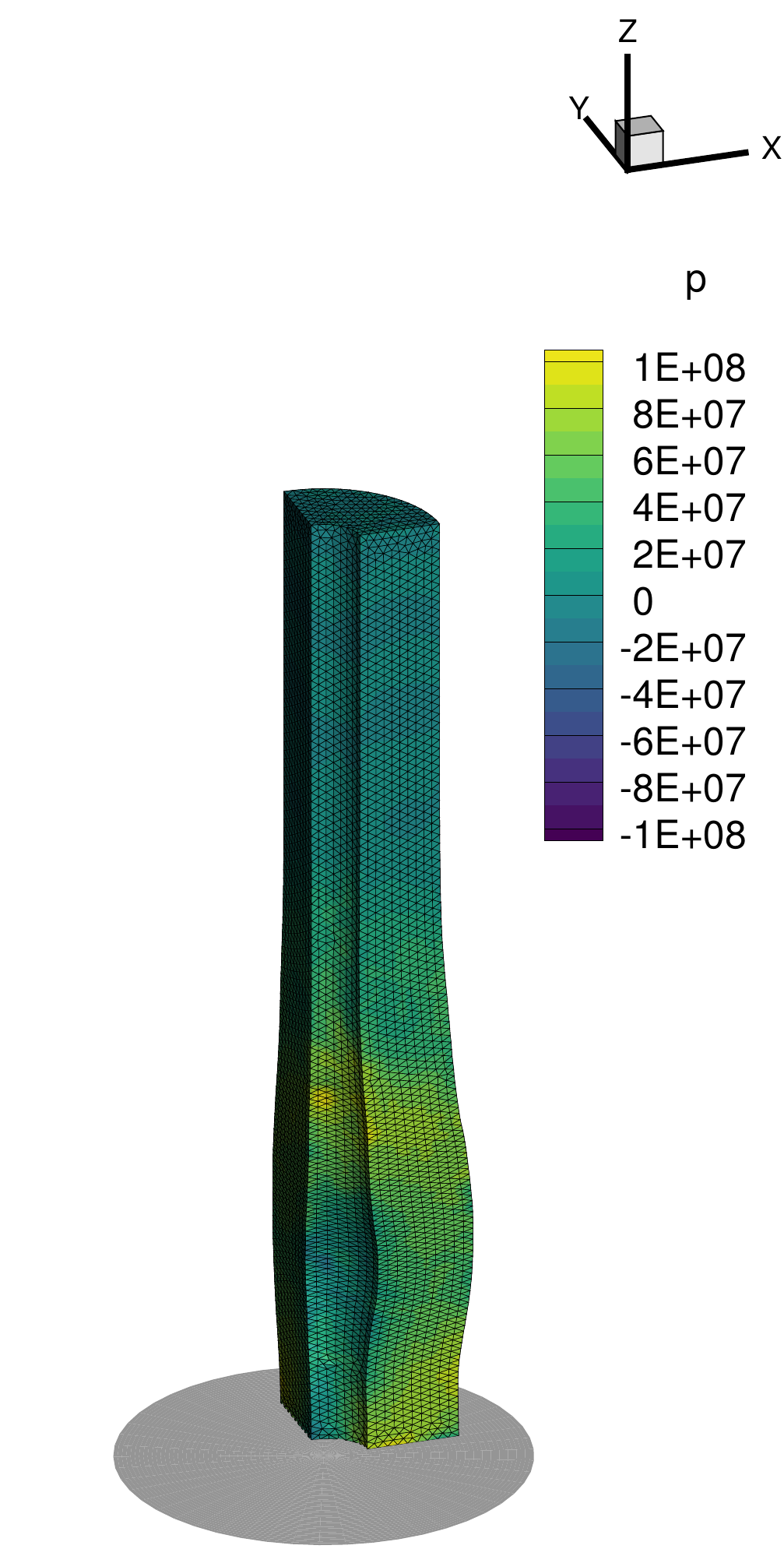} &
      \includegraphics[width=0.24\textwidth]{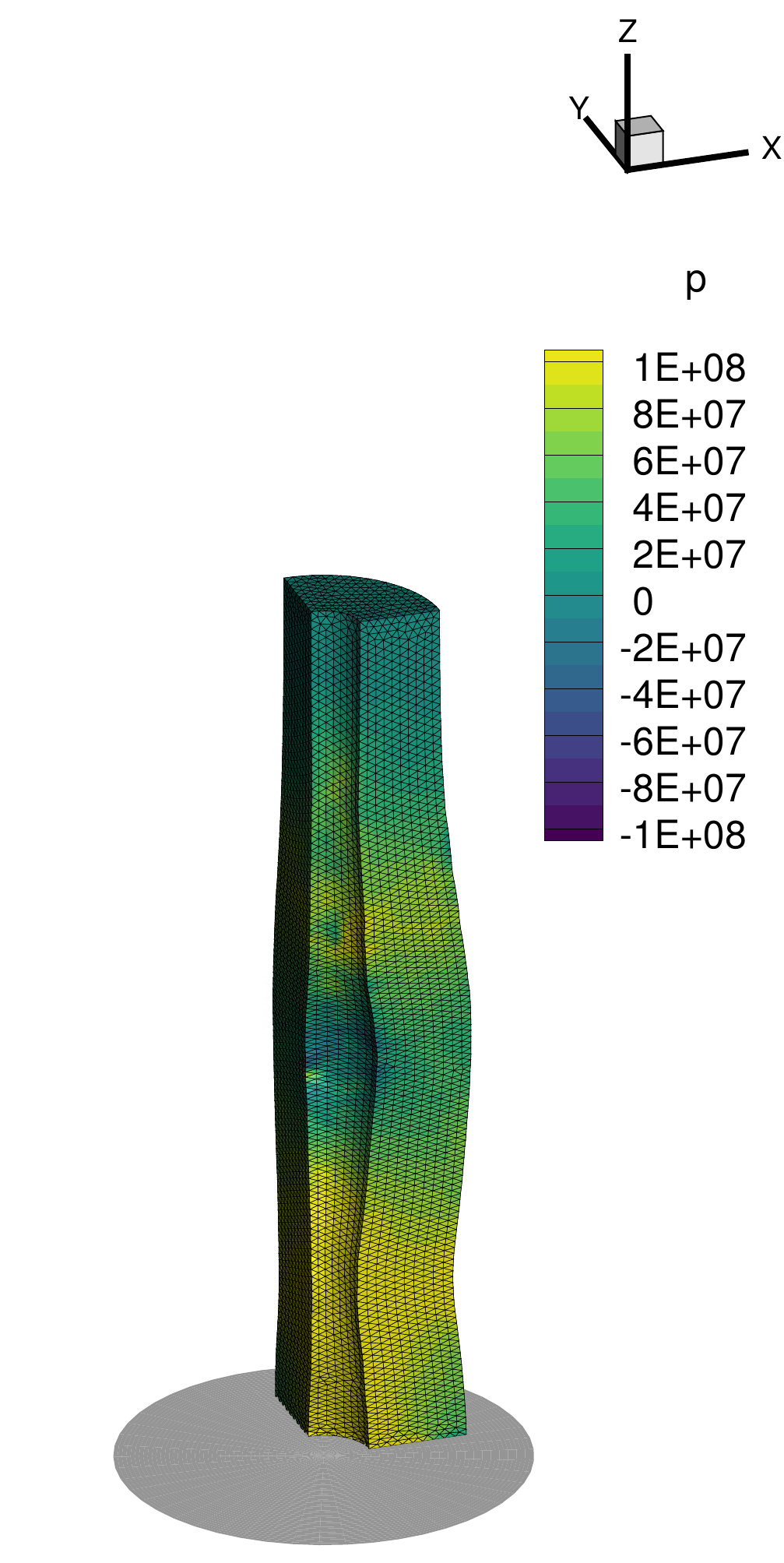} \\
      \includegraphics[width=0.24\textwidth]{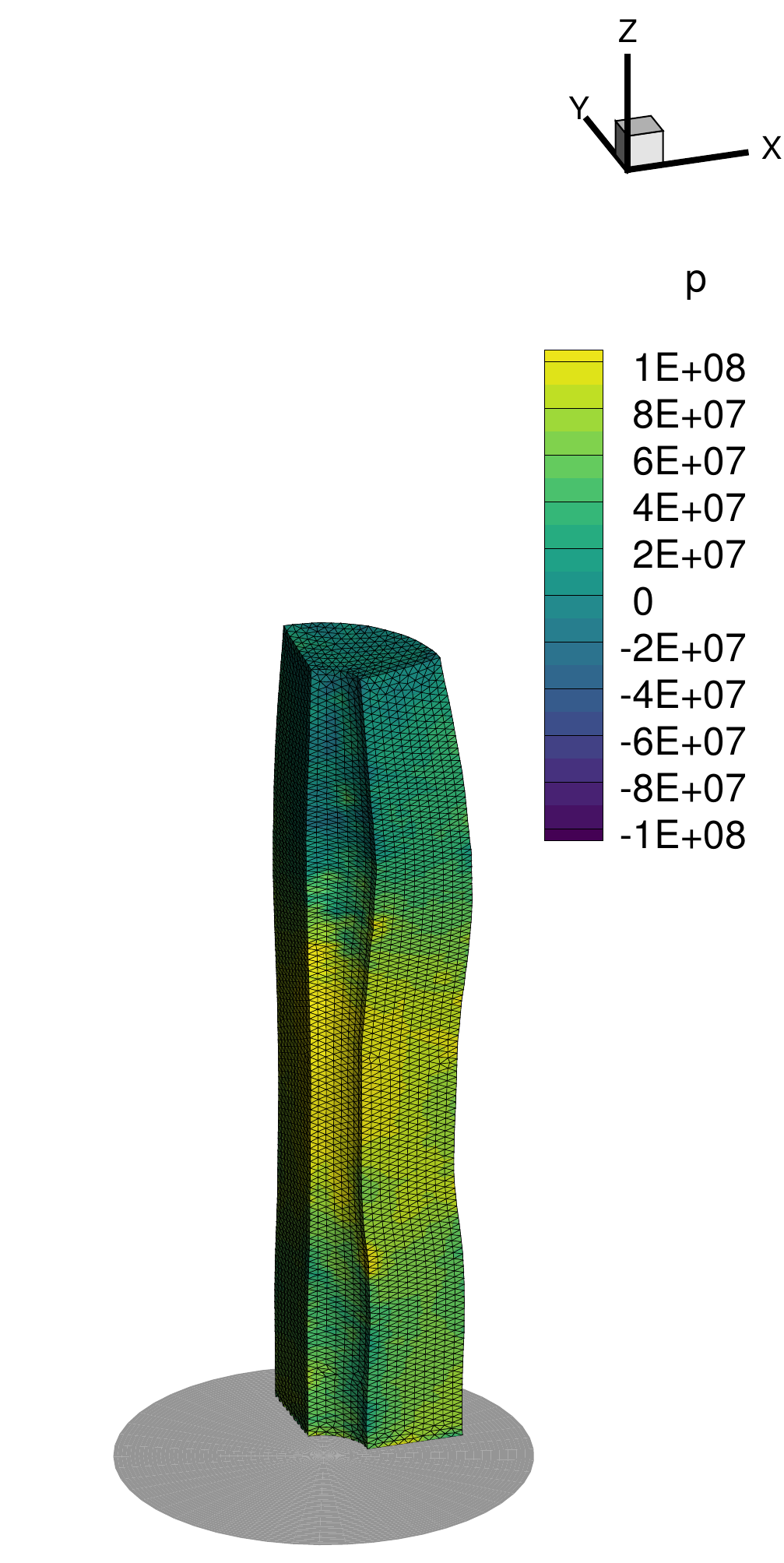} &       
      \includegraphics[width=0.24\textwidth]{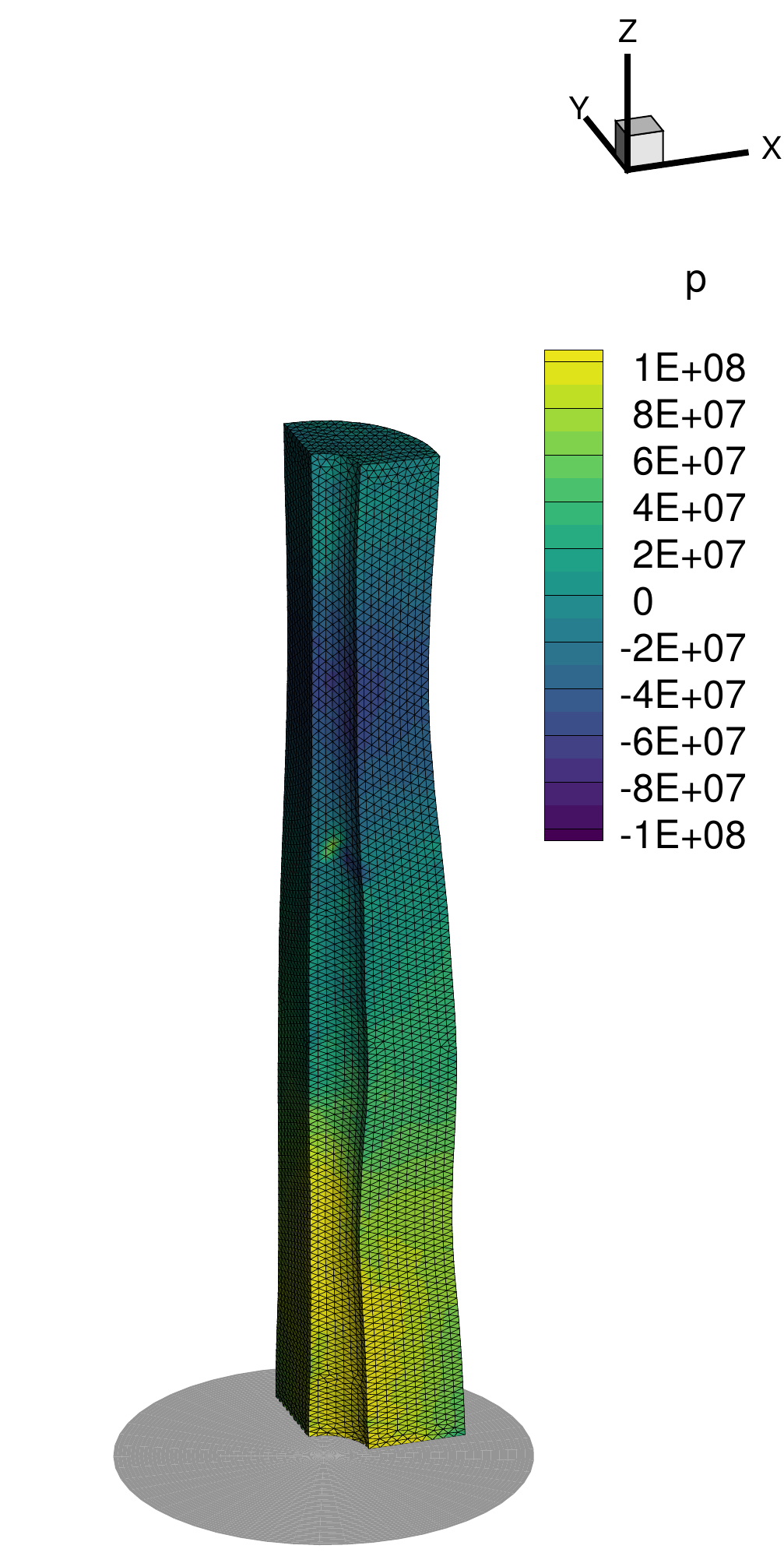} &
      \includegraphics[width=0.24\textwidth]{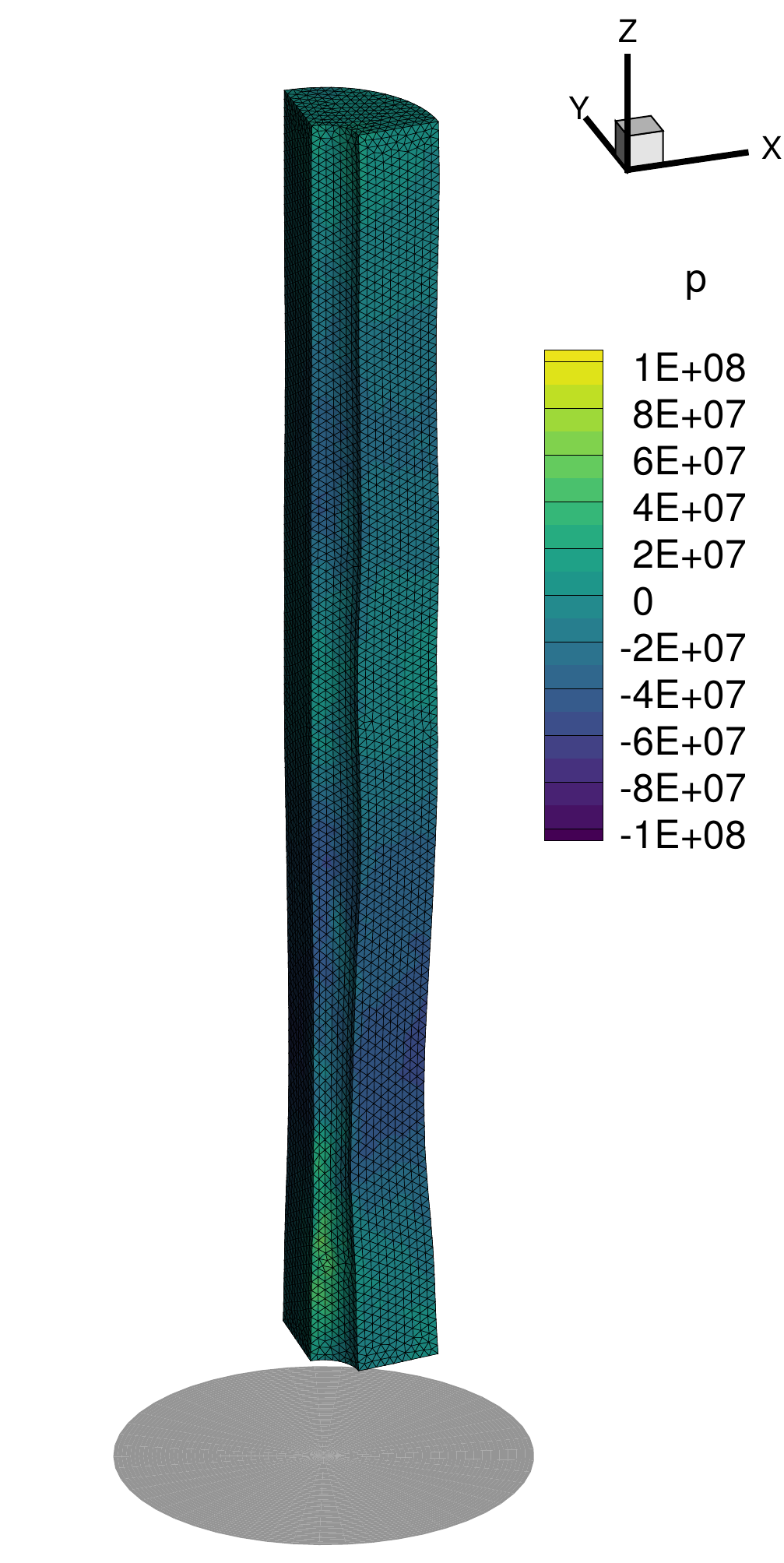} &     
      \includegraphics[width=0.24\textwidth]{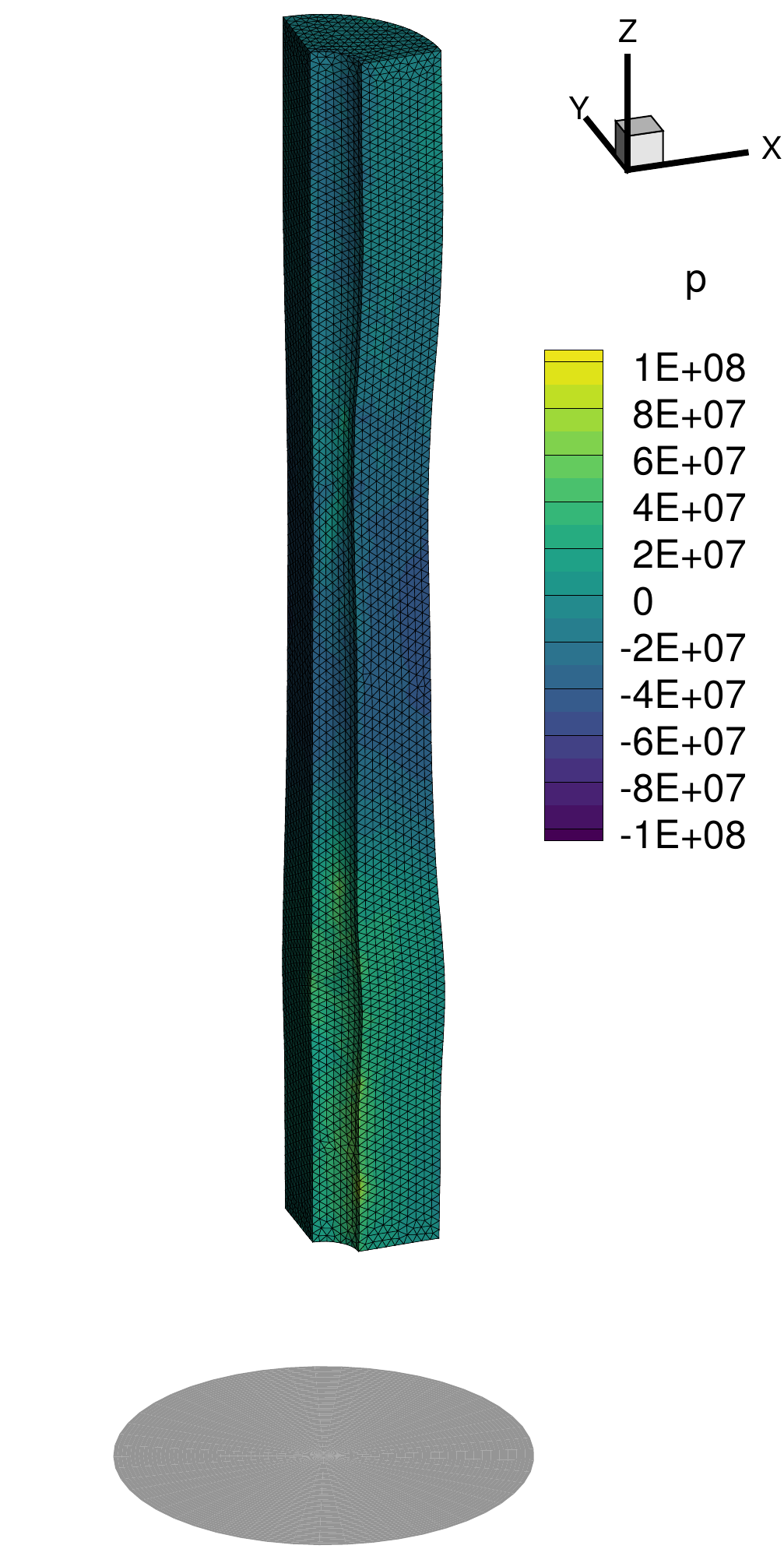} 
    \end{tabular} 
    \caption{Rebound of a hollow circular bar ---
      Time evolution of the deformation and pressure distribution at output times at times $t=50~\mu\text{s}$ then $75$, $100$, $125$, $150$, $200$, $300$ and the final time $t=325~\mu\text{s}$ (from top left to bottom right). }
    \label{fig.BarRebound3D}
  \end{center}
\end{figure}

\begin{figure}[!htbp]
  \begin{center}
    \begin{tabular}{cc} 
      \includegraphics[width=0.47\textwidth]{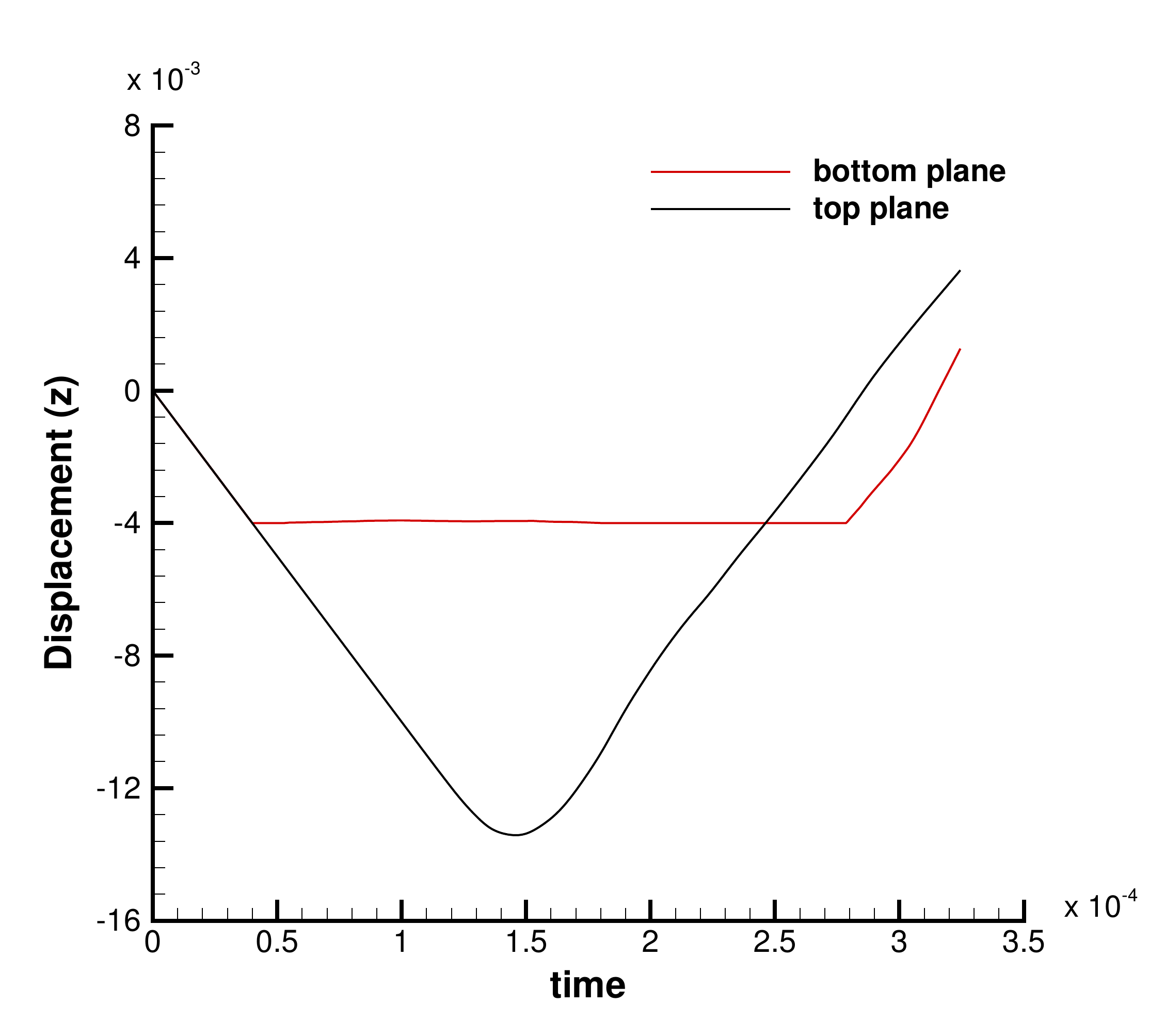} &       
      \includegraphics[width=0.47\textwidth]{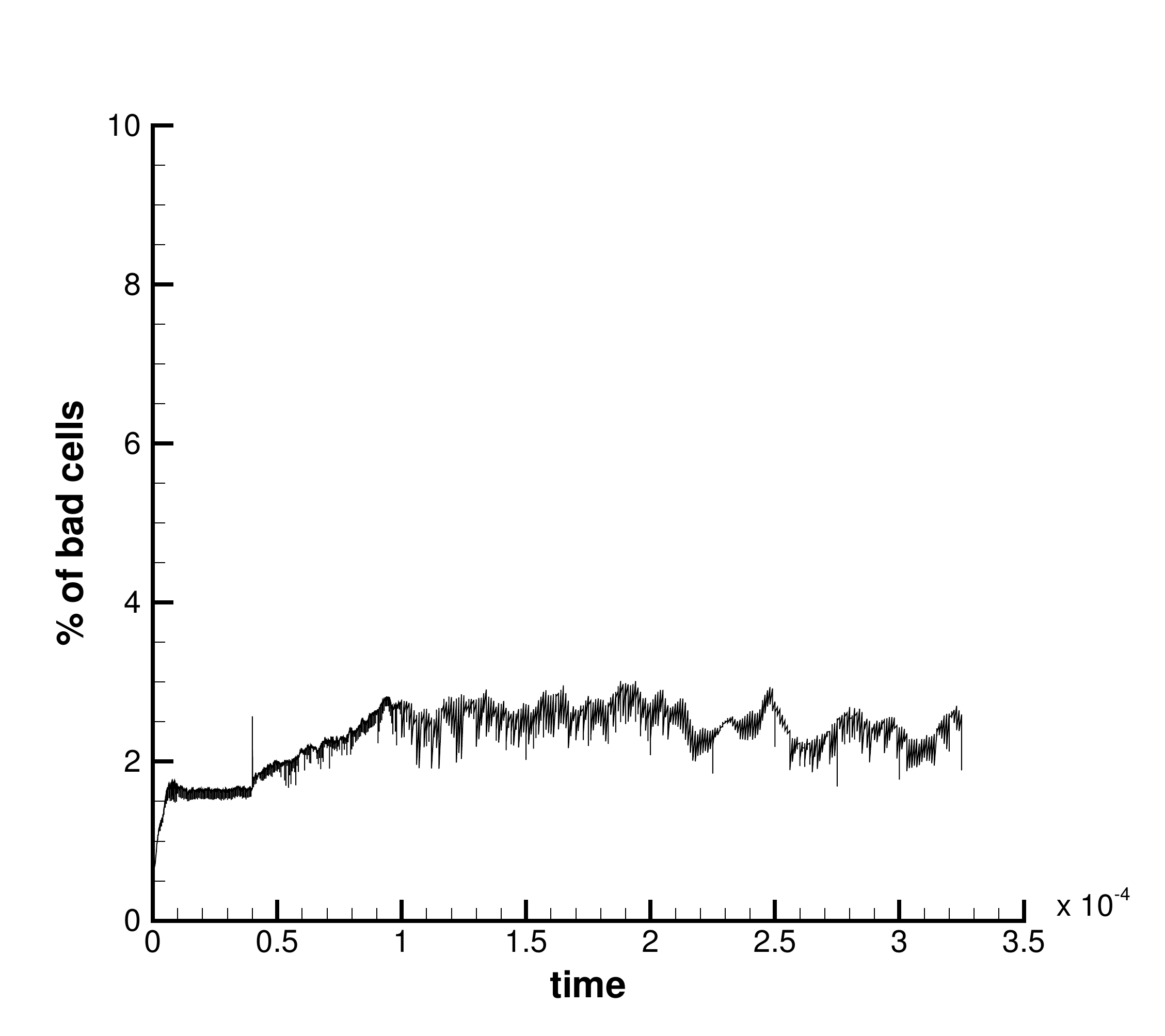}      
    \end{tabular} 
    \caption{Rebound of a hollow circular bar ---
      Time evolution of vertical displacement of the points on the top plane $\mathbf{x}_T=(1.6,0,32.4)\cdot 10^{-3}\text{m}$ and on the bottom plane $\mathbf{x}_B=(1.6,0,4)\cdot 10^{-3}\text{m}$ (left) and percentage of bad cells detected at each time step (right).}
    \label{fig.BarRebound3D_planes}
  \end{center}
\end{figure}

%\walter{setting from \cite{Haider_2018}}	
% Problem description
Taken from \cite{Haider_2018} as the 3D extension of a 2D contact problem found in \cite{Donea2003}, the impacting bar
test consists in the rebound of a hollow circular bar of outer diameter $6.4$~mm, inner diameter
$2$~mm and height $H= 32.4$~mm, see figure~\ref{fig:sketch_column}.
The bar impacts against a rigid friction-less wall with an initial velocity of $\vec{v}_0= (0, 0, -100)^t$~m/s
and the separation distance between the bar and wall is $4$~mm.
Before the impact time at $t=40$ $\mu$s the bar is on a ballistic flight.
Upon impact, the bar undergoes large compressive deformation until $t = 150$ $\mu$s when all the kinetic energy of the bar is converted into internal strain energy.
Afterwards, tensile forces develop and a bounce-off motion initiates in such a way that, 
at approximately $t \simeq 250$ $\mu$s, the bar completely detaches from the wall and moves upwards, still enduring internal milder deformations.
% Model
The neo-Hookean constitutive model is chosen with density $\rho^0 = 8930$~kg/m$^3$, Young's modulus $E = 585$~MPa and Poisson's ratio $\nu = 0.45$ and the final time is set to $326$ $\mu$s. \\
% Boundary conditions
The fixed wall is the $x-y$ plane and is considered as a restricted tangential displacement type BCs.
The rest of the material is subject to free-traction BCs. Special care must be paid to the points of the inner circle at the bottom of the bar.
Indeed for these points the BCs must evolve from free-traction to slip-wall BCs during the contact time up to
detachment. Specifically, free-traction BCs are used until the velocity of the nodes lying on the bottom face is downward pointing and the distance to the wall is greater than a prescribed tolerance of $10^{-12}$. As soon as the new node position would exceed the $z-$coordinate of the wall, i.e. $z=0$, the time step is modified in order to let the bar exactly hit the wall, then the boundary condition switches to slip wall type from the next time step on. Then, when the velocity of the bottom face nodes becomes upward pointing because of the rebound of the bar, as soon as the new node position would detach from the wall, the time step is again modified so that it exactly matches the time of detachment and finally the boundary condition changes again to free-traction for the rest of the simulation. 
% Numerical set up
One quarter of the hollow bar is meshed with $N_c=12254$ tetrahedra and a characteristics length of $1/128$.
In figure~\ref{fig.BarRebound3D} we present the time evolution of the deformation and pressure distribution (colors) at times $t=50~\mu\text{s}$ then $75$, $100$, $125$, $150$, $200$, $300$ and the final time $t=325~\mu\text{s}$.
The main behaviors and deformations are captured by the numerical simulations as compared to the results in \cite{Haider_2018}.
% ---- FIG ---------
% ---- FIG ---------
Following \cite{Haider_2018} (see Fig. 27), we present on the left panel of figure~\ref{fig.BarRebound3D_planes} the time evolution of vertical displacement of the points on the top $\mathbf{x}_T=(1.6,0,32.4)\cdot 10^{-3}\text{m}$ (black) and bottom $\mathbf{x}_B=(1.6,0,4)\cdot 10^{-3}\text{m}$ (red) planes. The general behavior is again qualitatively reproduced.
At last, on the right panel of figure~\ref{fig.BarRebound3D_planes}, we show the percentage of bad cells detected by the \aposteriori limiter and observe that, on average, less than $3\%$ demands limiting at each iteration.
This induces a rather efficient limiting procedure compared to classical \apriori slope limiters.
% ---- FIG ---------
% ---- FIG ---------

%
% TEST # 7 :Impact of a jelly-like droplet
%
\subsection{Impact of a jelly-like droplet} \label{ssec.JellyDrop}
\begin{figure}[!htbp]
  \begin{center}
    \begin{tabular}{cc} 
      \includegraphics[width=0.47\textwidth]{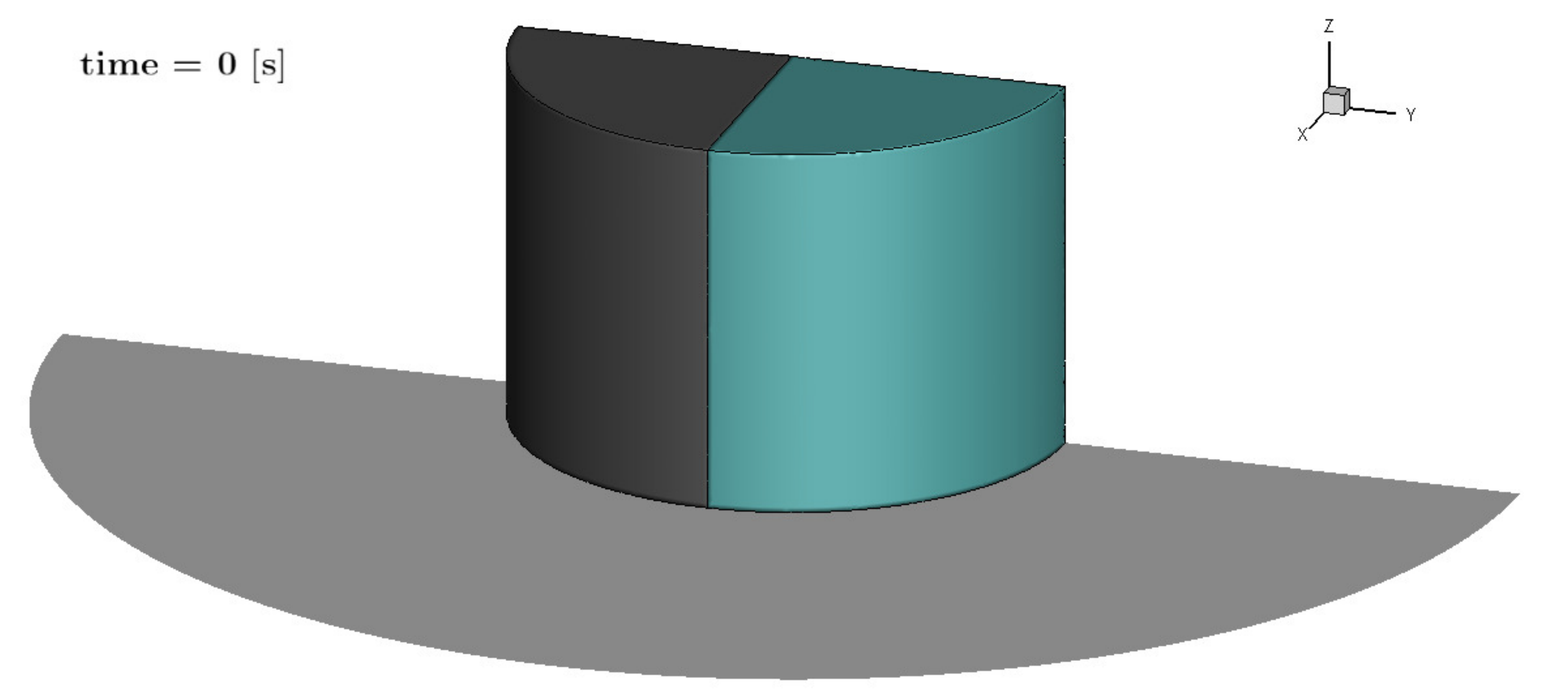} &       
      \includegraphics[width=0.47\textwidth]{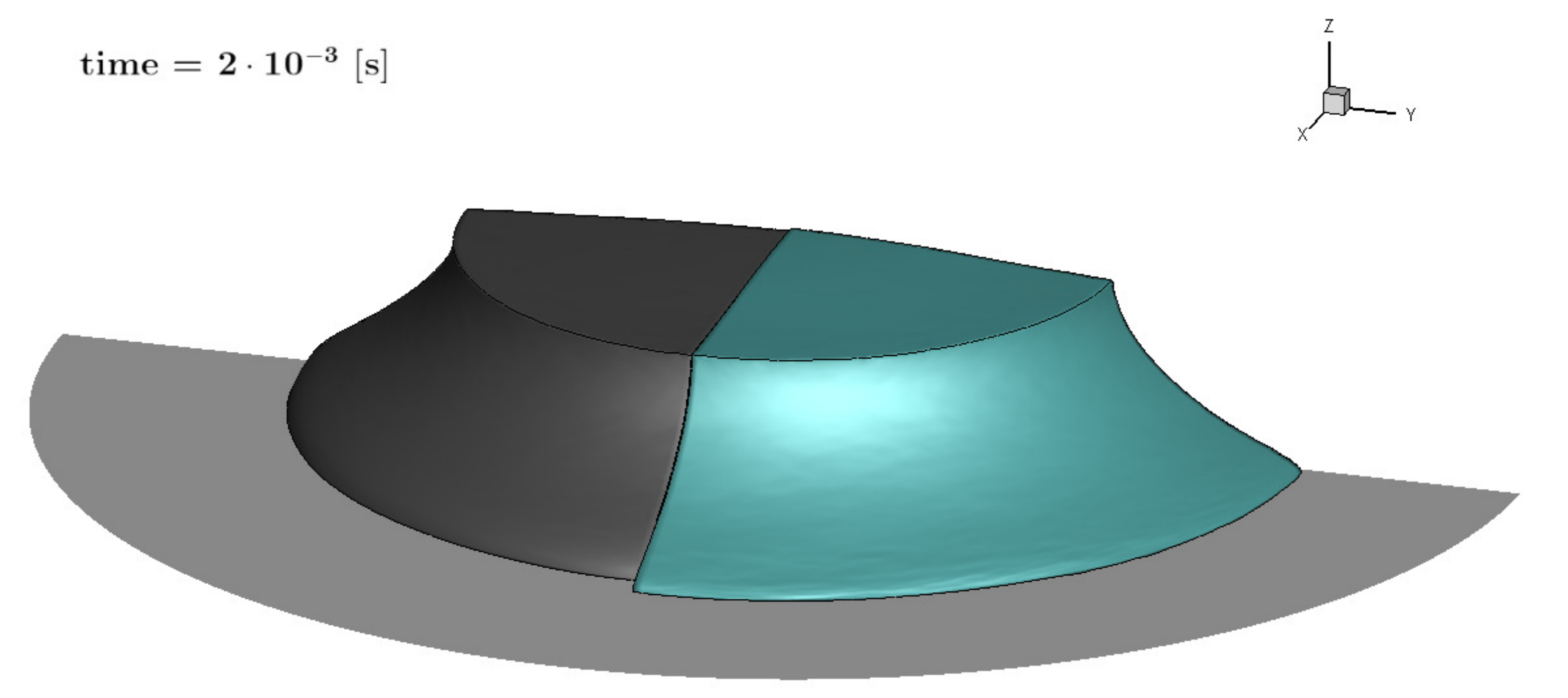} \\     
      \includegraphics[width=0.47\textwidth]{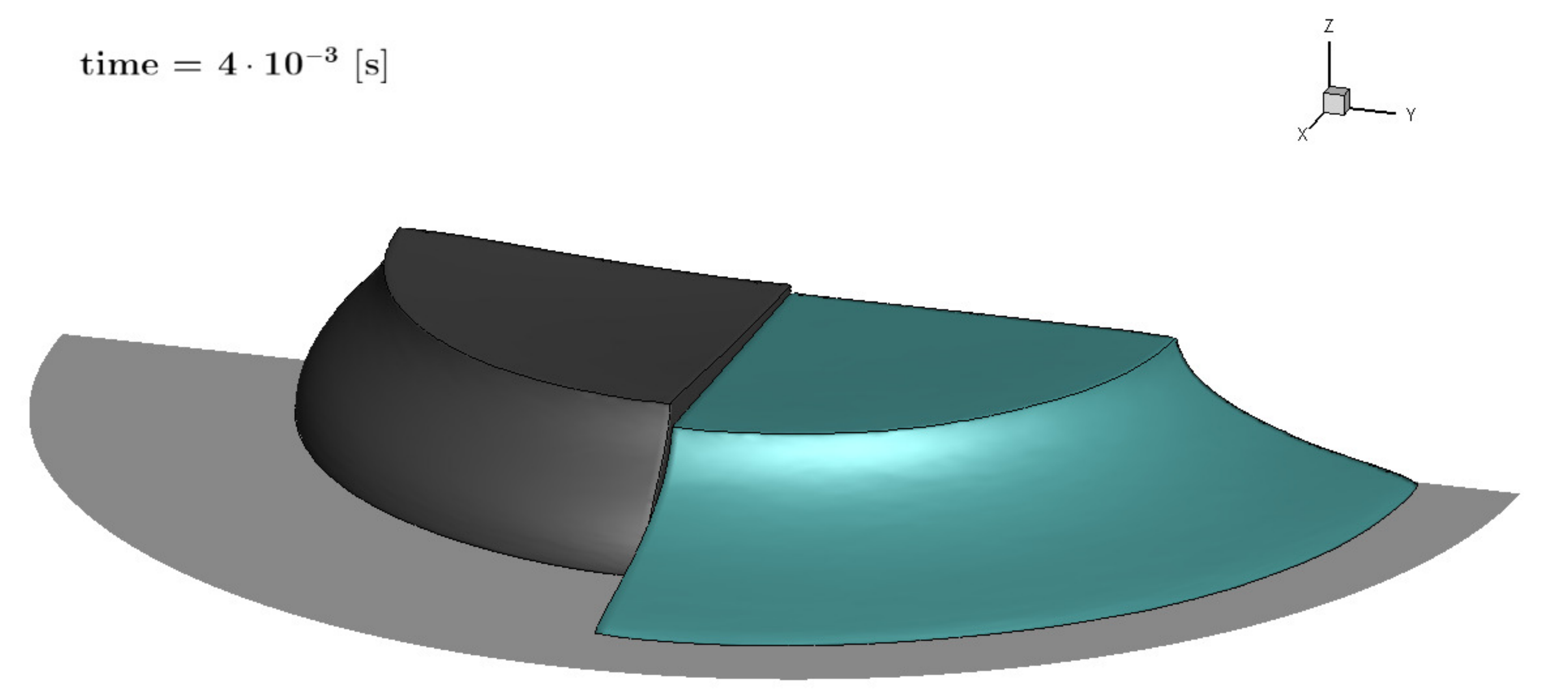} &       
      \includegraphics[width=0.47\textwidth]{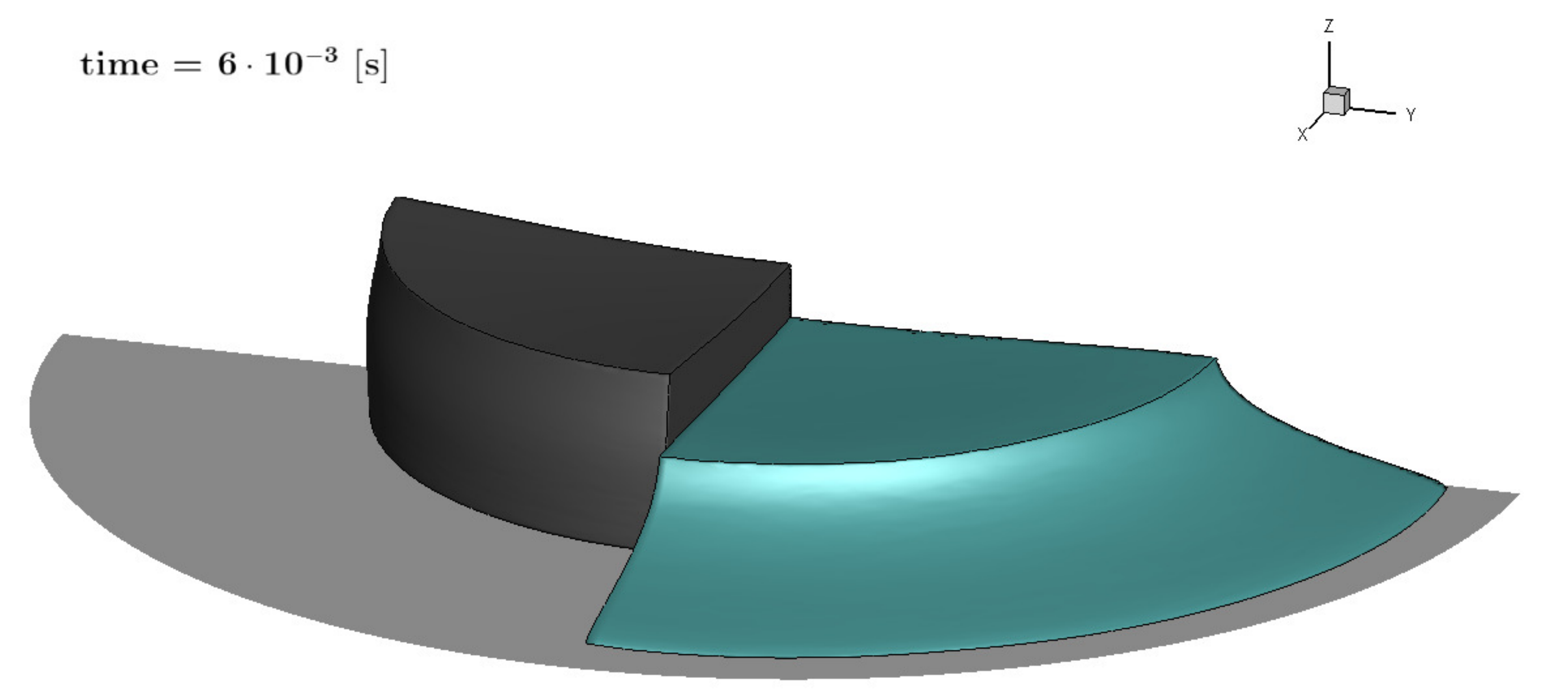} \\
      \includegraphics[width=0.47\textwidth]{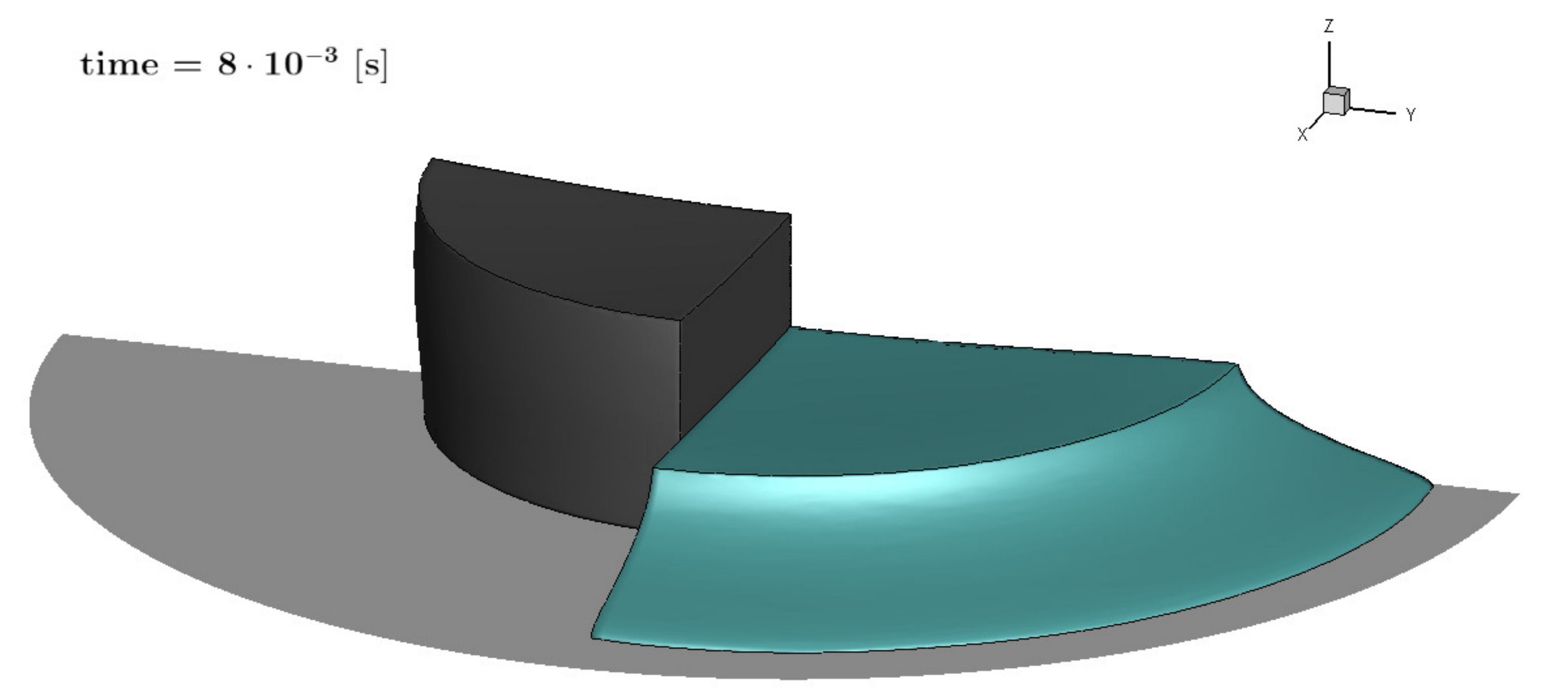} &       
      \includegraphics[width=0.47\textwidth]{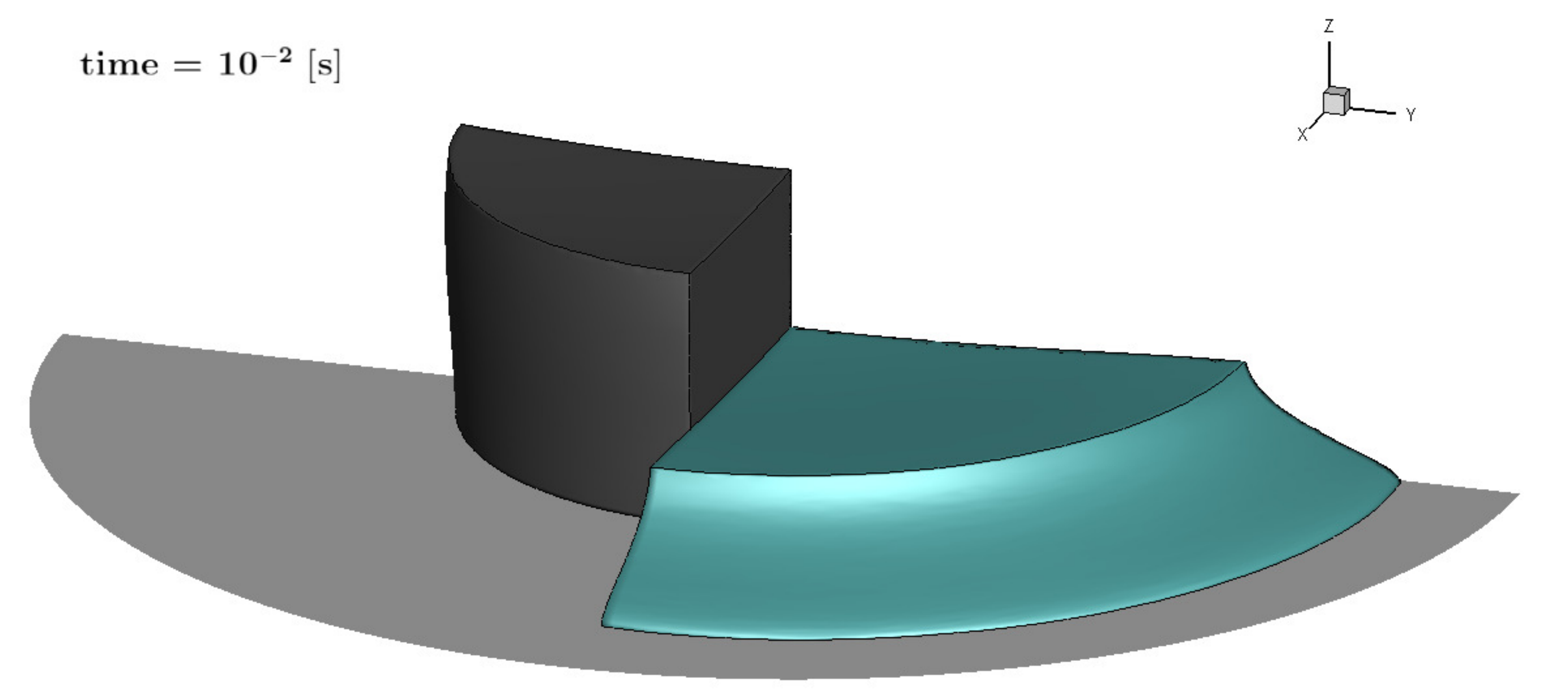} \\
    \end{tabular} 
    \caption{Impact of a jelly droplet with impact velocity $3~\text{m.s}^{-1}$ ---
      Time evolution of the droplet shape at different output times for neo-Hookean model ($a=-1$, petroleum shade) or non-linear one ($a=0$, black shade).}
    \label{fig.Jelly3D}
  \end{center}
\end{figure}

\begin{figure}[!htbp]
  \begin{center}
    \begin{tabular}{cc} 
      \includegraphics[width=0.47\textwidth]{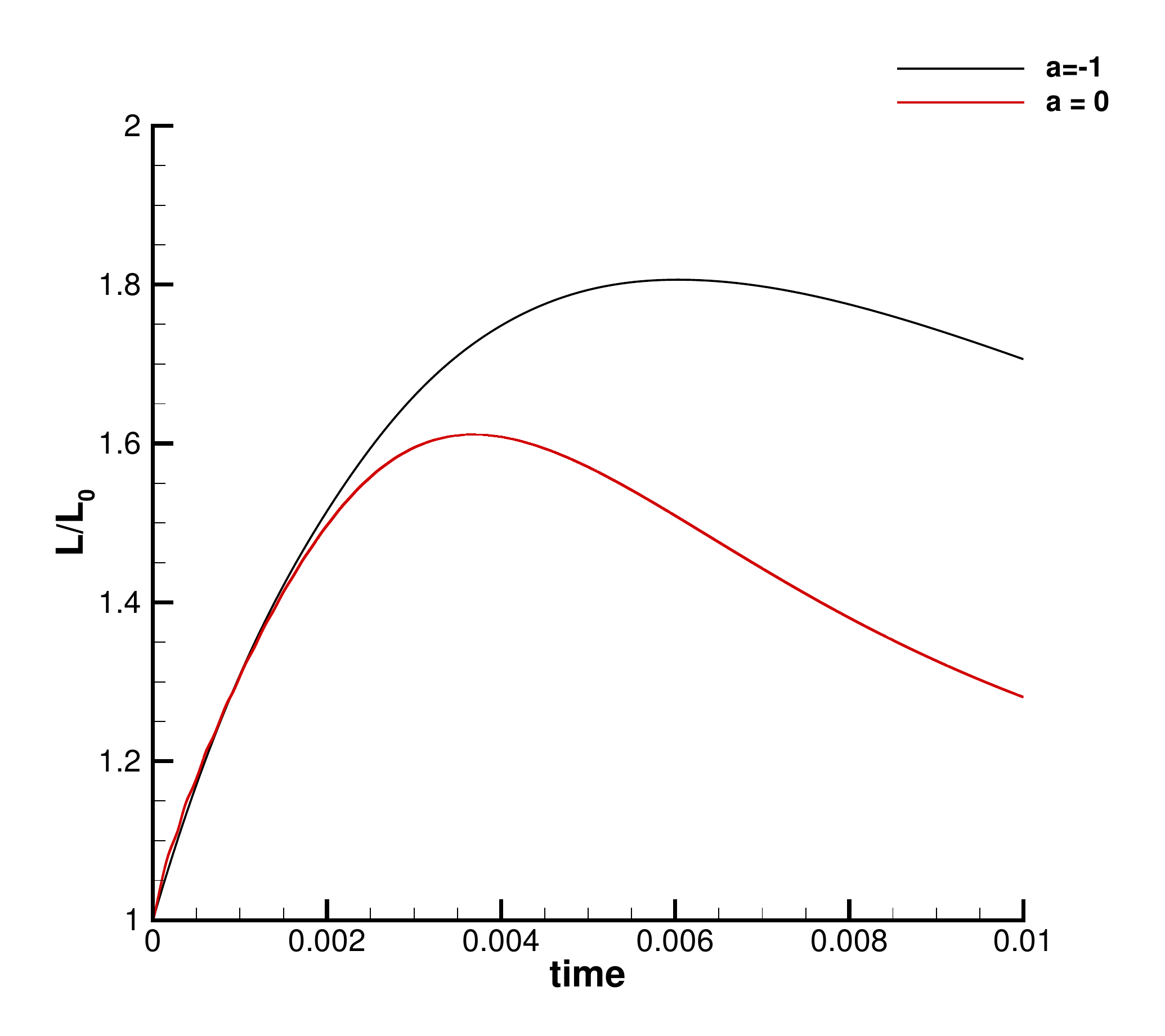} &       
      \includegraphics[width=0.47\textwidth]{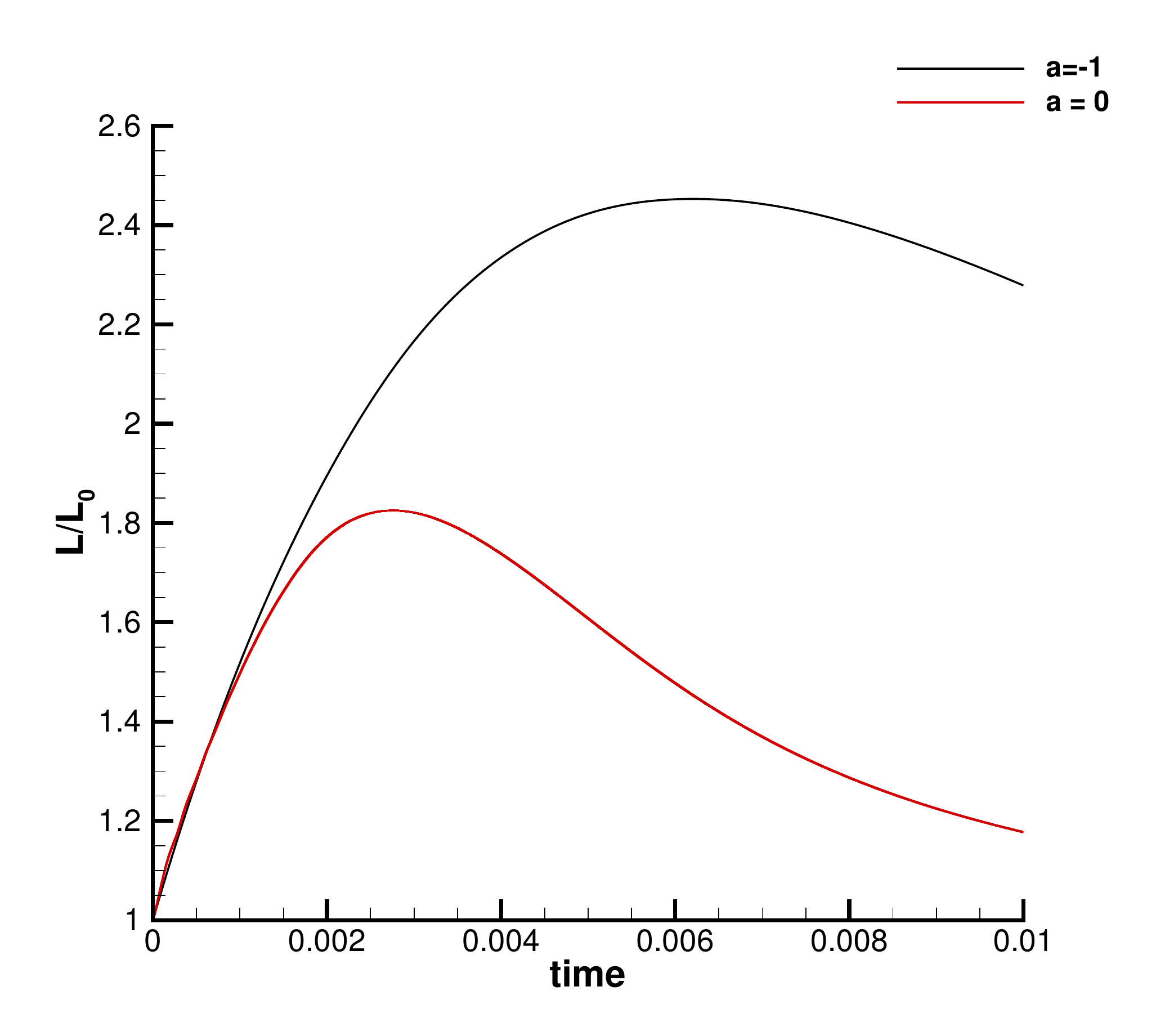} \\    
    \end{tabular} 
    \caption{Impact of a jelly droplet ---
      Time evolution of the maximum spreading of the droplet $L/L_0$ in the case neo-Hookean model ($a=-1$, black line) or non-linear one ($a=0$, red line) ---
      The impact velocity is $2~\text{m.s}^{-1}$ (left) and $3~\text{m.s}^{-1}$ (right).}
    \label{fig.Jelly_radius}
  \end{center}
\end{figure}

% Set-up
As a last test case we consider the impact of a jelly-like material onto a flat rigid horizontal surface,
inspired by the test in \cite{Hank2017}.
An initially cylinder of clay (bentonite) of diameter $L_0$ and height $h$ moves downward with velocity
$\vec{v}=(0,-v)~\text{m.s}^{-1}$, and material parameters $\gamma=2.2$, $p_\infty=10^6$, $\mu=85~\text{Pa}$, $\rho_0=1020$~kg/m$^3$,
Experiments of such impacts have been carried on in particular in \cite{luu_forterre_2009} on different types of surface
(smooth glass, hydrophobic).
In such situation we are interested in the final diameter of the impacting droplet $L$
and the experimental results show a quasi-linear behavior of the maximal spread factor with respect to the impact velocity.
Initially $L_0=12~\text{mm}$ and $h=8~\text{mm}$, and
two impact velocities are considered, $v=2$ and $3~\text{m.s}^{-1}$.
% Numerical setup
% - mesh
The numerical simulation considers a 3D polyhedral computational domain constituted by an approximation of $1/4$ of the initial bentonite cylinder $\Omega^0$ by a mesh made of $N_c=717396$ tetrahedra with characteristics length $1/100$.
% - material behaviors
Two constitutive laws are tested, namely the neo-Hookean model, $a=-1$, and the non-linear one $a=0$, see section~\ref{ssec:neo-hookean} for details.
% - BCs
Symmetry BCs are imposed for the $x=0$ and $y=0$ planes, while free-traction BCs are applied on the top and cylinder boundaries and slip wall type is prescribed on the bottom side.
% Results
In figure~\ref{fig.Jelly3D} are displayed the shapes of the material for successive times $t=2k\times 10^{-3}~\text{s}$
for $0\leq k \leq 5$ in the case of a $v=3~\text{m.s}^{-1}$ impact velocity.
The black shape corresponds to the non-linear model $a=0$, while the petroleum shape corresponds to a neo-Hookean one $a=-1$. They are put in respect to each other for comparison purposes. \\
Regardless of the constitutive model, i.e the value of $a$, the jelly-like material is compressed after the impact and deforms back and forth due to its elastic behavior.
As expected with the neo-Hookean model (petroleum shape) the spread of the droplet is much more pronounced and the droplet retrieves a cylinder-like shape slower compared to the non-linear model (black shape).
% ---- FIG ---------
% ---- FIG --------
In order to quantify this behavior we present in figure~\ref{fig.Jelly_radius} the
maximum spreading of the droplet, $L/L_0$, in the case $a=-1$ (black line) and $a=0$ (red line) for the two impact velocities. The neo-Hookean model produces faster and more pronounced elastic behaviors compared to the non-linear model which retrieves a ratio closer to one faster.
The experimental results in \cite{luu_forterre_2009} provide approximate values $2.25$ and $2.75$, respectively, while our simulations produce $1.8$ and $2.5$ in accordance to the numerical results in \cite{Hank2017}.
% ---- FIG ---------
% ---- FIG ---------	